\documentclass[pdflatex,sn-mathphys-num]{sn-jnl}
\usepackage{fullpage}
\usepackage{geometry}

\usepackage{amsmath,amssymb,amsfonts, amsthm}
\usepackage[]{algpseudocode}
\usepackage{graphicx} 
\usepackage{url}
\newtheorem{theorem}{Theorem}[section]

\newtheorem{assumption}[theorem]{Assumption}

\newtheorem{remark}[theorem]{Remark}
\newtheorem{case}{Case}[section]

\newcommand{\set}[1]{\left\{ #1 \right\}}
\newcommand{\lp}{\left(}
\newcommand{\rp}{\right)}

\newcommand{\Def}{\overset{\textbf{def}}{=}}
\newcommand{\RR}{\mathbb{R}}

\newcommand{\NN}{\mathbb{N}}

\newcommand{\eps}{\varepsilon}

\newcommand{\norm}[1]{\left\| #1 \right\|}

\newcommand{\+}[1]{\ensuremath{\mathbf{#1}}}

\newcommand{\sT}{\mathsf T}

\newcommand{\bigmid}{\bigl\vert}

\newcommand{\tup}[1]{\textup{(}#1\textup{)}}
\usepackage{amsopn}

\DeclareMathOperator{\ball}{Ball}

\newcommand{\sa}{\mathsf a}

\renewcommand{\vec}[1]{\mathbf{#1}}

\newcommand{\elements}[1]{\left[#1 \right]^\sT}

\usepackage{algorithm}%
\usepackage{textcomp}
\usepackage{xcolor}
\usepackage{graphics}
\usepackage{caption}

\usepackage{subcaption}
\usepackage{arydshln}
\usepackage{mathtools}
\usepackage{relsize}
\usepackage{BOONDOX-uprscr}
\usepackage{dirtytalk}
\usepackage{enumerate}  

\usepackage{color}
\definecolor{darkgreen}{rgb}{0,0.5,0}
\usepackage[colorinlistoftodos]{todonotes}

\usepackage{soul}
\setulcolor{red} 
\setstcolor{red} 
\sethlcolor{yellow} 

\newcommand{\edit}[1]{\textcolor{black}{#1}}

\usepackage{graphicx}%
\usepackage{multirow}%
\usepackage{amsmath,amssymb,amsfonts}%
\usepackage{amsthm}%
\usepackage{mathrsfs}%
\usepackage[title]{appendix}%
\usepackage{xcolor}%
\usepackage{textcomp}%
\usepackage{manyfoot}%
\usepackage{booktabs}%
\usepackage{algorithm}%
\usepackage{algorithmicx}%
\usepackage{algpseudocode}%
\usepackage{listings}%

\begin{document}

\title[Article Title]{Second-Order Time to Collision With Non-Static Acceleration}

\author[1]{\fnm{Hossein} \sur{Nick Zinat Matin}}\email{h-matin@berkeley.edu}
\equalcont{These authors contributed equally to this work.}
\author[1]{\fnm{Yuneil} \sur{Yeo}}\email{yuneily@berkeley.edu}
\equalcont{These authors contributed equally to this work.}

\author[1]{\fnm{Amelie} \sur{Ju-Kang Ngo}}\email{amelie\_ngo@berkeley.edu}

\author[2]{\fnm{Antonio~R.} \sur{Paiva}}\email{antonio.paiva@allstate.com}

\author[2]{\fnm{Jean} \sur{Utke}}\email{jutke@allstate.com}

\author*[1]{\fnm{Maria Laura} \sur{Delle Monache}}\email{mldellemonache@berkeley.edu}

\affil*[1]{\orgdiv{Civil and Environmental Engineering}, \orgname{University of California at Berkeley}, \orgaddress{\street{Sutardja Dai Hall}, \city{Berkeley}, \postcode{94720}, \state{California}, \country{USA}}}

\affil[2]{\orgdiv{D3 Analytics Innovation}, \orgname{Allstate Insurance Company}, \orgaddress{\street{3100 Sanders Road}, \city{Northbrook}, \postcode{60062}, \state{Illinois}, \country{USA}}}

\abstract{We propose a second-order time to collision (TTC) considering non-static acceleration and turning with realistic assumptions\edit{. This is equivalent to considering that} the steering wheel \edit{is held at} a fixed \edit{angle} with constant pressure on the gas or brake pedal \edit{and matches the well-known bicycle model}. Past works that use acceleration to compute TTC consider \edit{only longitudinally aligned} acceleration.

We additionally develop and present \edit{the Second-Order Time-to-Collision Algorithm using Region-based search (\edit{STAR})} to efficiently compute the proposed second-order TTC and overcome the current limitations of the existing built-in functions. The evaluation of the algorithm in terms of error and computation time is conducted through statistical analysis. 

Through numerical simulations and publicly accessible real-world trajectory datasets, we show that the proposed second-order TTC with non-static acceleration is superior at reflecting accurate collision times, especially when turning is involved.}



\maketitle

\section{Introduction}\label{sec:introduction}

Given the high use of vehicles in the transportation system, the public emphasizes the importance of ensuring safety in \edit{mixed-traffic} environments. In response to this vital need, many works evaluate the safety of mixed traffic using different methods. Since collisions and crashes are rare, the safety of mixed traffic cannot be \edit{assessed} fully based on historical collision records \cite{nikolaou2023exploiting}. As an alternative, many analysis utilize surrogate or indirect safety measures \edit{to} evaluate the collision risk, \edit{potential for }near-crashes, and the \edit{overall} safety \edit{in} mixed traffic without \edit{requiring} accident data \cite{wang2021review, johnsson2018search}. 

One of the most popular safety surrogate measures is time to collision (TTC)~\cite{arun2021systematic}. According to \cite{hayward1971near} and \cite{hou2015new}, TTC, or first-order TTC, is the time at which two vehicles would collide if vehicles were to move continuously \edit{in} their current direction and at current velocity. 

Due to the simplicity of calculating TTC based on measurements from sources like cameras, radar, LIDAR, and GPS, TTC has been \edit{studied in the context of} Advanced Driver Assistance Systems (ADAS), such as Forward Collision Warning Systems (FCWS) and Emergency Braking Systems (EBS) \cite{tak2018comparison}. These systems detect and estimate the potential collision risk in the future and take appropriate actions like triggering emergency brak\edit{ing} or \edit{warning} drivers to prevent and avoid collision \cite{tak2018comparison}. TTC \edit{has also} been \edit{studied for use} in \edit{Autonomous Vehicles~(AVs)} to ensure safety and reliability \cite{goudarzi2024collision} and to assess the safety impacts of AVs in traffic \cite{zhang2021evaluating}. 

Numerous studies on time to collision have explored extensions of the concept, \edit{aiming to account for different conditions or emphasize certain criteria.} \textit{Worst TTC} evaluates the worst-case collision time \edit{(i.e., smallest TTC)} among different trajectory pairs, offering a comprehensive risk assessment \cite{wachenfeld2016worst}. \textit{Time Exposed TTC} refers to the duration of time when the evaluated TTC is below a particular threshold, while \textit{Time Integrated TTC} is the total cumulative sum of all evaluated TTC \edit{values} below this threshold \cite{chin1991traffic, chin1997measurement}. Different works extend the concept of TTC for AVs and Vehicle-to-Vehicle~(V2V) Infrastructure systems \cite{lakhal2019risk}.

The definition of the first-order TTC makes it simple to calculate, but it has its limitations. The first-order TTC model assumes constant velocity, or no acceleration, which limits its applicability to straight trajectories \cite{venthuruthiyil2022anticipated}. This assumption fails when vehicles perform maneuvers such as turns at intersections or navigating roundabouts. 

Several works consider acceleration to overcome the limitations of the first-order TTC. \textit{Potential TTC} extends TTC by considering the deceleration of the leading vehicle, whereas the following vehicle has constant velocity \cite{wakabayashi2003traffic}. While the model offers a more nuanced collision assessment with the acknowledgment of the variability in motion, the model can be utilized only for the specific scenario. In \edit{\citet{brown2005adjusted}}, the author presents the concept of \textit{Adjusted Minimum TTC} that considers velocity, position, and acceleration in a one-dimensional case. \textit{Enhanced TTC} introduced in \edit{\citet{kiefer2005developing}} is another notable example of works that presents a new definition for TTC considering the case where a leading vehicle is decelerating. However, the concept is only applicable to the straight-line trajectory in a one-dimensional setting. 

The second-order TTC model, introduced by \edit{\citet{ward2015extending}}, also incorporates vehicle acceleration to address this limitation. However, \edit{that} work only consider\edit{s} a \say{static} acceleration vector, meaning its magnitude and direction \edit{are assumed to} remain constant over time \edit{when predicting future positions and determining whether they might collide}. This static acceleration assumption is unrealistic when vehicles are not traveling in a straight line, as the acceleration direction relative to velocity changes over time during turns. In addition, the second-order TTC \edit{computation} in \edit{\citet{ward2015extending}} requires the process of finding the loom points, the points on the ego vehicle closest to an agent vehicle \edit{of} interest. The process of finding the loom points would be complex as the points closest to an agent vehicle keep chang\edit{ing} in the estimated trajectories over time. 

\textbf{Focus and Contribution.} 
The objective of this paper, therefore, is to present the concept of the second-order TTC with more realistic assumptions, such as non-static acceleration and turning. 

Due to the complexity \edit{of} obtaining the exact formula of the second-order TTC, we introduce \edit{the Second-Order Time-to-Collision Algorithm using Region-based search (STAR)}, which computes the second-order TTC efficiently. Through statistical tests, we conclude that \edit{STAR} can produce accurate TTC with faster computation time compared to \edit{a} step-by-step simulation. 
To gain insights into the dynamics of TTC under different conditions, we employ numerical simulations to systematically investigate a range of scenarios using first-order and second-order TTC. These simulations are carefully crafted to reflect realistic conditions and parameter ranges encountered in intersections. 

In addition to numerical simulations, we incorporate real-world datasets to enhance the applicability of our study, ensuring that our findings are grounded in empirical observations and reflective of actual driving scenarios. The simulations prove that the second-order TTC is more effective in predicting time to collision than the first-order TTC.

The paper is organized as follows. Section~\ref{sec: concepts} presents the idea of the second-order TTC with non-static acceleration. Section~\ref{sec: algorithm} presents the algorithm for computation of the second-order TTC. Section~\ref{sec: evaluation} evaluates the proposed algorithm with its performances in terms of error and computation time. Section~\ref{sec: analysis} compares the first-order TTC with different simulated scenarios. Section~\ref{sec: real_data} additionally compares the performance of the first-order TTC and the second-order TTC through real-world data. 

\section{Second-Order TTC} \label{sec: concepts}
In this section, we present the concept of the second-order TTC with non-static acceleration \edit{and} provide a brief overview of the first-order TTC. 
\edit{\subsection{Notation}
We start by introducing some of the notations and definitions that will be used frequently in this work. 
\begin{itemize}
    \item We use boldface notations to represent the vectors and the elements of the vector are denoted by $\vec x = \elements{x_1, x_2} \in \RR^2$. 
    \item The inner product of two vectors $\vec x$ and $\vec y$ is denoted by $\vec x \centerdot \vec y$. 
    \item We define 
    \begin{equation*}
        \norm{\vec x}_2 \Def (x_1^2 + x_2^2)^{\tfrac 12} = (\vec x \centerdot \vec x)^{\tfrac 12}
    \end{equation*}
    \item Set $\mathbf Q$ denotes the set of intersection points between the trajectories of the vehicles. We also denote $\mathbf q \in \mathbf Q$ as each collision point. 
    \item Throughout the paper $\vec p(t)$, $\vec v(t)$ and $\sa(t)$ denote position, velocity and acceleration vectors at time $t$ respectively. 
    \item $C(\mathcal U; \RR)$ denotes space continuous functions from the domain $\mathcal U$ which take values in $\RR$ and $C^1(\mathcal U; \RR)$ denote the space of continuously differentiable functions from $\mathcal U$ with values in $\RR$.
    \item We use $\ball(\vec x_\circ, r) = \set{\vec x:  \norm{\vec x - \vec x_\circ}_2 < r}$ to denote an open ball centered at a point $\vec x_\circ$ and radius $r$. 
\end{itemize}}

\subsection{Dynamical model}
\edit{We start by introducing some assumptions on the underlying dynamics of vehicles. There have been various dynamical models that introduce and analyze the safety and stability of the second-order microscopic behavior both the Human-Driven Vehicles (HV) as well as (Connected) and Automated Vehicles (CAVs); we refer the readers to the incomplete list \cite{bando1995dynamical, wilson2011car, pipes1953operational, newell1961nonlinear, nick2022near, chou2024stability, 10384086, 9147363, 9147244, matin2024analytical, delle2019feedback}.
In addition, several models present multilane and lane-changing dynamics; see, for example, \cite{wu2017multi, kesting2007general, khelfa2023predicting}.}

\edit{Let $t \mapsto \vec \xi(t)$ and $t \mapsto \vec \zeta(t)$ denote the position and velocity vectors respectively. In particular, we consider
\begin{equation}\label{E:car_dynamic}
         (\dot{\vec \xi}, 
        \dot{\vec \zeta})^\sT = F(t,\vec \xi, \vec \zeta), \quad (\vec \xi(0), \vec \zeta(0)) = (\vec \xi_\circ , \vec \zeta_\circ)
\end{equation} 
where $B^\sT$ denotes a transpose of a vector $B$. The function $F: \RR_+ \times \RR^2 \times \RR^2 \to \RR^4$ is considered to be sufficiently smooth such that the solution of \eqref{E:car_dynamic} exists uniquely until the collision time $\mathcal T$ (see equations~\eqref{E:collision_first} and \eqref{E:collision_second}) and 
\begin{equation}\label{E:properties}
   \vec \zeta \in C^1([0, \mathcal T); (\ball(0, v_{max}))^2), \quad \vec{\dot \zeta} \in C([0, \mathcal T); (\ball(0, \sa_{\max}))^2),
\end{equation}
which implies physically feasible dynamics with bounded velocity and acceleration. Starting from any time $t = t_\circ$, we will approximate TTC by approximating the trajectory of the vehicles based on their possible behavior denoted by positions $\+{p}$, velocities $\+{v}$, and accelerations $\+{a}$ (See \eqref{E:1dttc_traj}, \eqref{E:2dttc_traj_linear}, \eqref{E:2dttc_traj_circle}, and \eqref{E:2dttc_traj_circle2}) and investigate the TTC from $t \ge t_\circ$. Details regarding the approximated trajectories of vehicles in both schemes are covered later.} 

\subsection{First-Order TTC}
First-order TTC is the simplest TTC with the assumption that vehicles maintain a constant velocity. Consequently, it is assumed that the vehicles continue moving along straight trajectories in the same velocity direction.

For first-order TTC, a positional vector at the center of a vehicle represents the position of a vehicle \+{p}. A vehicle is also assumed to be moving with a velocity \+{v}\edit{$_\circ$}. Thus, the motion of a vehicle can be described as:
\begin{equation}\label{E:1dttc_traj}
\begin{cases}
    \vec p(t)  \edit{= \+{p}(t; t_\circ)} = \+{p}_\circ + \+{v}\edit{_\circ} \edit{(t-t_\circ)}, \edit{\quad t \in [t_\circ, \mathcal T)} \\
    \edit{(\vec p(t_\circ),\vec v(t_\circ)) = (\vec p_\circ,\vec v_\circ) }\\
    \end{cases}
\end{equation} 
where $\+{p}(t)$ represents the position of a vehicle at time $t$ starting from
\begin{equation}\label{E:intial_condition}
        (\vec p_\circ, \vec v_\circ) = (\xi(t_\circ), \zeta(t_\circ))
\end{equation}
as defined in \eqref{E:car_dynamic}.

\edit{Using \eqref{E:1dttc_traj}}, the distance between vehicle $i$ and vehicle $j$ can be calculated by
\begin{equation}\label{E:d_ij_distance}
\begin{split}
    d_{ij}(t) &\edit{\Def} \lVert \+{p}_i(t) - \+{p}_j(t) \rVert_2, \edit{\quad t \in [t_\circ, \mathcal T)} \\
    \vspace{-1ex}
    &= \lVert (\+{p}_{i,\circ} - \+{p}_{j,\circ}) + (\+{v}_{i,\edit{\circ}} - \+{v}_{j,\edit{\circ}})\edit{(t-t_\circ)}\rVert_2\\   
\end{split}\end{equation}
where $\+{p}_i(t)$ and $\+{p}_j(t)$ denotes the position of vehicle $i$ and vehicle $j$ respectively. 

\begin{figure}[ht!]
    \centering
  \includegraphics[width=0.8\linewidth]{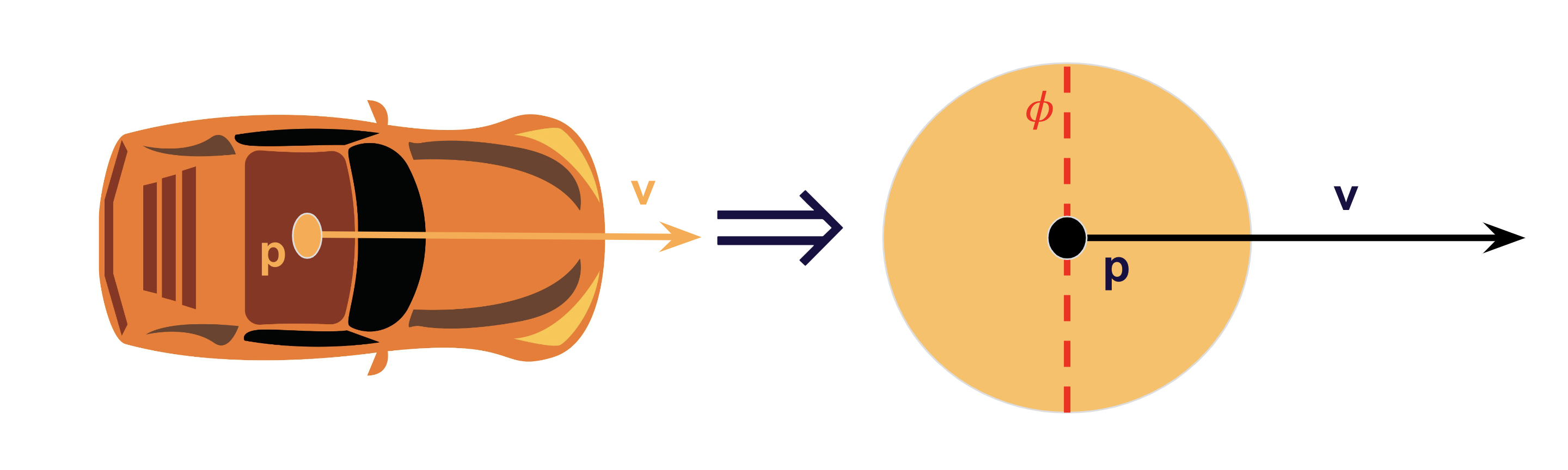}
    \caption{\edit{Representation of position and velocity of a vehicle. Vehicles' dimension is approximated by circles with dimension $\phi$.}}
    \label{fig:transition_to_circle}
\end{figure}

To simplify the treatment of \edit{the} vehicle's dimensions, we approximate their shape as a circle. If we designate the center of the vehicle as the reference position, a collision occurs when the distance between the projected positions is less than the sum of the radius of both circles. In the scenario without acceleration, this can be determined by solving for $d_{ij}^2(t) = \phi^2$, where $\phi$ is the diameter of the circles if they have the same radius. Figure \ref{fig:transition_to_circle} shows \edit{the representation of the position and the velocity of a vehicle and }the approximation of \edit{a} vehicle's shape as a circle. 

\edit{In addition, from \eqref{E:d_ij_distance}} the time $t^*$ when \edit{$d_{ij}(t) =\phi$} is \edit{given by}
\edit{
\begin{equation} \label{eq:1d_ttc}
\begin{cases}
t^* = t_\circ + \frac{-\lp (\+{p}_{i,\circ} - \+{p}_{j,\circ})\centerdot (\+{v}_{i,\edit{\circ}} - \+{v}_{j,\edit{\circ}}) \rp \pm \sqrt z}{\lVert \+{v}_{i,\edit{\circ}} - \+{v}_{j,\edit{\circ}} \rVert_2^2} \\
z \Def  \lp (\+{p}_{i,\circ} - \+{p}_{j,\circ}) \centerdot(\+{v}_{i,\edit{\circ}} - \+{v}_{j,\edit{\circ}})\rp^2 - \lVert \+{v}_{i,\edit{\circ}} - \+{v}_{j,\edit{\circ}} \rVert_2^2 \cdot \lp\norm{\+{p}_{i,\circ} - \+{p}_{j,\circ} }_2^2 - \phi^2 \rp . \end{cases}
\end{equation}}

\edit{Therefore considering \eqref{eq:1d_ttc},} first-order TTC 
\begin{equation}\label{E:collision_first}
\begin{split}
    \edit{\mathcal T = \mathcal T_F = \mathcal T_F(t_\circ)} &\edit{\Def} \edit{\inf}\set{t\edit{\geq t_\circ}: d_{ij}(t) - \phi \edit{\leq} 0}\\ &\edit{\equiv} \edit{\inf}\set{t\edit{\geq t_\circ}: \lVert (\+{p}_{i,\circ} - \+{p}_{j,\circ}) + (\+{v}_{i,\edit{\circ}} - \+{v}_{j,\edit{\circ}})\edit{(t-t_\circ)}\rVert_2 - \phi \edit{\leq} 0},
\end{split}\end{equation}
\edit{(see Remark \ref{R:difference_TTC}) has three possible cases as follows:}
\edit{\begin{itemize}
\item[(i)] If $z \edit{ < 0}$, then there are no real solutions and therefore no collision. 
\item[(ii)] If $z=0$, then the collision happens at one point and the collision time can be calculated by
\begin{equation*} \edit{\mathcal T_F} = t_\circ + \frac{-\lp (\+{p}_{i,\circ} - \+{p}_{j,\circ})\centerdot (\+{v}_{i,\edit{\circ}} - \+{v}_{j,\edit{\circ}}) \rp }{\lVert \+{v}_{i,\edit{\circ}} - \+{v}_{j,\edit{\circ}} \rVert_2^2}.\end{equation*}
\item[(iii)] Otherwise, there would be two real solutions, \edit{for which} the smallest \textit{positive} solution \edit{is considered as} first-order TTC.\end{itemize}}

\subsection{Second-Order TTC}\label{sec:ttc2}

\begin{figure} [ht!]
    \centering
  \includegraphics[width=0.6\linewidth]{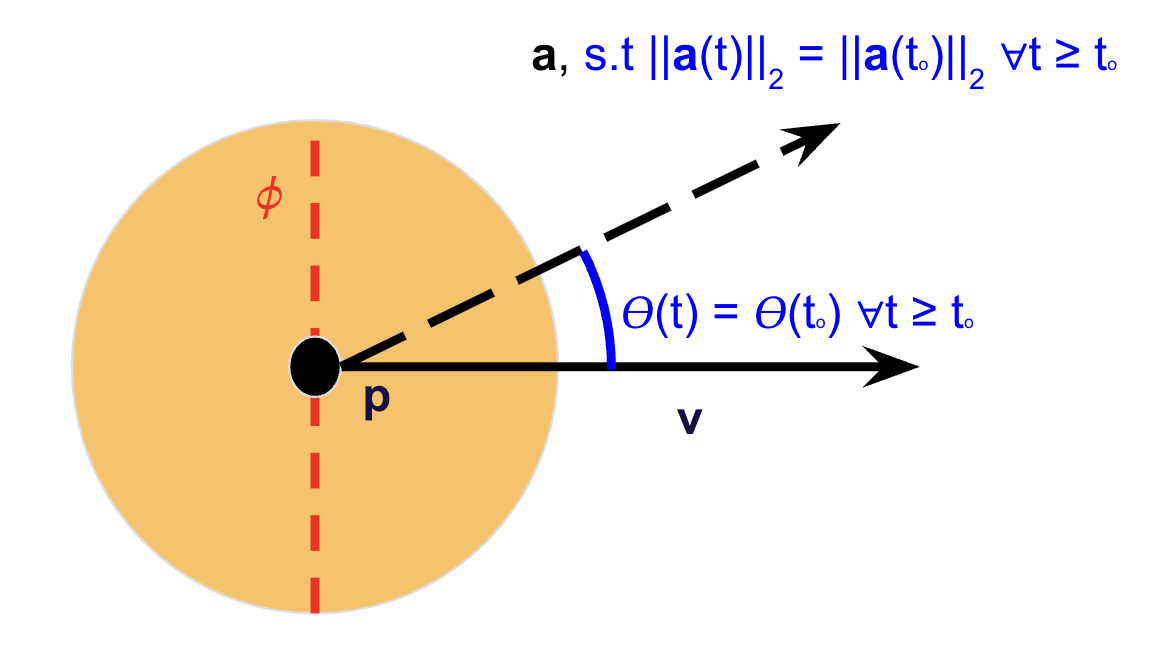}
    \caption{A vehicle with diameter $\phi$ is illustrated. The two assumptions behind the Second-Order TTC computation \edit{are} shown in blue.}
    \label{fig:assumptions}
\end{figure}

Due to the assumption of a constant velocity, the first-order TTC is unrealistic as it does not consider acceleration or deceleration. This also means that vehicle interactions involving non-straight trajectories (e.g., when turning) are not accurately characterized. We consider second-order TTC with non-static acceleration to provide a more realistic TTC. Figure \ref{fig:assumptions} shows two realistic assumptions made to compute the second-order TTC with non-static acceleration. 

The two \edit{standing} assumptions in this case are: 
\begin{enumerate}[\edit{\bf {A}}1.]
  \item \edit{The relative angle $\theta(t)$ between the acceleration $\sa(t)$ and the velocity vectors $\vec{v(t)}$ (see Figure \ref{fig:assumptions}) is constant. This implies that a driver holding the steering wheel in a fixed position.}
  \item \edit{The magnitude of the acceleration vector\edit{, $\norm{\sa(t)}_2$,} is constant, which is associated with holding the same pressure on the gas or brake pedal.}
\end{enumerate}

These assumptions correspond to a constant turn rate and acceleration \edit{magnitude} kinematic model~\cite{lefevre2014survey}, also known as the ``bicycle model." This decouples the changes in velocity and yaw rate simplifying the analysis.

Since the \emph{direction} of the acceleration vector is always changing over time, it is more convenient to use longitudinal ($\edit{\sa_f}$) and lateral ($\edit{\sa_s}$) acceleration components. The longitudinal acceleration, $\edit{\sa_f}$, affects the vehicle's speed while lateral acceleration, $\edit{\sa_s}$, determines how fast the vehicle is turning. 

To compute $\edit{\sa_f}$ and $\edit{\sa_s}$, we set $\edit{(\+{p}_\circ,\+{v}_\circ,\+{a}_\circ) = (\xi(t_\circ), \zeta(t_\circ), \dot \zeta(t_\circ))} $ and \edit{we consider the unit vector} 
\edit{\begin{equation} \label{E:vstar} 
\+v_* = \frac{\edit{\+{v}_\circ}}{\norm{\vec v_\circ}_2} = \elements{v_{*x}, v_{*y}}.
\end{equation}}
\edit{In addition, the orthogonal vector to $\+{v}$ will be denoted by $\+{v}_\perp = \elements{-v_{*y}, v_{*x}}$}. The longitudinal and lateral accelerations then can be determined as $\edit{\sa_f} = \edit{\+{a}_\circ}\centerdot\+{v}_{*}$ and $\edit{\sa_s} = \edit{\+{a}_\circ}\centerdot\+{v}_\perp$. A vehicle is turning left when $\edit{\sa_s} > 0$, while a vehicle is turning right when $\edit{\sa_s} < 0$. Depending on the magnitude of $\edit{\sa_s}$, \edit{we can categorize the predicted trajectory of the vehicle:}

\edit{\paragraph{$\bullet$ Straight line.}} If the magnitude of $\edit{\sa_s}$ is zero or very small, then one can assume that no turning in the predicted trajectory will be involved and that a vehicle is likely to move in a straight line. In other words, the direction of the velocity vector and the direction of the acceleration vector are the same \edit{(i.e., they are collinear)}. Thus, the vehicle travels in a straight line with its position, velocity, and acceleration predicted according to the standard equations\edit{; for all $t \in [t_\circ, \mathcal T)$,} 
\begin{equation}\label{E:2dttc_traj_linear}
\begin{cases}
    \+{p}(t) \edit{= \+{p}(t; t_\circ)} = \edit{\+{p}_\circ} + \edit{\+{v}_\circ} \edit{(t-t_\circ)} + \frac{1}{2} \edit{\+{a}_\circ} \edit{(t-t_\circ)}^2, \\
    \+{v}(t) \edit{= \+{v}(t; t_\circ)} = \edit{\+{v}_\circ} + \edit{\+{a}_\circ} \edit{(t-t_\circ)}, \\
    \+{a}(t) \edit{\equiv \+{a}_\circ} \\
    \edit{(\vec p(t_\circ),\vec v(t_\circ),\vec a(t_\circ))= (\vec p_\circ,\vec v_\circ,\vec a_\circ)},
\end{cases}
\end{equation}
where the initial condition $\vec p_\circ$ and $\vec v_\circ$ are defined as in \eqref{E:intial_condition} and $\vec a_\circ \Def \dot \zeta(t_\circ)$.
\edit{\paragraph{$\bullet$ Circular.}} If the magnitude of $\edit{\sa_s}$ is not zero and is not small, the constant turn rate and acceleration assumptions mean that \edit{a} vehicle will be moving according to a circular trajectory. Assuming no slip, the radius $r$ and initial angular velocity $\omega_0$ in radians per second along the circular trajectory are given by,
\begin{equation*}
r = \frac{\lVert \edit{\+{v}_\circ} \rVert_2^2}{|\edit{\sa_s}|}, \quad \omega_0 = \frac{\edit{\sa_s}}{\lVert \edit{\+{v}_\circ} \rVert_2}.
\end{equation*}
With $r$ and $\omega_0$, the center of the circular trajectory $\+{c}$ is
\begin{equation*} \+{c} = \edit{\+{p}_\circ} + r \cdot \text{sgn}(\edit{\sa_s}) \cdot \+{v}_\perp \end{equation*}
and the initial phase $\alpha_0$ is
\begin{equation*} \alpha_0 = \text{sgn}(p_{\edit{\circ, y}} - c_y) \cdot \arccos \left( \frac{p_{\edit{\circ, x}} - c_x}{r} \right) \end{equation*}
where the $x$ and $y$ subscripts refer to the component values of \edit{the} respective vector\edit{, making $p_{\edit{\circ, x}}$ and $p_{\edit{\circ, y}}$ to be components of $\+{p}_\circ$.} The $\mathrm{sgn}$ function is defined as 
\begin{equation*}\text{sgn}(z) = \begin{cases} 1, & \text{if } z \geq 0 \\ -1, & \text{if } z < 0 \end{cases}. \end{equation*}
We now define the predicted circular trajectory as
\begin{equation}\label{E:2dttc_traj_circle}
\begin{cases}
     \+{p}(t) \edit{= \+{p}(t; t_\circ)} = \+{c} + r \begin{bmatrix} \cos(\alpha_0 + \omega(t) \edit{(t-t_\circ)}) \\ \sin(\alpha_0 + \omega(t) \edit{(t-t_\circ)}) \end{bmatrix} = \edit{\+{c} + r \begin{bmatrix} \cos(\alpha(t)) \\ \sin(\alpha(t)) \end{bmatrix}}, \edit{\quad t \in [t_\circ, \mathcal T)}\\
    \edit{ \vec p(t_\circ) = \vec p_\circ }.
\end{cases}
\end{equation}
where $\omega(t)$ consider\edit{s} the longitudinal acceleration along the direction of travel\edit{, which is defined as \begin{equation*}
\omega(t) = \begin{cases} 
    \min{\set{0, \omega_0 + \frac{\edit{\sa_f}}{2r} \edit{(t-t_\circ)}}}, \quad\text{if } \omega_0 < 0, \\
    \max{\set{0, \omega_0 + \frac{\edit{\sa_f}}{2r} \edit{(t-t_\circ)}}}, \quad\text{if } \omega_0 > 0. 
\end{cases}
\end{equation*}}
\edit{$\omega(t)$ is defined in such a way that there is no reversal of direction in the circular trajectory and that the vehicle does not move backward, which happens when $\text{sgn}(\sa_f) \neq \text{sgn}(\omega_0)$ or when $\omega(t)$ crosses zero and changes the sign.} The velocity and the acceleration of a vehicle in the predicted circular trajectories are
\edit{\begin{equation}\label{E:2dttc_traj_circle2}
    \begin{cases}
        \+{v}(t) \edit{= \+{v}(t; t_\circ)} = \begin{cases}
            r \begin{bmatrix} -\sin(\alpha(t))(\omega_0 + \frac{\edit{\sa_f}}{r} \edit{(t-t_\circ)}) \\ \cos(\alpha(t))(\omega_0 + \frac{\edit{\sa_f}}{r} \edit{(t-t_\circ)})\end{bmatrix}, \quad\text{if } \omega(t) \neq 0\\
            0, \quad\text{if } \omega(t) = 0\\
        \end{cases} \\ 
        \+{a}(t) \edit{= \+{a}(t; t_\circ)} = \begin{cases}
            r \begin{bmatrix} -\cos(\alpha(t))(\omega_0 + \frac{\edit{\sa_f}}{r} \edit{(t-t_\circ)})^2 - \sin(\alpha(t))(\frac{\edit{\sa_f}}{r})\\ -\sin(\alpha(t))(\omega_0 + \frac{\edit{\sa_f}}{r} \edit{(t-t_\circ)})^2 + \cos(\alpha(t))(\frac{\edit{\sa_f}}{r})\end{bmatrix}, \quad\text{if } \omega(t) \neq 0\\
            0, \quad\text{if } \omega(t) = 0\\
        \end{cases}\\
         \edit{ (\vec v(t_\circ),\vec a (t_\circ)) = (\vec v_\circ,\vec a_\circ)}
    \end{cases}
\end{equation}}
\edit{for }$\edit{ t \in [t_\circ, \mathcal T)}$. 

After we determine $\+{p}_i(t)$ and $\+{p}_j(t)$ for vehicle $i$ and vehicle $j$, the distance between centers of two vehicles is $d_{ij}(t) = \lVert \+{p}_i(t) - \+{p}_j(t) \rVert_2$. We now compute second-order TTC, $\edit{\mathcal T_S}$, as the earliest time \edit{from $t_\circ$} when $d_{ij}(t) - \phi \leq 0$ where $\phi$ denotes the diameter of the circles. That is (see Remark \ref{R:difference_TTC})
\begin{equation}\label{E:collision_second}
    \edit{\mathcal T_S = \mathcal T_S(t_\circ)}\edit{\Def} \edit{\inf}\set{t \edit{\geq t_\circ}: d_{ij}(t) - \phi \edit{\leq} 0}.
\end{equation}
    
Note that, because of the periodicity of the circular equations, there are an infinite number of roots to the equation $d_{ij}(t) - \phi = 0$. \edit{By the definition \eqref{E:collision_second}, $\mathcal T_S$ }corresponds to the smallest positive root \edit{greater than $t_\circ$}.

Although $\edit{\mathcal T_S}$ is well defined to be the smallest positive root to the equation $d_{ij}(t) - \phi = 0$, the process of finding the smallest positive root with classical solvers is difficult. Bounded root finding solvers like \textit{brentq} provide a root that is not necessarily the smallest positive root. In addition, the function needs to have different signs at the interval bounds for the bounded root-finding solvers to provide the root. Due to the limitation of the available built-in function for finding the smallest root among different roots, we develop a new algorithm.

\section{\edit{Second-Order Time-to-Collision Algorithm using Region-based search (STAR)}} \label{sec: algorithm}
In this section, we present \edit{the Second-Order Time-to-Collision Algorithm using Region-based search (STAR)} for second-order TTC. \edit{In this section, we elaborately describe the process of finding the potential regions where depending on the trajectory of the vehicles the collision can happen. In addition, we analytically show that (see Theorem \ref{lemma1} below) searching these neighborhoods is sufficient to find the collision time between the vehicles. By choosing these neighborhoods sufficiently small, we show that our proposed algorithm would be computationally efficient in comparison to the state of the art.}

\edit{\begin{theorem} \label{lemma1}
    Considering vehicles driving according to the dynamics \eqref{E:car_dynamic}. Let us fix the time $t_\circ$ and assume the trajectories of the vehicles on $t \ge t_\circ$ follow the dynamics \eqref{E:2dttc_traj_linear} or \eqref{E:2dttc_traj_circle}. If the collision time $\edit{\mathcal T_S} < \infty$\edit{, defined as in \eqref{E:collision_second}}, there exists \edit{$t_\circ \le \check t < \edit{\mathcal T_S} $}, such that $d_{ij}(t) - \phi$ is strictly decreasing \edit{on $[\check t, \mathcal T_S)$}. 
\end{theorem}}
\begin{proof} 
\edit{For the simplicity of the proof we ignore the dimension $\phi$ of the vehicles \tup{treating them as particles}. The proof remains the same for $\phi >0$.}
\edit{By the definition \eqref{E:collision_second} of $\edit{\mathcal T_S}$ we have that 
\begin{equation} \label{E:no_collision}
d_{ij}(t)>0, \quad t \in [t_\circ, \mathcal T_S). \end{equation}
By the assumption $\mathcal T_S < \infty$ and the continuity of $t \mapsto d_{i,j}(t)$, we have that 
\begin{equation}\label{E:converg_collision}
    \liminf_{t \nearrow \mathcal T_S} d_{ij}(t) =0.
\end{equation}
Let us assume by contradiction to the statement of the theorem that for any time $\check t \ge t_\circ$, there exists $\bar t$ and $\hat{\mathcal T}_S$ with $\bar t < \hat{\mathcal T}_S \in [\check t, \mathcal T_S)$ such that 
\begin{equation} \label{E:dotd_pos}
    \dot d_{ij}(t) \ge 0, \quad t \in [\bar t, \hat {\mathcal T}_S].
\end{equation}
For the rest of the proof, we consider different cases depending on the trajectory of the vehicles (linear or circular). Suppose first that the trajectories of both vehicles are along the straight lines (see Figure \ref{fig:linear_motion}). In this case, under the assumption \eqref{E:dotd_pos} the direction $\vec v_i(t)$ and/or $\vec v_j(t)$ should change on $[\bar t, \hat{\mathcal T}_S]$ which contradicts  the dynamics \eqref{E:2dttc_traj_linear}. Therefore, using \eqref{E:converg_collision}, in the case of linear motion, the claim follows.\\
Now, we consider the case where one vehicle has linear motion and another one moves along a circular pattern (see Figure \ref{fig:linear_circ_motion}) and when the two vehicles follow the circular pattern (see Figure \ref{fig:circ_circ_motion}).  }
\begin{figure}[ht!]
    \centering
    \includegraphics[width=0.5\linewidth]{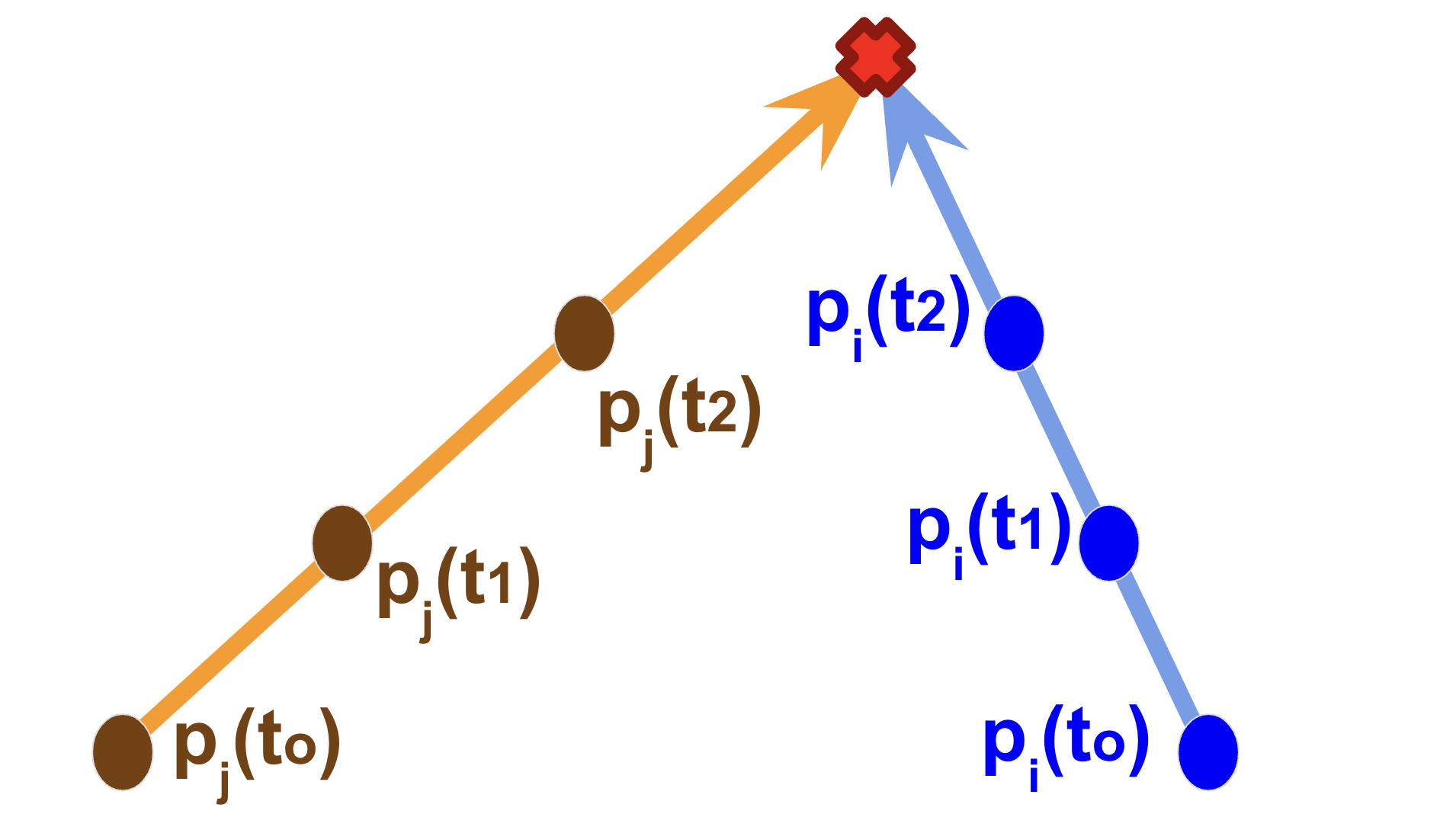}
    \caption{\edit{In the case of linear trajectories, if the collision happens the distance between the vehicles will be decreasing along a linear trajectory.}}
    \label{fig:linear_motion}
\end{figure}
\begin{figure} [ht!]
    \centering
    \begin{subfigure}[b]{0.31\textwidth}
        \centering
        \includegraphics[width=\linewidth]{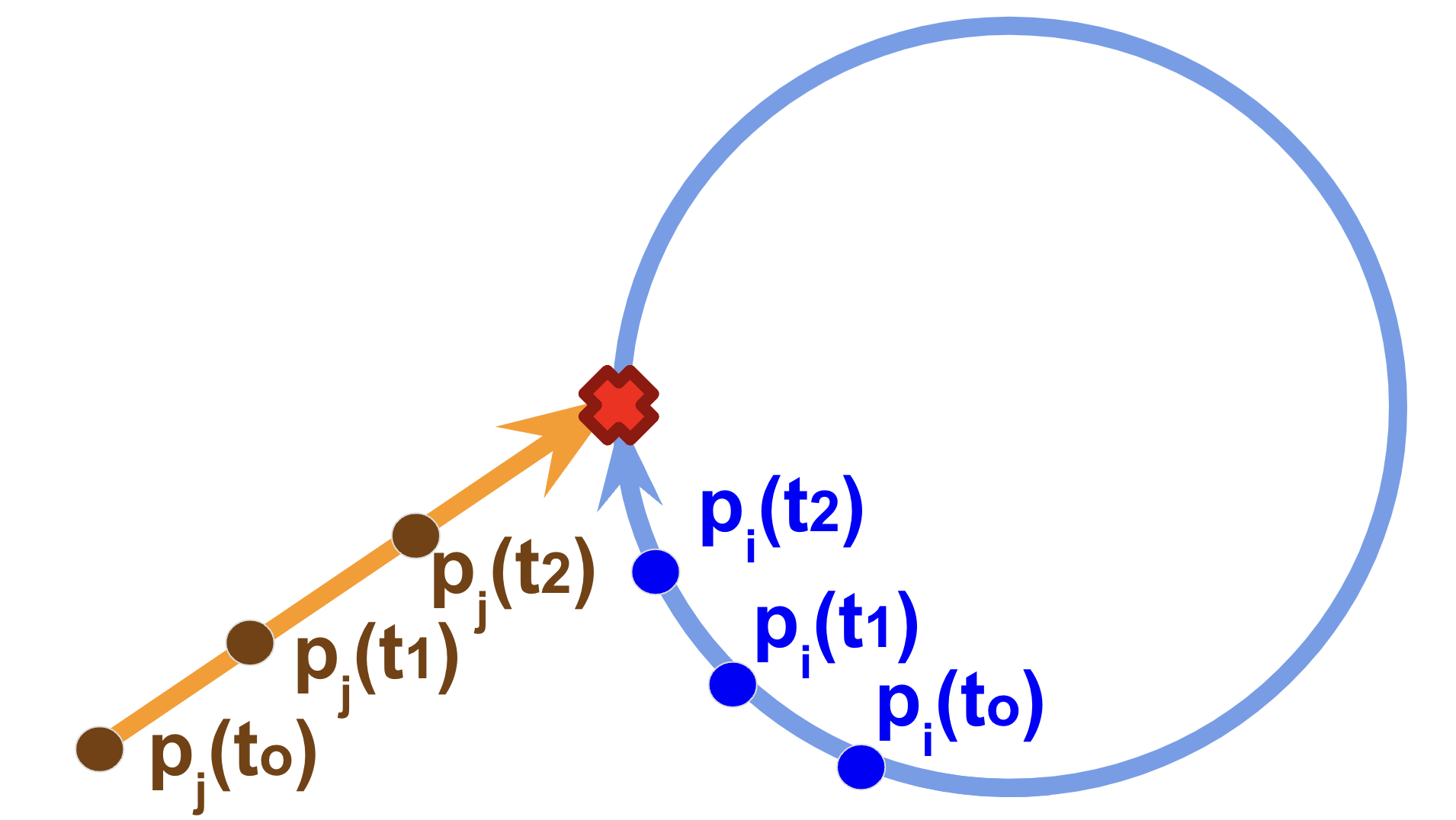}
        \caption{}
        \label{fig:collision_lc1}
    \end{subfigure}
    \begin{subfigure}[b]{0.333\textwidth}
        \centering
        \includegraphics[width=\linewidth]{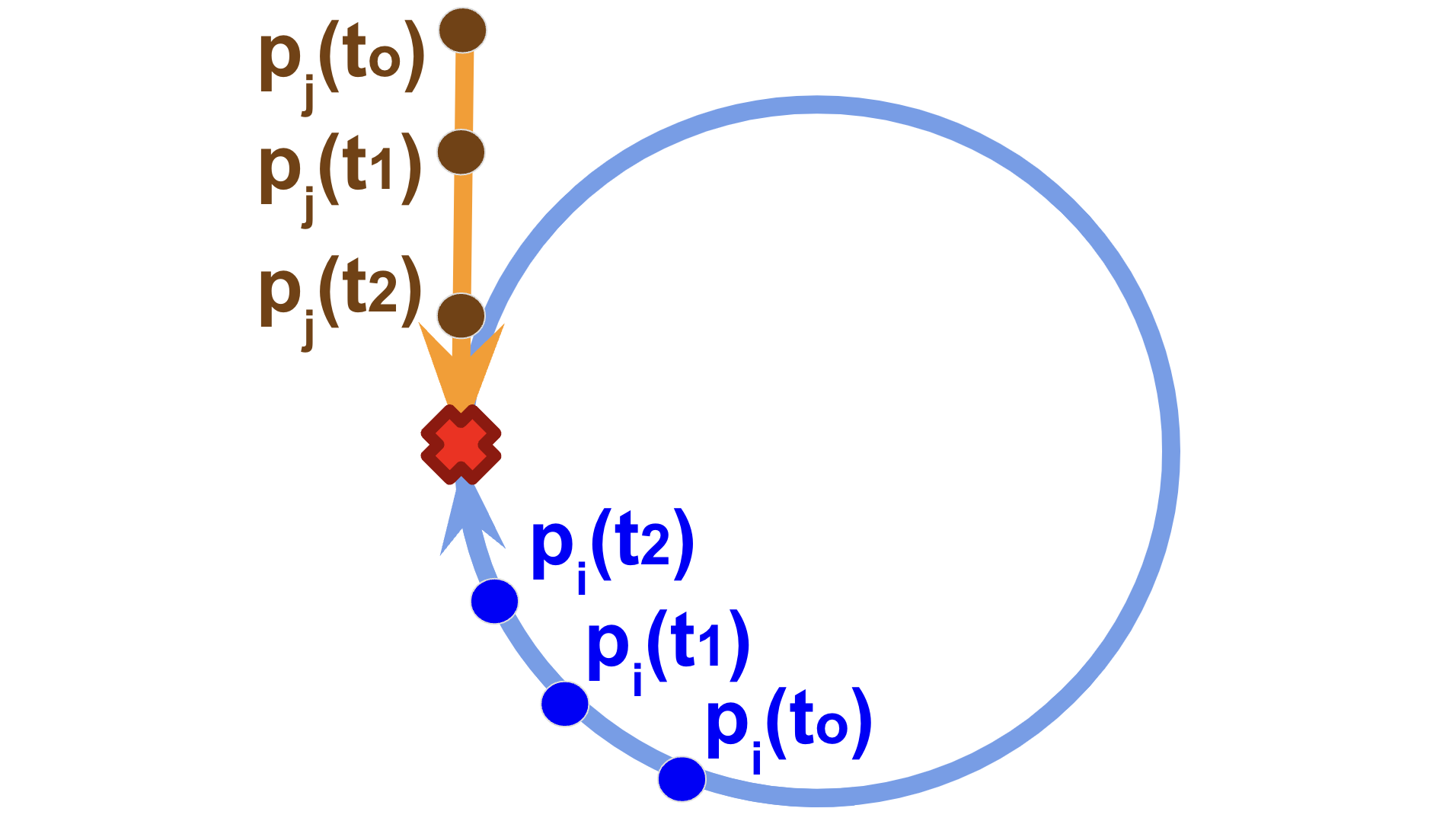}
        \caption{}
        \label{fig:collision_lc2}
    \end{subfigure}
    \begin{subfigure}[b]{0.333\textwidth}
        \centering
        \includegraphics[width=\linewidth]{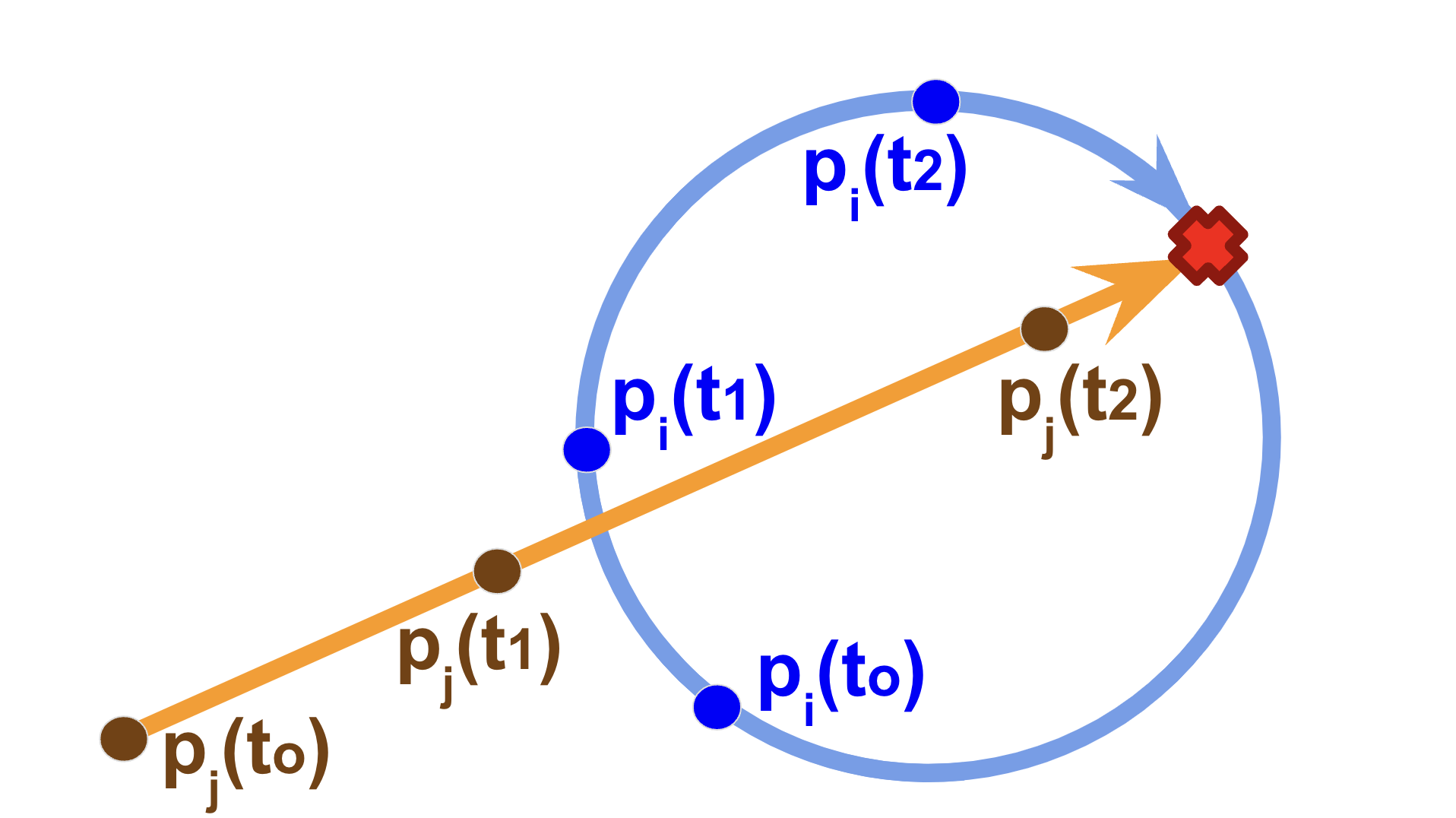}
        \caption{}
        \label{fig:collision_lc3}
    \end{subfigure}
    \caption{\edit{Possible Collision Cases when one vehicle's estimated trajectory is linear while the other vehicle's estimated trajectory is circular.}}
    \label{fig:linear_circ_motion}
\end{figure}

\begin{figure} [ht!]
    \centering
    \begin{subfigure}[b]{0.325\textwidth}
        \centering
        \includegraphics[width=\linewidth]{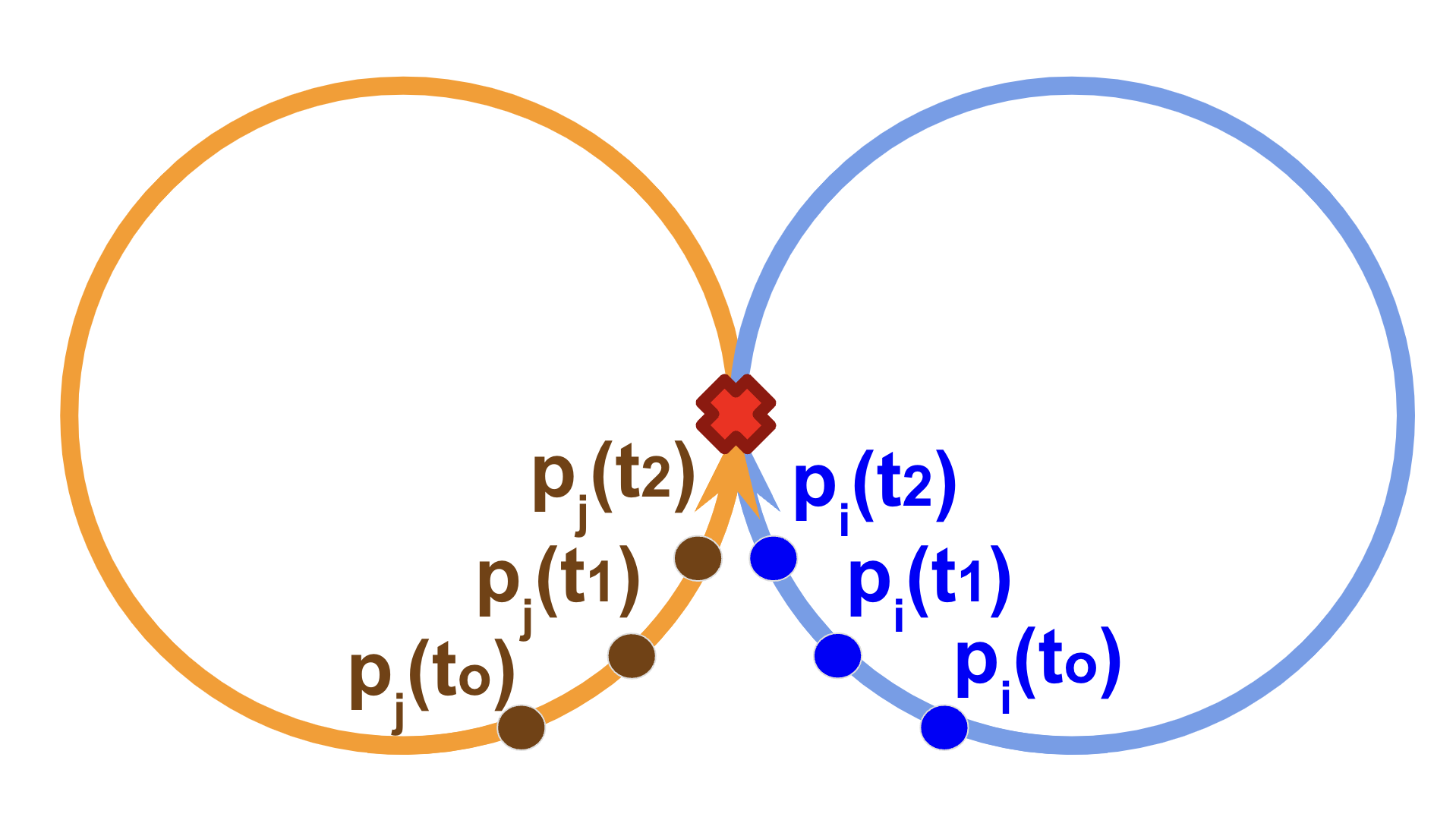}
        \caption{}
        \label{fig:collision_cc1}
    \end{subfigure}
    \begin{subfigure}[b]{0.325\textwidth}
        \centering
        \includegraphics[width=\linewidth]{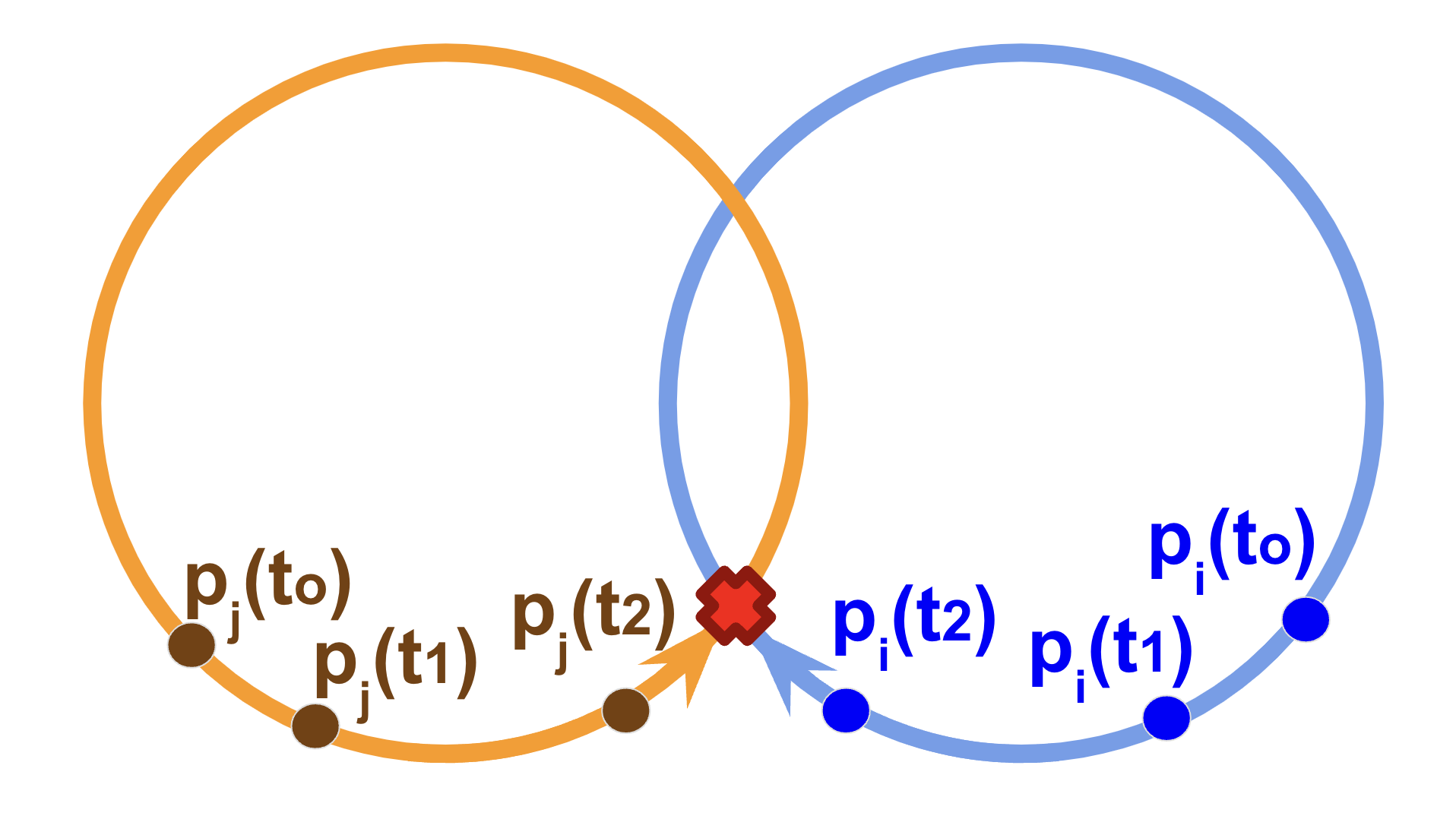}
        \caption{}
        \label{fig:collision_cc2}
    \end{subfigure}
    \begin{subfigure}[b]{0.325\textwidth}
        \centering
        \includegraphics[width=\linewidth]{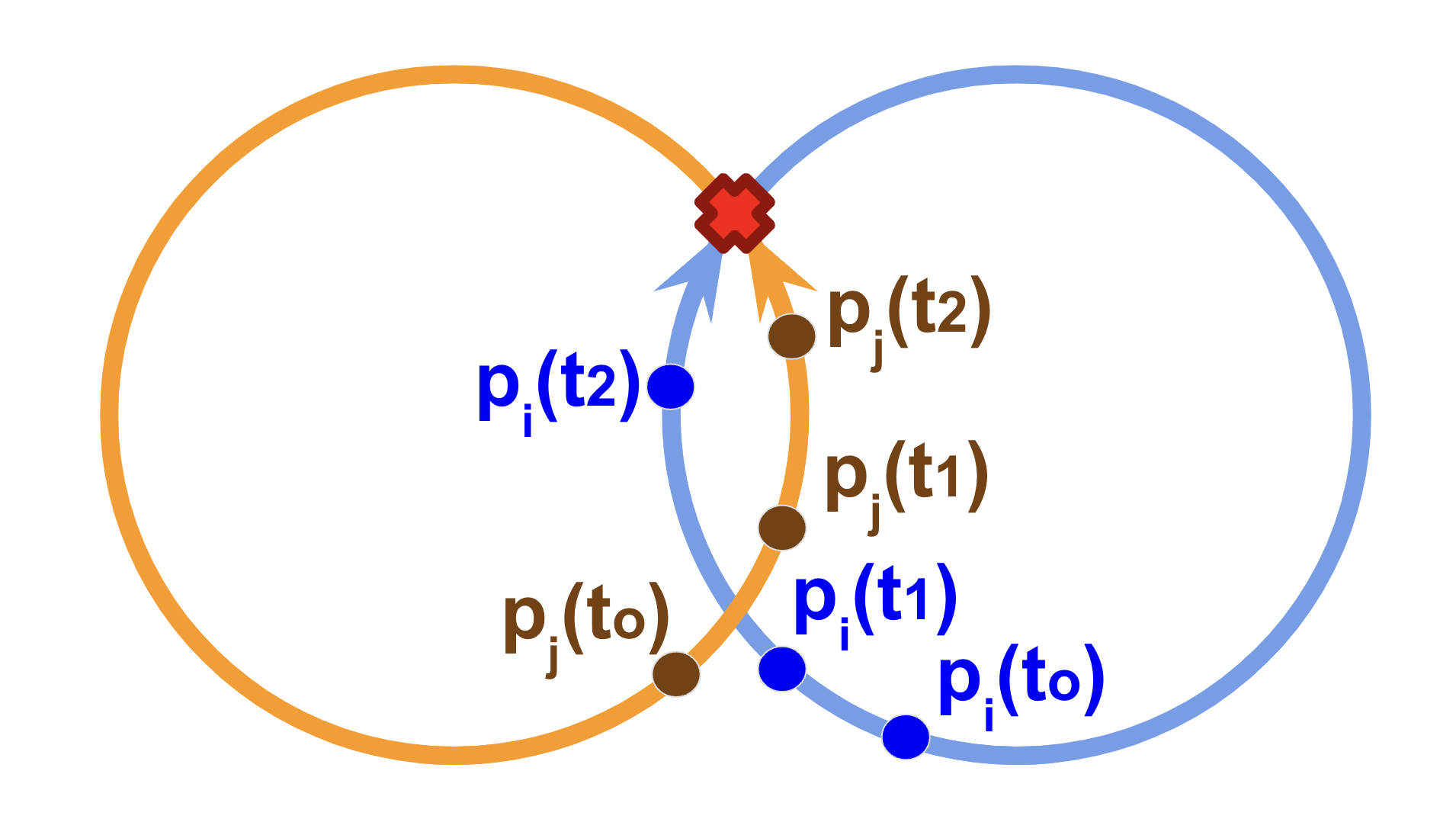}
        \caption{}
        \label{fig:collision_cc3}
    \end{subfigure}
    \caption{\edit{Possible Collision Cases when both vehicles' estimated trajectories are circular.}}
    \label{fig:circ_circ_motion}
\end{figure}
\edit{For both of these cases, the assumption \eqref{E:dotd_pos} implies that there exists a sequence $\set{\mathcal I_m \Def [t_m, t_{m + 1}): m \in \NN}$ of intervals with $t_m \to \mathcal T_S$ as $m \to \infty$ such that 
\begin{equation} \label{E:dir_change} 
\dot d_{ij}(t) \bigmid_{\mathcal I_m} \ge 0, \quad \dot d_{ij}(t) \bigmid_{\mathcal I_{m + 1}} < 0, \quad \forall m \in \NN
\end{equation}
Now, considering the assumption $\mathcal T_S < \infty$, Eq. \eqref{E:dir_change} implies that the $\norm{\vec a}_2$ (the magnitude of the acceleration) can grow unboundedly. This, however, contradicts the properties \eqref{E:properties} of the dynamical model \eqref{E:car_dynamic}. Therefore, \eqref{E:dir_change} and consequently \eqref{E:dotd_pos} violate the properties of the underlying dynamical model and hence cannot hold. Putting all together, the claim of the theorem follows.}
\end{proof}

\edit{\subsection{Algorithmic Structure}}
\edit{In this section, we will elaborate on the details of the proposed algorithm. Algorithm STAR employs the Runge-Kutta-Fehlberg method (RK45) introduced in \cite{fehlberg1969low}, due to its adaptive stepsize control of the local truncated errors \cite{hairer10solving}. Such an adaptive stepsize control in RK45 guarantees accuracy and efficiency in terms of computational complexity \cite{ince1956integration}.}

\noindent \edit{Considering the efficient search region proved in Theorem~\ref{lemma1} and the practical efficiency of RK45, \edit{STAR} shows accurate, precise, and efficient detection of the second-order TTC.}

\edit{Algorithm \ref{alg:cap} shows the pseudocode of the proposed schemes which uses \edit{the} \textit{solve-ivp} function available in the Python environment.}

\begin{algorithm}[!ht]
\caption{\edit{Second-Order Time-to-Collision Algorithm using Region-based search}}\label{alg:cap}
\begin{algorithmic}
\Require Equations $\+{p}_i(t)$, $\+{v}_i(t)$, $\+{a}_i(t)$, $\+{p}_j(t)$, $\+{v}_j(t)$, and $\+{a}_j(t)$, small tolerance $\epsilon$, \edit{size of the search region (radius) $\eps$}, \edit{a diameter of approximation} $\phi$.
\Ensure Second-Order TTC\edit{, $\mathcal T_S$}
\State \edit{Initialize an empty list of intersection points $\+{Q}$}

\State \edit{\textbf{1. Find Intersection Points} (Section \ref{sec: intersection_point})}
\State \edit{Based on $\+{p}_i(t)$, $\+{p}_j(t)$, and $\phi$, Compute all intersection points and append them to $\+{Q}$}

\State \edit{\textbf{2. Define Search Region and Compute Entry and Exit Times} (Section \ref{sec: search_region})}
\For{\edit{each vehicle at each intersection point $\vec q \in \+{Q}$}}
    \State \edit{Define a search region around $\vec q$.}
    \State \edit{Compute entry time and exit time into the search region of $\vec q$ for each vehicle}
\EndFor
\State \edit{Sort $\+{Q}$ based on the minimum entry time regardless of the vehicle.}
\State \edit{\textbf{3. Search for Collision} (Section \ref{sec: search_ttc})}
\State \edit{Initialize $\edit{\mathcal T_S}$ $\gets \infty$}
\For{\edit{each $\vec q \in \+{Q}$ in sorted order}}
    \State \edit{Evaluate collision chances through entry and exit times on a search region of $\vec q$}
    \State \edit{Set a search time interval $\mathcal I$ on a search region of $\vec q$.}
    \For{\edit{each $t$ $\in \mathcal I$}} RK45
        \If{\edit{$d_{ij}(t) \leq \phi + \epsilon$}}
            \State \edit{$\mathcal{T}_S \gets t$ and \textbf{stop the search}}
        \ElsIf{\edit{Termination Criteria Met}}
            \+{break}
        \EndIf
    \EndFor
\EndFor
\State \Return $\edit{\mathcal T_S}$
\end{algorithmic}
\end{algorithm}

\edit{Before we start with the details of the proposed algorithm sketched in Sections \ref{sec: intersection_point}, \ref{sec: search_region}, and \ref{sec: search_ttc}, we need to define a few notations to explain the concepts. 
Let us set the constant time horizon for searching the collision time to be denoted by $\edit{\mathfrak T}$ and we find $\edit{\mathcal T_S \le \mathfrak T}$. In other words, $\frak T$ is the predefined threshold for the algorithmic search of the time to collision.} While $\edit{\mathcal T_S}$ can be found by searching the entire range of $\left[\edit{t_\circ}, \edit{t_\circ +}\edit{\mathfrak T}\right]$, we present a more efficient approach to search for $\edit{\mathcal T_S}$ \edit{under the following consideration:}
\edit{\begin{assumption}\label{ass:circle_one_time}
 Any vehicle only follows a circular trajectory for at most one cycle. This is simply because the circular trajectory reflects the turning dynamics of the vehicles.
\end{assumption}
}

\edit{Based on this assumption, we identify the travel time of each vehicle.}

\edit{\begin{case}
    A vehicle moving in a circular trajectory.
\end{case}}
\edit{Then, the travel time $\tilde{t}$ is the time at which either a vehicle completes one cycle and returns to the original position or when $\omega(t) = 0$. The period of a vehicle can be computed by setting \edit{$\+{p}(t_\circ)$ = $\+{p}_\circ$} = $\+{p}(\tilde{t})$; \edit{i.e., after one cycle, the vehicle will return to the original location $\vec p$.} Therefore, the travel time $\tilde{t}$, the smallest non-negative real number among all possible values, is}
\edit{\begin{equation} \label{E:t_tilde}
    \tilde{t} = 
    \begin{cases}
        \frac{-\omega_0 \pm \sqrt{\omega_0^2 \pm 4(2\pi \frac{\edit{\sa_f}}{2r})}}{\frac{\edit{\sa_f}}{r}} & \text{if } \text{sgn}(\sa_f) = \text{sgn}(\omega_0), \sa_f \neq 0 \\ 
        \min{\set{\frac{-2r\omega_0}{\sa_f}, \frac{-\omega_0 \pm \sqrt{\omega_0^2 \pm 4(2\pi \frac{\edit{\sa_f}}{2r})}}{\frac{\edit{\sa_f}}{r}}}} & \text{if } \text{sgn}(\sa_f) \neq \text{sgn}(\omega_0), \sa_f \neq 0 \\ 
        \max{\set{\pm \frac{2\pi}{\omega_0}}} & \text{if }\sa_f = 0\\ 
    \end{cases}
\end{equation}}

\edit{\begin{case}
    A vehicle moving in a linear trajectory. 
\end{case}}
\edit{In this case, we set the travel time to be either $\mathfrak T$ or the time when $\vec{v}(t) = \elements{0,0}$; that is, the time until the vehicle stops due to deceleration. Therefore, the travel time $\tilde{t}$ is defined as}

\edit{\begin{equation} \label{E:t_tilde_linear}
    \tilde{t} = 
    \begin{cases}
        \frac{\norm{\vec{v}_\circ}}{\norm{\sa}}& \text{if } \sa = -k \vec{v}_\circ, \text{for some $k > 0$}, \\ 
        \mathfrak T & \text{else.}\\ 
    \end{cases}
\end{equation}}
\edit{The first condition in \eqref{E:t_tilde_linear} explains the case where the $\sa$ and $\vec v_\circ$ are in opposite directions and therefore the vehicle eventually might stop at some time.}

Based on the \edit{travel time} of vehicle $i$ and vehicle $j$, the time interval for the search of TTC is $\left[\edit{t_\circ}, \edit{t_\circ+}\tilde{T}\right]$ where 
\edit{\begin{equation} \label{E:col_search}
\tilde{T} \Def \min \set{\tilde{t}_i,\tilde{t}_j, \edit{\mathfrak T}}.
\end{equation}}

\edit{\subsubsection{Intersection Points Computation} \label{sec: intersection_point}}
\edit{In this section, we find the intersection points of two trajectories. In fact, these intersection points will be used in the next section to determine the potential regions where a collision could happen. }

\edit{Given $\edit{(\+{p}_\circ,\+{v}_\circ,\+{a}_\circ)}$ of vehicle $i$ and vehicle $j$, we determine the type of trajectory, either circular or linear trajectory, and corresponding $\+{p}(t)$, $\+{v}(t)$, and $\+{a}(t)$ for each vehicle based on Section \ref{sec:ttc2}. By considering both vehicle $i$ and vehicle $j$, we get $\+{p}_i(t)$, $\+{v}_i(t)$, $\+{a}_i(t)$, $\+{p}_j(t)$, $\+{v}_j(t)$, and $\+{a}_j(t)$. Then three possible combinations of trajectories can be found:}
\edit{\begin{itemize} 
\item Circle-Circle trajectories, 
\item Line-Line trajectories, 
\item Circle-Line trajectories. 
\end{itemize}}


\edit{\paragraph{A. Circle-Circle Trajectories.}} \edit{In the case of circle-circle trajectories, we have the following cases in terms of the possibility of collision. }
\edit{\begin{case} For $\norm{\vec{c}_i - \vec{c}_j} > r_i + r_j + \phi$, then $\mathbf Q = \varnothing$. 
\end{case}}
\edit{\begin{case} For $\norm{\vec{c}_i - \vec{c}_j} \le  r_i + r_j + \phi$ a collision can happen.
\end{case}}
\edit{In this case, we first determine the midpoint of the line segment connecting the centers of two circular trajectories which is}

\edit{\begin{equation} \label{eq:mid_point_circle}
    \frak{m}_{ij} = \vec{c}_i + \frac{r_i^2 - {r_j}^2 + \norm{\vec{c}_i - \vec{c}_j}_2^2}{2\norm{\vec{c}_i - \vec{c}_j}_2^2} \vec{c}_j - \vec{c}_i.
\end{equation}}

\edit{Given $\pmb{\kappa} = \vec{c}_j - \vec{c}_i = \elements{\kappa_x,\kappa_y}$, we find $\pmb{\kappa}_\perp$, the perpendicular vector to $\pmb{\kappa}$, which is}

\edit{\begin{equation} \label{eq:perpendicular vector_circle}
\pmb{\kappa}_\perp = \frac{\elements{-\kappa_y,\kappa_x}}{\norm{\pmb{\kappa}}_2}.
\end{equation}}

\edit{The intersection points are located at} 

\edit{\begin{equation} \label{eq:circle_intersections}
\mathbf Q = \set{\frak{m}_{ij} \pm \sqrt{r_i^2 - \left(\frac{r_i^2 - {r_j}^2 + \norm{\pmb{\kappa}}_2^2}{2\norm{\pmb{\kappa}}_2}\right)^2} \pmb{\kappa}_\perp}.
\end{equation}}

\edit{\begin{case} For $\vec{c}_i = \vec{c}_j$ and $r_i = r_j$, two circular trajectories are completely overlapped. In this case, whenever two vehicles encounter for the first time, a collision happens. This collision time can be found explicitly. 
\end{case}}

\edit{\paragraph{B. Line-Line Trajectories.}}
\edit{When both vehicles are determined to move in linear trajectories, we first begin with formulating the linear equation in a Cartesian coordinate system describing the linear trajectory moving along the direction of the velocity $\vec{v_\circ} = \elements{v_{\circ,x},v_{\circ,y}}$ at an initial position $\vec{p_\circ} = \elements{p_{\circ,x},p_{\circ,y}}$. The linear equation in a Cartesian coordinate system is}
\edit{\begin{equation*}
    \begin{cases}
        y = \frac{v_{\circ,y}}{v_{\circ,x}} x + \left( p_{\circ,y} - \frac{v_{\circ,y}}{v_{\circ,x}} p_{\circ,x} \right)& \text{if } v_{\circ,x} \neq 0, \\ 
        x = p_{\circ,x} & \text{if } v_{\circ,x} = 0.
    \end{cases}
\end{equation*}}
\edit{Here, the first equation formulates the non-vertical linear trajectory while the second equation explains the vertical trajectory of vehicles.} 

\edit{After determining the linear equation for both vehicles' trajectories, we now find the intersection points. It should be noted that, taking the direction of the linear trajectories, i.e., $\text{sgn}(v_{\circ,x})$, is a key point in finding proper intersections as it represents the motion of the vehicles starting from the initial position. We denote the initial position $\vec{p}_{i,\circ} = \elements{p_{i,\circ,x},p_{i, \circ,y}}$ and the initial velocity $\vec{v}_{i,\circ} = \elements{v_{i,\circ,x},v_{i, \circ,y}}$ of vehicle $i$ (similarly for the $j$-th vehicle, $\vec{p}_{j,\circ}$ and  $\vec{v}_{j,\circ}$). 
The following cases explore the possible collision points in linear-linear trajectories.}
\edit{\begin{case} Both vehicles are moving in a vertical line trajectory in a Cartesian coordinate system, i.e., $v_{i,\circ,x} = 0$ and $v_{j,\circ,x} = 0$. In this case, either two trajectories coincide for which the collision would happen at their first encounter and can be found explicitly, or they do not intersect at all. \end{case}}

\edit{\begin{case} One of the trajectories is non-vertical, i.e., $v_{i,\circ,x}$ or $v_{j,\circ,x}$ is non-zero.
\end{case}}
\edit{Then the intersection point is at}
\edit{\begin{equation*}
   \begin{split}
        &\mathbf Q =  \set{\elements{p_{j,\circ,x},\frac{v_{i,\circ,y}}{v_{i,\circ,x}} p_{j,\circ,x} + \left( p_{i,\circ,y} - \frac{v_{i,\circ,y}}{v_{i,\circ,x}} p_{i, \circ,x} \right)}}\quad \text{if } v_{i,\circ,x} \neq 0, v_{j,\circ,x} = 0, \\
        &\mathbf Q =  \set{\elements{p_{i,\circ,x},\frac{v_{j,\circ,y}}{v_{j,\circ,x}} p_{i,\circ,x} + \left( p_{j,\circ,y} - \frac{v_{j,\circ,y}}{v_{j,\circ,x}} p_{j, \circ,x} \right)}}\quad \text{if } v_{i,\circ,x} = 0, v_{j,\circ,x} \neq 0. 
    \end{split}
\end{equation*}}

\edit{\begin{case}  None of the vehicles are moving along the vertical linear trajectory in a Cartesian coordinate system. \end{case}}
\edit{In this case, the intersection point can be found at}
\edit{\begin{equation*}
    \mathbf Q = \set{\elements{
    \frac{b_j - b_i}{m_i - m_j}, \, m_i \frac{b_j - b_i}{m_i - m_j} + b_i
    }}
\end{equation*}}
\edit{where}
\edit{\begin{align*}
    m_i = \frac{v_{i,\circ,y}}{v_{i,\circ,x}}, \quad
    m_j = \frac{v_{j,\circ,y}}{v_{j,\circ,x}}, \quad
    b_i = p_{i,\circ,y} - m_i p_{i,\circ,x}, \quad
    b_j = p_{j,\circ,y} - m_j p_{j,\circ,x}.
\end{align*}}

\edit{\paragraph{C. Circle-Line Trajectories.}} \edit{Now, let us assume that vehicle $i$ is determined to be moving along the linear trajectory while vehicle $j$ is determined to be moving along the circular trajectory with center $\vec{c}_j = \elements{c_{j,x}, c_{j,y}}$ and radius $r_j$. }
\edit{\begin{case}
    The line trajectory is a non-vertical line.
\end{case}}
\edit{In this case, the intersection points of two trajectories are defined by} 
\edit{\begin{equation*}
    \begin{split}
        &\mathbf Q = \varnothing, \quad \text{\scriptsize{if $B^2 - 4AC < 0$}} \\
        &\mathbf Q = \set{\elements{ \frac{-B}{2A}, \frac{v_{i,\circ,y}}{v_{i,\circ,x}} \frac{-B}{2A} + \left( p_{i,\circ,y} - \frac{v_{i,\circ,y}}{v_{i,\circ,x}} p_{i,\circ,x}\right)}},\quad \text{\scriptsize{if $B^2 - 4AC = 0$}} 
    \end{split}
\end{equation*}}
\edit{\begin{multline*}
\mathbf Q = \set{\elements{\frac{-B \pm \sqrt{B^2 - 4AC}}{2A}, \frac{v_{i,\circ,y}}{v_{i,\circ,x}} \frac{-B \pm \sqrt{B^2 - 4AC}}{2A} + \left( p_{i,\circ,y} - \frac{v_{i,\circ,y}}{v_{i,\circ,x}} p_{i,\circ,x}\right)}},\\
\text{\scriptsize{if $B^2 - 4AC = 0$}}
\end{multline*}}
\edit{where} 
\edit{\begin{align*}
A &= 1+\left(\frac{v_{i,\circ,y}}{v_{i,\circ,x}}\right)^2, \\
B &= 2\left(\left(\frac{v_{i,\circ,y}}{v_{i,\circ,x}}\right) (p_{i,\circ,y} - \frac{v_{i,\circ,y}}{v_{i,\circ,x}} p_{i,\circ,x}-c_{j,y})-c_{j,x}\right), \\
C &= c_{j,x}^2+(p_{i,\circ,y} - \frac{v_{i,\circ,y}}{v_{i,\circ,x}} p_{i,\circ,x}-c_{j,y})^2-r_j^2). 
\end{align*}}
\edit{\begin{case}
   The linear trajectory is a vertical line. 
\end{case}}
\edit{For this case, the intersection points can happen at \begin{equation*} \mathbf Q = \set{\elements{p_{i,\circ,x},c_{j,y} \pm \sqrt{r_j^2-(p_{i,\circ,x}-c_{j,x})^2}}}.
\end{equation*}}

\edit{It is, again, important to make sure the direction of the line is considered when finding the intersection points in the case of circle-line trajectories.}

\edit{\subsubsection{Search Region and Boundary Hitting Times} \label{sec: search_region}}
\edit{In Section \ref{sec: intersection_point}, the intersection points between any two trajectories are calculated. In this section, using such intersection points, we proceed to find regions in which the collision between vehicles will be searched by our proposed algorithm. In particular, given the intersection point $\vec q \in \+Q$, a search region is defined as the neighborhood $\ball(\vec q, \eps)$ for a sufficiently small radius $\eps$ considering the \textit{trajectory of only one vehicle}. By investigating the \textit{interaction of another vehicle} with such search region the collision can occur if both vehicles can be positioned in the same neighborhood in a specific time period. Consequently, the time-to-collision can be efficiently determined by searching only at the moments when each vehicle enters and exits the search region.}

\edit{We define the entry time to the search region of the intersection point $\vec q$ as follows:}
\edit{\begin{equation*}
\begin{split}
    T_{\vec q}^{\text{enter}} \Def \inf \set{t \ge t_\circ: \mathfrak{d}(\vec p(t), \vec q) < \eps}
\end{split}
\end{equation*}}
\edit{where $\frak d(\vec p(t), \vec q)$ defines the distance between the two points along the trajectory $t \mapsto \vec p(t)$ and calculations will be elaborated in \eqref{E: t_enter_lienar} and Remark \ref{R:distance}. In other words, $T_{\vec q}^{\text{enter}}$ represents the first time that a vehicle approaches $\ball(\vec q, \eps)$. Similarly, the exit time from this neighborhood can be defined as}
\edit{\begin{equation*}
\begin{split}
    T_{\vec q}^{\text{exit}} \Def \inf \set{t \ge T_{\vec q}^{\text{enter}}: \mathfrak{d}(\vec p(t), \vec q) \ge \eps}
\end{split}
\end{equation*}}
\edit{We now find $T_{\vec q}^{\text{enter}}$ and $T_{\vec q}^{\text{exit}}$ for each type of trajectory; i.e. linear and circular.}
\edit{\paragraph{$\bullet$ Linear Trajectory.}} \edit{Finding the intersection points between $\ball(\vec q, \eps)$ and the linear trajectory, $\vec q'$ and/or $\vec q''$, is by $\vec{v}_*$ (defined as in \eqref{E:vstar}), $\vec q'$ = $\vec q - \eps \vec{v}_*$ and $\vec q''$ = $\vec q + \eps \vec{v}_*$. Then $T_{\vec q}^{\text{enter}}$ from the initial position $\vec{p_\circ}$ is as follows.}

\edit{\begin{equation} \label{E: t_enter_lienar}
T_{\vec q}^{\text{enter}} = 
    \begin{cases}
         \min{\set{\tilde{T}, \frac{-(\vec{v}\centerdot\vec{\sa})\pm \sqrt{(\vec{v}\centerdot \vec{\sa})^2 - 2||\sa||^2((\vec{p_\circ}-\vec q')\centerdot \vec{\sa})}}{||\sa||^2}}},& \text{if } \norm{\sa} \neq 0 \\
         \min{\set{\tilde{T}, \frac{((\vec{p_\circ}-\vec q')\centerdot \vec v)}{||\vec {v}||^2}}},& \text{if } \norm{\sa} = 0.    \end{cases}
\end{equation}}

\edit{Similarly, $T_{\vec q}^{\text{exit}}$ are defined as}

\edit{\begin{equation*}
T_{\vec q}^{\text{exit}} =
    \begin{cases}
         \min{\set{\tilde{T}, \frac{-(\vec{v}\centerdot\vec{\sa})\pm \sqrt{(\vec{v}\centerdot \vec{\sa})^2 - 2||\sa||^2((\vec{p_\circ}-\vec q'')\centerdot \vec{\sa})}}{||\sa||^2}}}& \text{if } \norm{\sa} \neq 0, \\
        \min{\set{\tilde{T}, \frac{((\vec{p_\circ}-\vec q'')\centerdot \vec v)}{||\vec {v}||^2}}}& \text{if } \norm{\sa} = 0.   
    \end{cases}
\end{equation*}}

\edit{\paragraph{$\bullet$ Circular Trajectory.}} \edit{For the description of the process of finding the search region and search time interval for the circular trajectory, let us consider a fixed time $t \in [t_\circ, \tilde T]$, the position $p(t)$ of a vehicle and a fixed point $\vec{o} = [o_x,o_y]$ on a circular trajectory with center $\vec{c}$ and radius $r$. We begin with finding the angle $(\vec p(t), \vec o) \mapsto \Theta(t) \in [0, \pi]$ between $\vec{o}$ and $\vec{p}(t)$,} 
\edit{\begin{equation} \label{eq:theta}
\Theta(\vec p(t), \vec o)= \arccos{\frac{(p_{x}(t)-c_{x})(o_x-c_{x})+(p_{y}(t)-c_{y})(o_y-c_{y})}{\sqrt{(p_{x}(t)-c_{x})^2+(p_{y}(t)-c_{y})^2}\sqrt{(o_x-c_{x})^2+(o_y-c_{y})^2}}}. 
\end{equation}}

\edit{As $\Theta(\vec p(t), \vec o)$ represents the smallest angle between $\vec{o}$ and $\vec{p}(t)$ without the consideration of the vehicle's moving direction, we find the angle, $\vartheta(t)$, between two points based on $\vec{v}_\circ = [v_{\circ,x},v_{\circ,y}]$. }
\edit{\begin{equation} \label{eq:real_theta}
    \vartheta(\vec p(t), \vec o) = 
        \begin{cases} 
        \Theta(\vec p(t), \vec o) & \text{if } (\vec{p}(t)-\vec{c}) \times (\vec{o}-\vec{c})) > 0, \vec{v}_\circ \centerdot (\vec{p}_{\perp}-\vec{c}_{\perp}) > 0 \\
        2\pi - \Theta(\vec p(t), \vec o) & \text{if } (\vec{p}(t)-\vec{c}) \times (\vec{o}-\vec{c})) > 0, \vec{v}_\circ \centerdot (\vec{p}_{\perp}-\vec{c}_{\perp}) < 0  \\
        \Theta(\vec p(t), \vec o) & \text{if } (\vec{p}(t)-\vec{c}) \times (\vec{o}-\vec{c})) < 0, \vec{v}_\circ \centerdot (\vec{p}_{\perp}-\vec{c}_{\perp}) < 0  \\
        2\pi - \Theta(\vec p(t), \vec o) & \text{if } (\vec{p}(t)-\vec{c}) \times (\vec{o}-\vec{c})) < 0,  \vec{v}_\circ \centerdot (\vec{p}_{\perp}-\vec{c}_{\perp}) > 0
    \end{cases}
\end{equation}}
\edit{where $\vec{p}_{\perp} = \elements{- p_{\circ, y}, p_{\circ,x}}$ and $\vec{c}_{\perp} = \elements{- c_{y}, c_{x}}$. The cross product, $(\vec{p}(t)-\vec{c}) \times (\vec{o}-\vec{c}))$, determines the relative physical location of $\vec{p}(t)$ with respect to $\vec{o}$, either counter-clockwise (positive) or clockwise (negative). The dot product, $\vec{v}_{\circ} \centerdot (\vec{p}_{\perp}-\vec{c}_{\perp})$ determines if the vehicle is moving in the counter-clockwise direction (positive) or clockwise direction (negative) on the circular trajectory. }
\edit{\begin{remark}
    In \eqref{eq:real_theta}, while the velocity vector $\vec v$ changes in time, due to the definition \eqref{E:t_tilde}  determining whether the vehicle is moving clockwise or counter-clockwise can be simply done by merely considering the initial direction $\vec v_\circ$.
\end{remark}}
\edit{\begin{remark} \label{R:distance} For a fixed time $t \in [t_\circ, \tilde T]$, while we do not directly use the concept of the distance $\frak d(\vec p(t), \vec q)$ to calculate the hitting times, we can explicitly calculate such distance based on the value of $\vartheta(t)$. In particular, given $\vartheta(t)$ from \eqref{eq:real_theta}, we have that $\mathfrak{d}(\vec p(t), \vec q) = r \vartheta(\vec p(t), \vec q)$. 
\end{remark}}
\edit{Now we can find $T_{\vec q}^{\text{enter}}$ and $T_{\vec q}^{\text{exit}}$ when a vehicle moves along the circular trajectory. We first find $\vec q' = \elements{q_x', q_y'}$ and/or $\vec q'' = \elements{q_x'', q_y''}$, the intersection point between $\ball(\vec q, \eps)$ with $\vec q = \elements{q_x, q_y}$ and the circular trajectory with center $\vec c = \elements{c_x,c_y}$ and radius $r$, using \eqref{eq:mid_point_circle}, \eqref{eq:perpendicular vector_circle}, and \eqref{eq:circle_intersections}. Using \eqref{eq:theta} and \eqref{eq:real_theta}, we find the $\vartheta_{\eps^-}$ and $\vartheta_{\eps^+}$, the smallest and the largest angle between the initial position $\vec p (t_\circ)$ and $\ball(\vec q, \eps)$ respectively while considering the direction of the movement.}

\edit{Then, $T_{\vec q}^{\text{enter}}$ is, }
\edit{\begin{equation*} 
    T_{\vec q}^{\text{enter}} = 
    \begin{cases}
    \min{\set{\tilde{T}, \frac{-\omega_0 \pm \sqrt{\omega_0^2 \pm 4(\frac{\edit{\sa_f}}{2r} \vartheta_{\eps^-} )}}{\frac{\edit{\sa_f}}{r}}}} & \text{if } \sa_f \neq 0 ,\\ 
    \min{\set{\tilde{T}, \pm \frac{\vartheta_{\eps^-}}{\omega_0}}} & \text{if } \sa_f = 0 ,
    \end{cases}
\end{equation*}}
\edit{while $T_{\vec q}^{\text{exit}}$ is}
\edit{\begin{equation*}
    T_{\vec q}^{\text{exit}} = 
    \begin{cases}
    \min{\set{\tilde{T}, \frac{-\omega_0 \pm \sqrt{\omega_0^2 \pm 4(\frac{\edit{\sa_f}}{2r} \vartheta_{\eps^+} )}}{\frac{\edit{\sa_f}}{r}}}} & \text{if } \sa_f \neq 0, \\ 
    \min{\set{\tilde{T}, \pm \frac{\vartheta_{\eps^+}}{\omega_0}}} & \text{if } \sa_f = 0 .
    \end{cases}
\end{equation*}}
\edit{We repeat the process of finding the time each vehicle enters and exits the search region of each intersection point. Then we sort the intersection points in $\+ Q$ based on the minimum entering time into the search region regardless of which vehicle. }


\edit{\subsubsection{Determining Time to Collision} \label{sec: search_ttc}}

\edit{Let us fix point $\vec q \in \vec Q$ as the first possible collision point between vehicles. In this section, we will find the collision time of vehicles colliding at $\vec q$ (the process of finding the TTC for other points in $\vec Q$ would follow the same procedure). More precisely, suppose vehicle $i$ enters the $\ball(\vec q, \eps)$ at $T_{i, \vec q}^{\text{enter}}$ and exits at $T_{i, \vec q}^{\text{exit}}$ (similarly for vehicle $j$). Therefore, the search time interval for determining TTC can be defined by:}
\edit{\begin{equation*}
   \mathcal I \Def \begin{cases}
        [t_\circ + T_{j, \vec q}^{\text{enter}}, t_\circ + T_{i, \vec q}^{\text{exit}}] & \text{if } T_{i, \vec q}^{\text{enter}} < T_{j, \vec q}^{\text{enter}} < T_{i, \vec q}^{\text{exit}}, \\
        [t_\circ + T_{i, \vec q}^{\text{enter}}, t_\circ + T_{j, \vec q}^{\text{exit}}] & \text{if } T_{j, \vec q}^{\text{enter}} < T_{i, \vec q}^{\text{enter}} < T_{j, \vec q}^{\text{exit}},\\
        \text{Aborting the search at $\vec q$} &, \text{otherwise}.
    \end{cases}
\end{equation*}}
\edit{We use RK45 to search TTC over the interval $\mathcal I$. Then the following search and termination criteria would be evaluated: } 
\edit{\paragraph{Searching Criteria.}} \edit {If at a time $\tau \in \mathcal I$, the condition}
\edit{\begin{equation*}
    d(\vec p_i(\tau), \vec p_j(\tau)) \leq \phi, 
\end{equation*}}
\edit{is satisfied, then $\mathcal T_S \approx \tau$ would be considered as an approximation of time to collision and the algorithm will be terminated.}
\edit{\paragraph{Termination Criteria.}} \edit{If the condition
\begin{equation*} d(\vec p_i(t), \vec p_j(t)) > \phi,
\end{equation*}} 
    \edit{and one of the following holds}
    \edit{\begin{itemize}
        \item $\mathfrak{d}(\vec p_i(t), \vec q) > \eps$, or 
        \item $\mathfrak{d}(\vec p_j(t), \vec q) > \eps$, or 
        \item $\dot d(\vec p_i(t), \vec p_j(t)) \geq 0$,
    \end{itemize}}
\noindent \edit{then the search will be terminated as the collision cannot happen at intersection point $\vec q$, and we move on to the next intersection point in $\+ Q$. If the search is terminated based on the termination criteria at every intersection point in $\+ Q$, we set $\mathcal T_S = \infty$ as the imminent collision is not expected in the near future.}


\edit{
\begin{remark}\label{R:difference_TTC}
    In the definitions \eqref{E:collision_first} and \eqref{E:collision_second}, the TTC is defined considering the starting point $t = 0$. In evaluating the performance of the proposed algorithms in the following section, for the purpose of comparison, we calculate the TTC starting from any fixed point $t_\circ$. 
\end{remark}
}

\section{Evaluation of the Algorithm} \label{sec: evaluation}

We now assess the performance of \edit{STAR}. For the evaluation of \edit{STAR}, we use the Apple M2 Max Processor with 38 GPU Core, 12 CPU Core, and 96GB Memory. The implementation of \edit{STAR} and the step-by-step simulation used for evaluation are CPU-bound and are not parallelized. We perform 1001 trials to evaluate the performance of the algorithm by looking at the algorithm error and computation time. We randomly selected initial conditions as \edit{$p_{\circ,x}, p_{\circ,y} \in (-20,20)$, $v_{\circ,x}, v_{\circ,x} \in (-1,1)$; $\sa_{\circ,x}, \sa_{\circ,x} \in (-0.1, 0.1)$} for \edit{$\+{p}_\circ = \elements{p_{\circ,x}, p_{\circ,y}}$, $\+{v}_\circ = \elements{v_{\circ,x}, v_{\circ,y}}$, and $\+{a}_\circ = \elements{\sa_{\circ,x}, \sa_{\circ,y}}$} of both vehicle $i$ and vehicle $j$. Throughout the paper, we set $\phi$ to be 5. The maximum time interval for the search of TTC is $\edit{\mathfrak T}=$100 for the algorithm evaluation. 
\subsection{Algorithmic Error}
We first compute the error between the second-order TTC computed from the step-by-step simulation and \edit{STAR} for each trial. Step-by-step simulation checks for the collision at every time step without the use of the numerical method. The timestep size of the step-by-step simulations is 1E-5. For trials with $\edit{\mathcal T_S}$ $\neq \infty$, we subtract the TTC computed from the step-by-step simulation from the TTC computed by \edit{STAR}. Table \ref{tab:alg_error} presents the absolute error of \edit{STAR} across \edit{$\frak n = $} \edit{140} trials compared to the step-by-step simulations with different step sizes. 
\begin{table} [!ht]
    \centering
    \begin{tabular}{|c|c|c|c|c|c}\hline
         Step Size & Absolute Mean Error $\frak{\bar{x}}$ & Standard Deviation $\frak s$\\ \hline
         1E-5 & \edit{2.927E}-6 & \edit{3.320E-6}\\ \hline
    \end{tabular}
    \caption{Algorithmic Error of \edit{STAR} across \edit{140} trials.}
    \label{tab:alg_error}
\end{table}

\begin{figure} [ht!]
    \centering
    \includegraphics[width=0.6\linewidth]{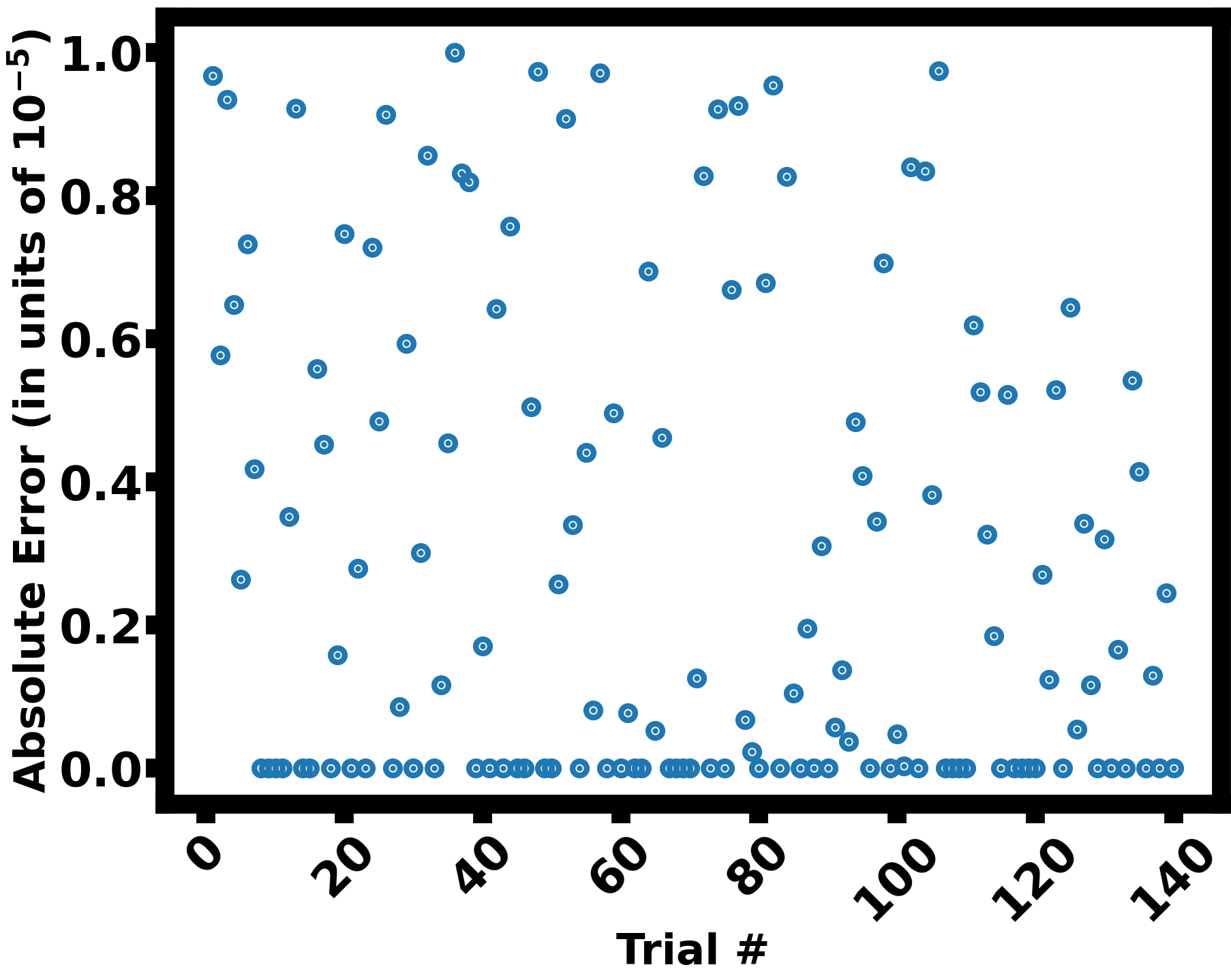}
    \caption{Scatter plot of the absolute algorithmic error of \edit{STAR} across \edit{140} trials compared to the step-by-step simulation with a step size of 1E-5. Each point represents the error of each trial where the x-axis \edit{indicates} the trial number and the y-axis indicates the error.}
    \label{fig:error_alg}
\end{figure}

Figure \ref{fig:error_alg} shows the plot of algorithmic error of \edit{STAR} in all \edit{140} trials compared to the step-by-step simulations.

From Figure \ref{fig:error_alg}, we can visually confirm that most errors are very close to 0 with the magnitude of errors in all \edit{140} trials being less than \edit{the} corresponding stepsize of the step-by-step simulation. Since all errors in all \edit{140} trials are less than the stepsize, we conclude that the error is acceptable and that the algorithm is capable of producing accurate second-order TTC. 

The results in Table~\ref{tab:alg_error} clearly show that the mean absolute error is smaller than half of the stepsize, which is the error upper bound of the step-by-step simulations. This is also supported by our observations from Figure~\ref{fig:error_alg}. We now perform the hypothesis test to check the magnitude of the absolute algorithmic error with respect to the timestep size of the step-by-step simulation. As there is no population mean and no population standard deviation, we need to perform a one-sample t-test. The hypotheses are: 






\begin{itemize}
    \item Null Hypothesis (\(H_0\)): \(\mu \geq \mu_0 = \text{1E-5}\)
    \item Alternative Hypothesis (\(H_1\)): \(\mu < \mu_0 = \text{1E-5}\)
\end{itemize}

\noindent where $\mu$ is the population mean of the sample and $\mu_0$ is the hypothesized population mean. The t-statistic is calculated as 
\begin{equation*}
\frak{t} = \frac{\frak{\bar{x}} - \mu_0}{\frak s / \sqrt{\edit{\frak n}}} = \frac{(\edit{2.927E-6}) - (1E-5)}{(\edit{3.320E-6}) / \sqrt{\edit{140}}} \approx \edit{-25.21}.
\end{equation*}
The corresponding p-values are approximately 0 for all stepsizes. Therefore, \edit{based on the hypothesis test, one} can claim that the absolute algorithm error is \edit{smaller} than \edit{the} corresponding stepsize. $\hfill\diamond$

Therefore, we can conclude that \edit{STAR} is capable of producing accurate TTC values.$\hfill\diamond$


\subsection{Computation Time}

In addition, we record the running time of each trial to compute the average running time of the \edit{STAR}. It is important to check this because if the running time of \edit{STAR} is greater than the step-by-step simulation, it weakens the reason to use \edit{STAR} to compute the second-order TTC or $\edit{\mathcal T_S}$. Table \ref{tab:alg_runtime} presents the average running time of \edit{STAR} and the average running time of the step-by-step algorithm.

\begin{table} [!ht]
    \centering
    \begin{tabular}{cr@{$\,\pm\,$}lc}\toprule
         & \multicolumn{2}{c}{Running time (sec)} & Average speedup \\ \midrule
         \edit{STAR} & \edit{0.0028} & \edit{0.0047} & 1$\times$ \\
         Simulation (Stepsize: 1E-2) & \edit{0.0388} & \edit{0.0285} & \edit{14}$\times$ \\
         Simulation (Stepsize: 1E-3) & \edit{0.3974} & \edit{0.3140} & \edit{142}$\times$ \\
         Simulation (Stepsize: 1E-5) & \edit{36.3991} & \edit{26.5099} & \edit{13000}$\times$ \\\bottomrule
    \end{tabular}
    \caption{Overall running time for estimating the second-order TTC using \edit{STAR} or step-by-step simulations with stepsize of 1E-2, 1E-3, and 1E-5. Means and standard deviations are obtained over 1001 trials.}
    \label{tab:alg_runtime}
\end{table}
Table \ref{tab:alg_runtime} shows the performance difference between the step-by-step simulation and \edit{STAR} in terms of computation time \edit{across $\frak n = $ 1001 trials}. Based on the average time and its standard deviation, \edit{STAR} outperforms the step-by-step simulation as \edit{STAR} produces the second-order TTC 2000 times faster. 

\begin{figure} [ht!]
    \centering
    \begin{subfigure}[b]{0.47\textwidth}
        \centering
        \includegraphics[width=\linewidth]{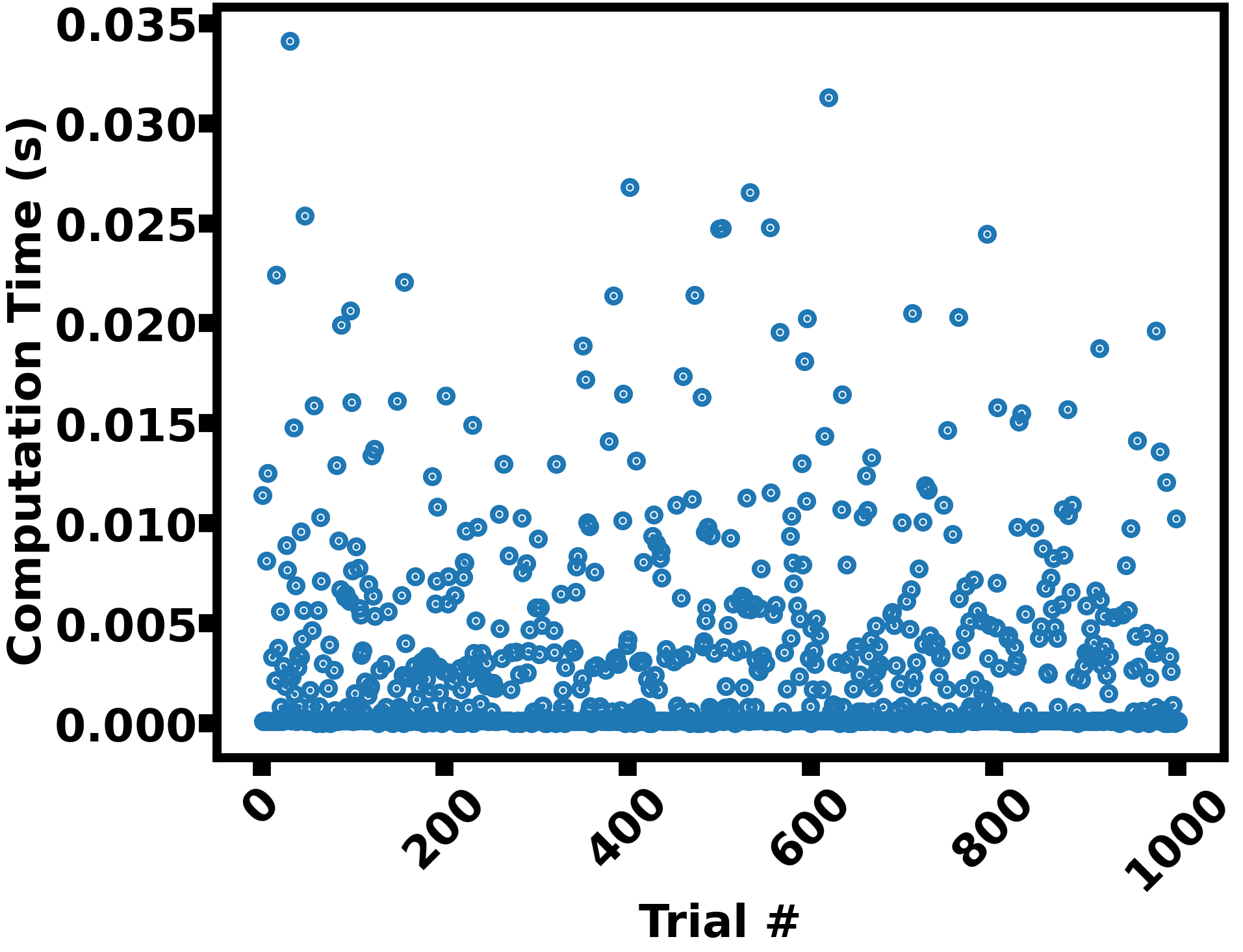}
        \caption{\edit{STAR}}
        \label{fig:comptime_alg}
    \end{subfigure}
    \begin{subfigure}[b]{0.47\textwidth}
        \centering
        \includegraphics[width=\linewidth]{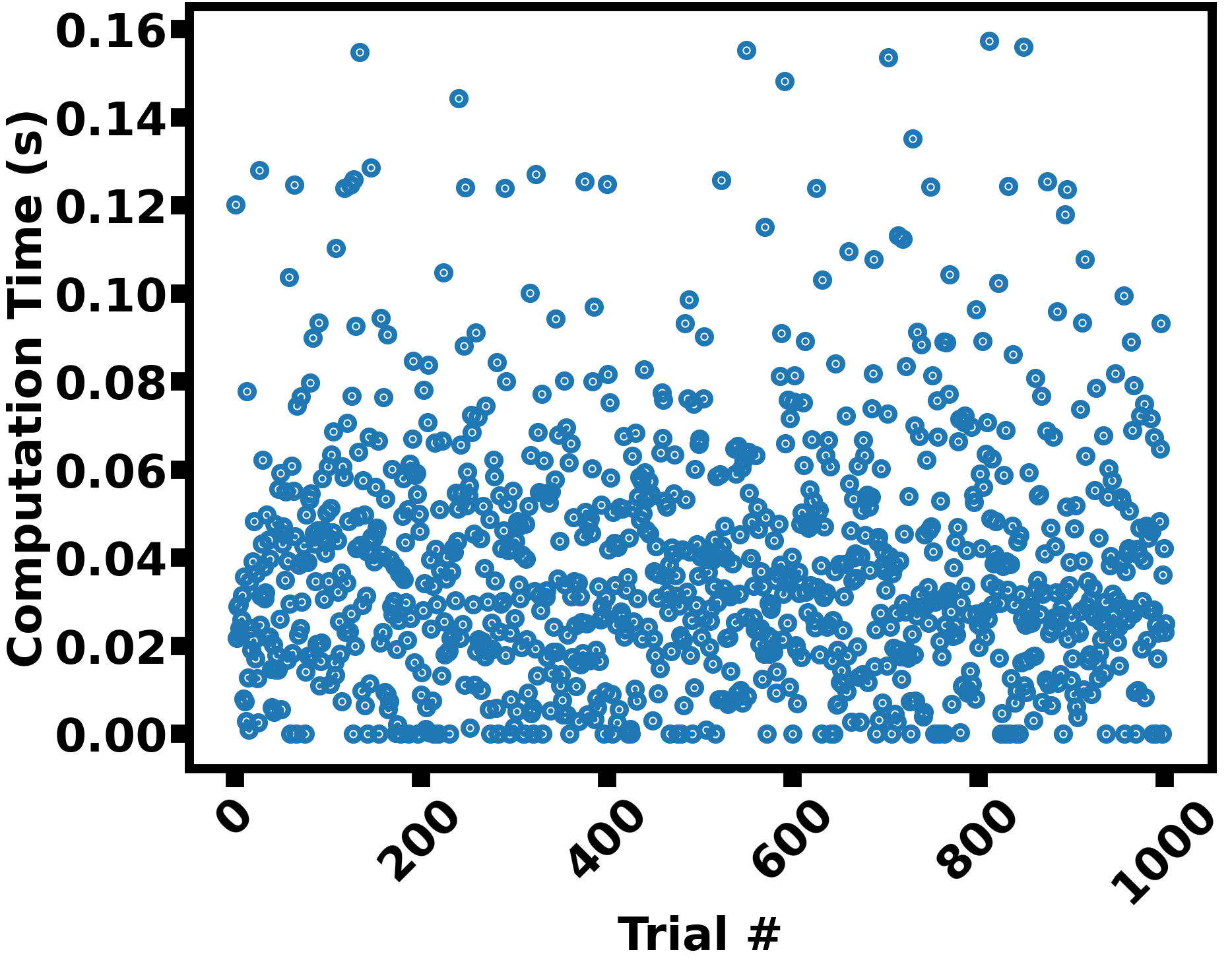}
        \caption{Step-by-Step Simulation(Stepsize: 1E-2)}
        \label{fig:comptime_sim1}
    \end{subfigure}
    \begin{subfigure}[b]{0.47\textwidth}
        \centering
        \includegraphics[width=\linewidth]{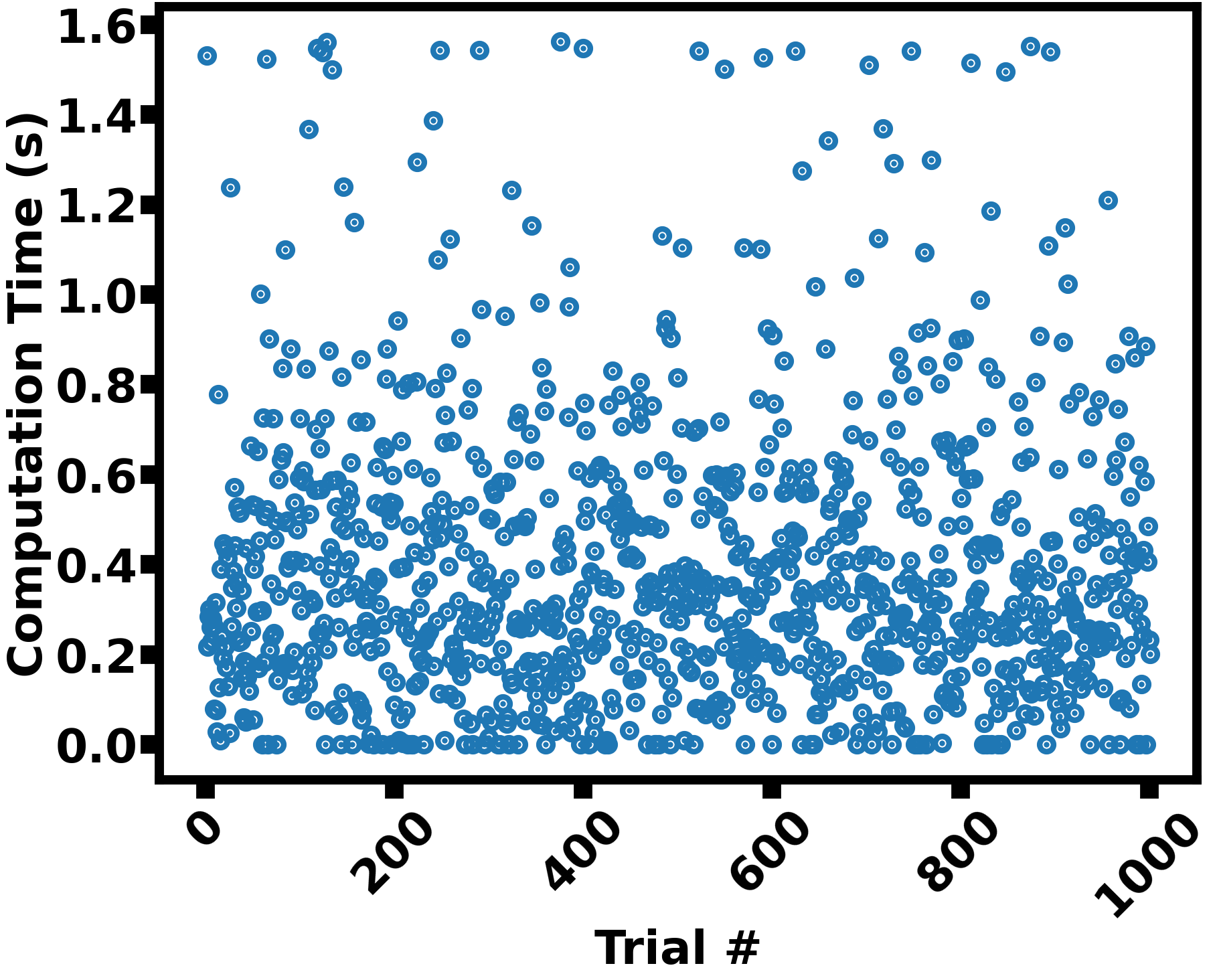}
        \caption{Step-by-Step Simulation(Stepsize: 1E-3)}
        \label{fig:comptime_sim2}
    \end{subfigure}
    \begin{subfigure}[b]{0.47\textwidth}
        \centering
        \includegraphics[width=\linewidth]{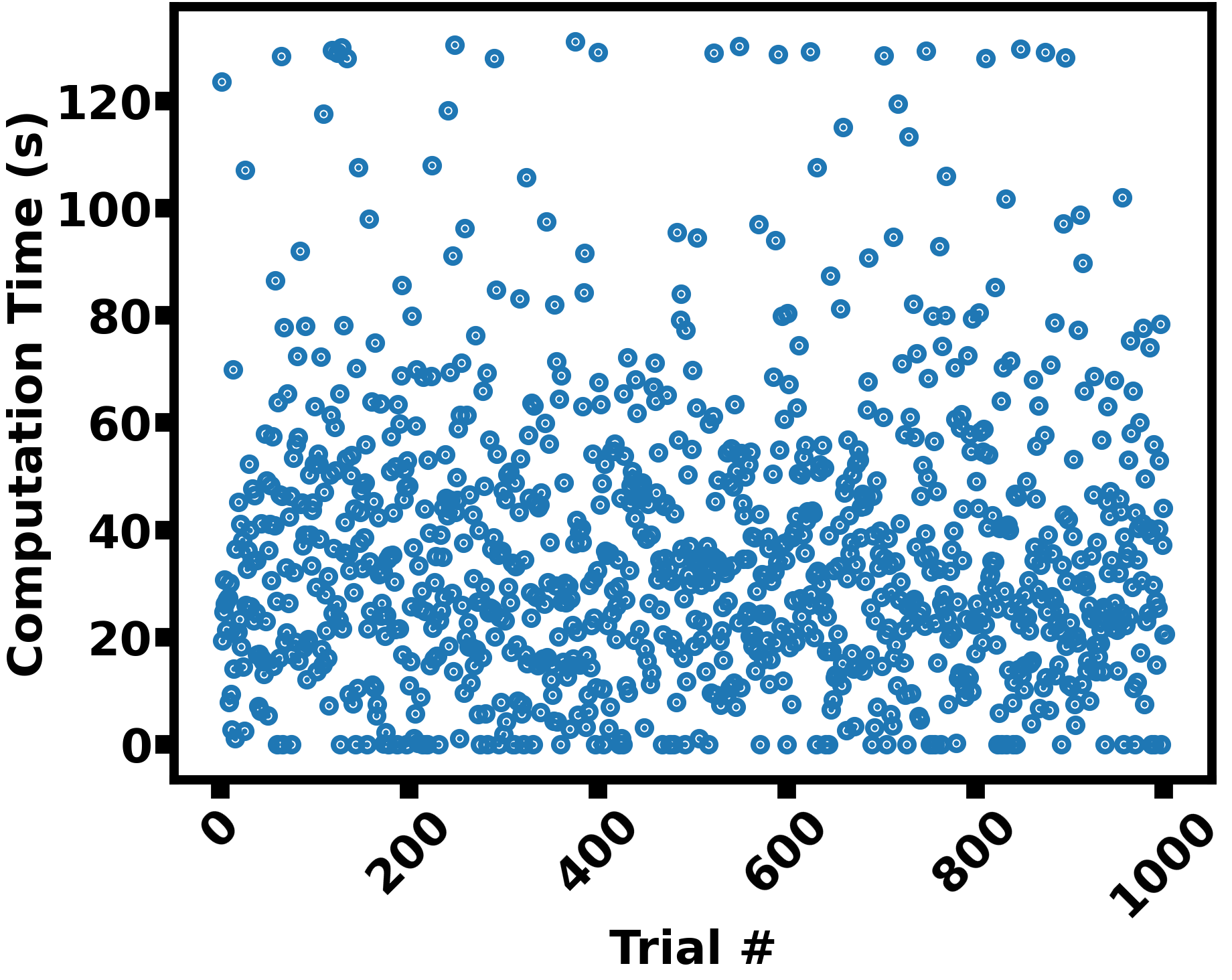}
        \caption{Step-by-Step Simulation(Stepsize: 1E-5)}
        \label{fig:comptime_sim3}
    \end{subfigure}
    \caption{Computation time of \edit{STAR} and the step-by-step simulations in stepsize of 1E-2, 1E-3, and 1E-5. Each point represents the computation time of either \edit{STAR} or step-by-step simulation where the x-axis \edit{indicates} the trial number and the y-axis indicates the computation time in seconds.}
    \label{fig:comptime}
\end{figure}

We additionally visually compare the computation time of \edit{STAR} and the step-by-step simulation by plotting the computation time of both approaches for each trial. Figure \ref{fig:comptime} shows the computation time of \edit{STAR} (Figure \ref{fig:comptime_alg}) and the computation time of the step-by-step simulation (Figure~\ref{fig:comptime}(b)-(d)).

From Figure \ref{fig:comptime}, we can visually confirm that the step-by-step simulation requires more computation time compared to \edit{STAR} for all 1001 trials. When making the comparison between \edit{STAR} and the step-by-step simulation with \edit{a} step size of 1E-5 \edit{in} Figure \ref{fig:comptime}, it is important to keep \edit{in} mind that the y-axis representing the computation time of the step-by-step simulation is in units of $10^{2}$ while it is not for \edit{STAR}. 

\begin{figure} [ht!]
    \centering
    \begin{subfigure}[b]{0.303\textwidth}
        \centering
        \includegraphics[width=\linewidth]{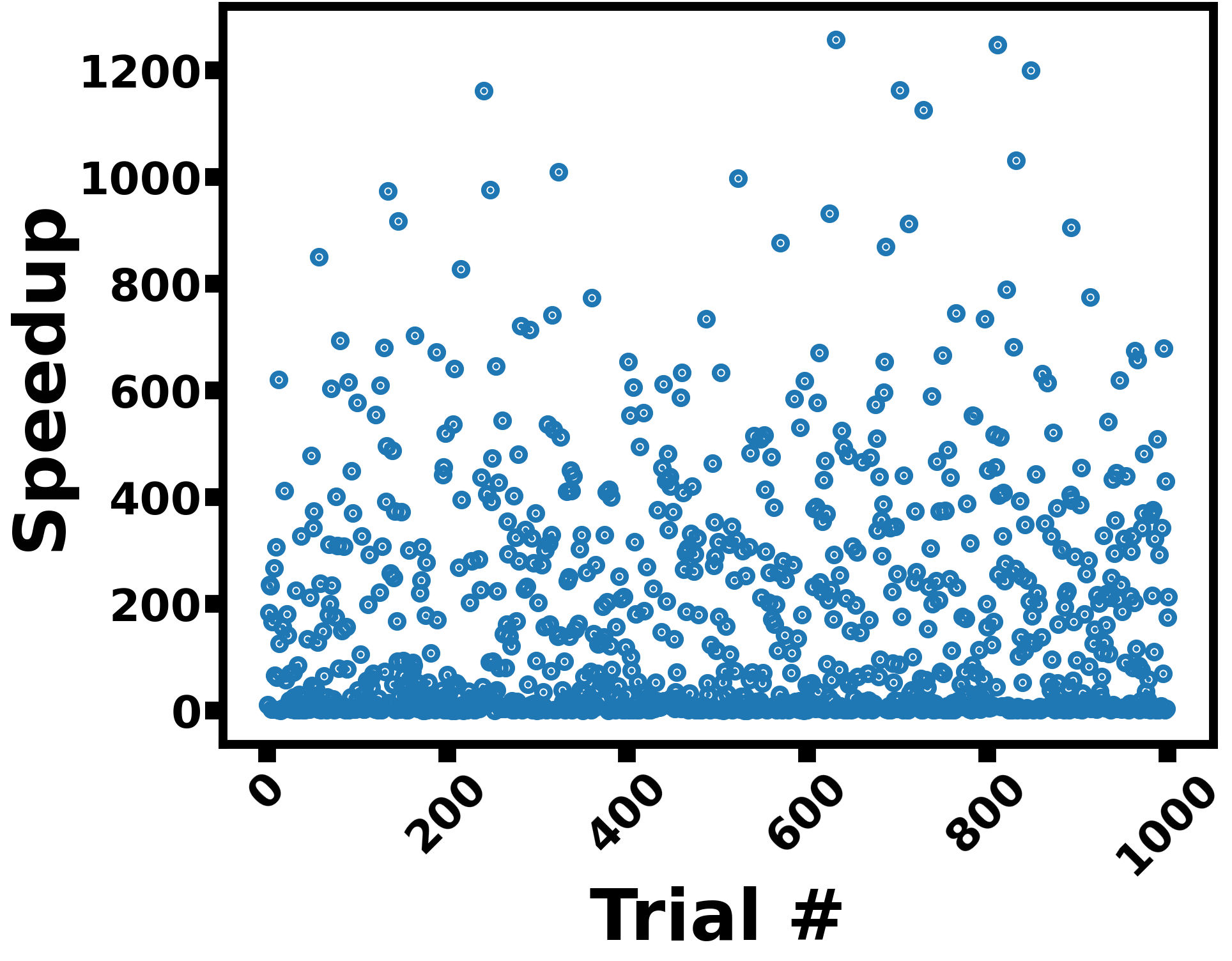}
        \caption{Step-by-Step Simulation (Stepsize: 1E-2)}
        \label{fig:comptime_speedup_sim1}
    \end{subfigure}
    \begin{subfigure}[b]{0.31\textwidth}
        \centering
        \includegraphics[width=\linewidth]{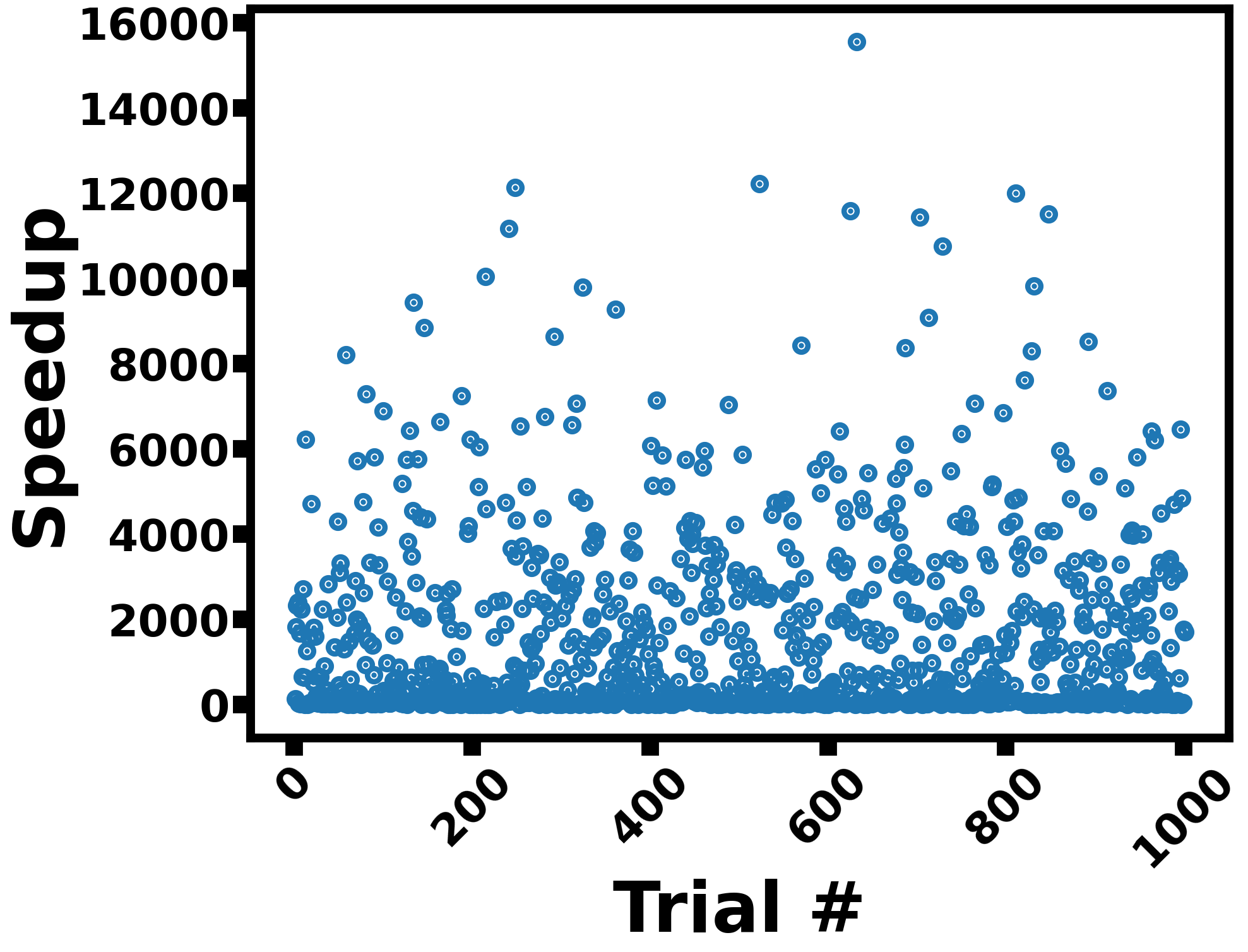}
        \caption{Step-by-Step Simulation (Stepsize: 1E-3)}
        \label{fig:comptime_speedup_sim2}
    \end{subfigure}
    \begin{subfigure}[b]{0.31\textwidth}
        \centering
        \includegraphics[width=\linewidth]{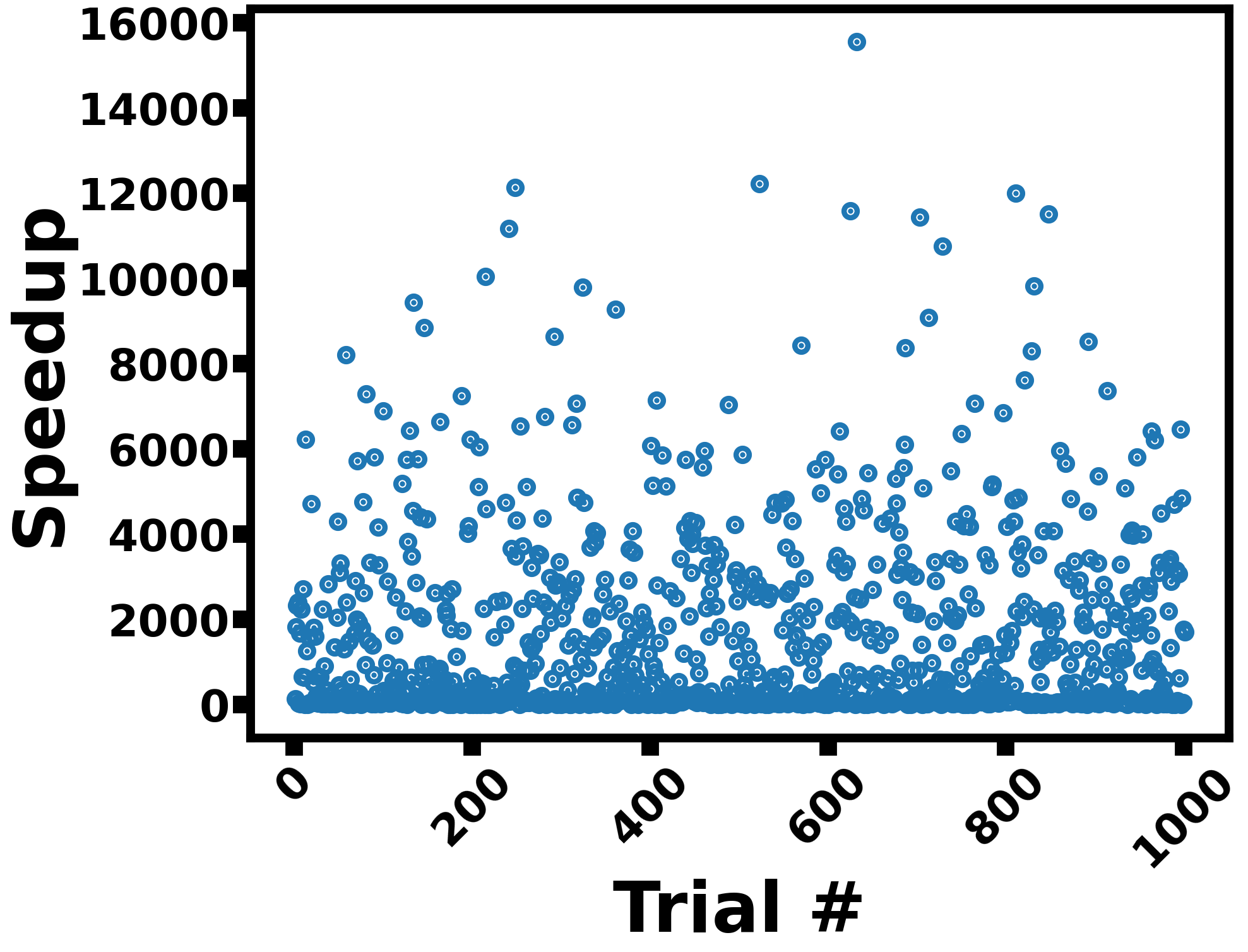}
        \caption{Step-by-Step Simulation (Stepsize: 1E-5)}
        \label{fig:comptime_speedup_sim3}
    \end{subfigure}
    \caption{Speedup due to \edit{STAR} compared to step-by-step simulation across 1001 trials. Each point represents the speedup of each trial. Speedup is defined as the computation time of step-by-step simulation divided by the computation time of \edit{STAR}. The x-axis \edit{indicates} the trial number and the y-axis indicates the number of times speedup.}
    \label{fig:comptime_comparison}
\end{figure}

We additionally find the computation time difference between the step-by-step simulation and \edit{STAR} by subtracting the computation time of \edit{STAR} \edit{from} the computation time of the step-by-step simulation for each trial. \edit{The l}arger the difference, \edit{the} larger the difference between the step-by-step simulation and \edit{STAR} in terms of the computation time. Figure \ref{fig:comptime_comparison} shows the speedup in computation time due to \edit{STAR} compared to step-by-step simulation.

Figure \ref{fig:comptime_comparison} supports the observation made in Figure \ref{fig:comptime} that the computation of the second order TTC almost always takes more time through step-by-step simulation compared to \edit{STAR}. The speedup is mostly greater than 0 for all 1001 trials, especially when the stepsize for the step-by-step simulation gets smaller. The speedup is less pronounced for the 1E-2~sec stepsize, but this is achieved by a corresponding increase in the error. At that stepsize and larger, one risks missing sudden \edit{near-crash} vehicle interactions, especially those occurring at high speeds. However, \edit{STAR} achieves its computation speed without making those trade-offs.

In addition to visual observation from Figure \ref{fig:comptime} and \ref{fig:comptime_comparison}, we conduct the statistical tests, a two-sample t-test, to prove that \edit{STAR} produces the TTC faster than the step-by-step simulation.
\begin{itemize}
    \item Null Hypothesis (\(H_0\)): \(\mu_1 \geq \mu_2 \)
    \item Alternative Hypothesis (\(H_1\)): \(\mu_1 < \mu_2\)
\end{itemize}
The t-statistic for the two-sample t-test is 
\begin{equation*}
\frak t = \frac{\frak{\bar{x}_1} - \frak{\bar{x}_2}}{\sqrt{\frac{\frak{s_1^2}}{\edit{\frak n}_1} + \frac{\frak{s_2^2}}{\edit{\frak n}_2}}}  = \begin{cases}
    \frac{\edit{0.0028}s - \edit{0.0388}s}{\sqrt{\frac{(\edit{0.0047}s)^2}{1001} + \frac{(\edit{0.0285}s)^2}{1001}}} = \edit{-39.432}, \quad \text{Stepsize = 1E-2}\\
    \frac{\edit{0.0028}s - \edit{0.3974}s}{\sqrt{\frac{(\edit{0.0047}s)^2}{1001} + \frac{(\edit{0.3140}s)^2}{1001}}} = \edit{-39.755}, \quad \text{Stepsize = 1E-3}\\
    \frac{\edit{0.0028}s - \edit{36.3991}s}{\sqrt{\frac{(\edit{0.0047}s)^2}{1001} + \frac{(\edit{26.5099}s)^2}{1001}}} = \edit{-43.438}, \quad \text{Stepsize = 1E-5}\\
\end{cases}
\end{equation*} 
Based on the t-statistic, the p-value is approximately 0 for all stepsizes. We, therefore, reject the null hypothesis and state that \edit{STAR} can produce the TTC value faster than the simulation. $\hfill\diamond$


\section{Analysis Through Numerical Simulation} \label{sec: analysis}


In this section, we simulate two vehicles in various realistic settings at an intersection. For each scenario, we compare the performance of the first-order TTC (1D-TTC) and the second-order TTC (2D-TTC) in terms of how well it captures the risk of collision correctly. 

Given the setup of a scenario, we compute the first-order TTC and the second-order TTC based on the initial condition of vehicles in a scenario. In addition to the TTC value under both schemes, we plot the estimated trajectories of vehicles \edit{based on initial status at time $t_\circ$} and plot the $\textit{d}_{ij}-\phi$ \edit{along} the estimated trajectories of vehicles under both schemes. If the collision is to happen \edit{on estimated trajectories of vehicles}, $\textit{d}_{ij}-\phi$ is less than or equal to 0. \edit{We additionally provide t}he color gradient map\edit{, showing} the location of vehicles \edit{along the estimated trajectory from $t_\circ$}. \edit{On the color gradient map, i}f the estimated trajectories of two vehicles are close to each other and have similar colors, a potential collision might be occurring.

\edit{For the numerical simulations, w}e especially look at the scenarios where turning is involved as both the first-order TTC and the second-order TTC perform similar\edit{ly} in scenarios with both vehicles moving in straight line trajectories. Each scenario is simulated for 10 seconds \edit{where we compute TTC under both schemes at each timestep using estimated trajectories over a 20-second horizon based on positions, velocities, and accelerations of vehicles at that timestep.}


\subsection{Scenario 1: Both Vehicles Turning at an Intersection}

\begin{figure} [ht!]
    \centering
  \includegraphics[width=0.9\linewidth]{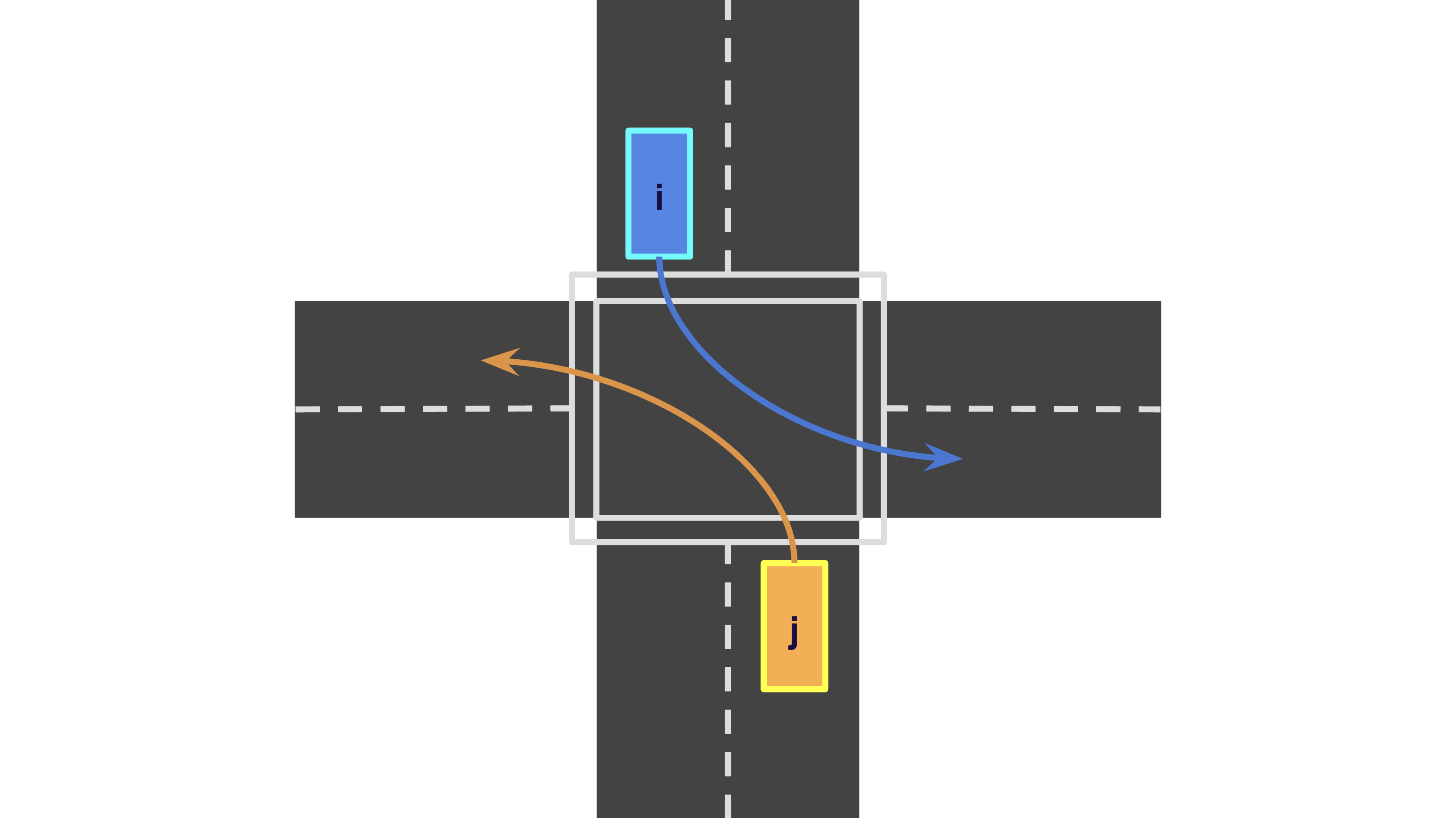}
    \caption{Scenario 1 illustration}
    \label{fig:scen1}
\end{figure}

In this scenario, two vehicles are stopped at an intersection on opposing sides going in opposite directions. When the traffic light turns green, they both turn left, as seen in Figure \ref{fig:scen1}. 

As depicted in Figure \ref{fig:scen1}, a collision between vehicles is not expected although vehicles might come close to each other while turning at the intersection. This corresponds to a scenario where two vehicles safe\edit{l}y negotiate turning at an intersection with minimal safety risk.

\begin{figure} [ht!]
    \centering
    \includegraphics[width=0.49\linewidth]{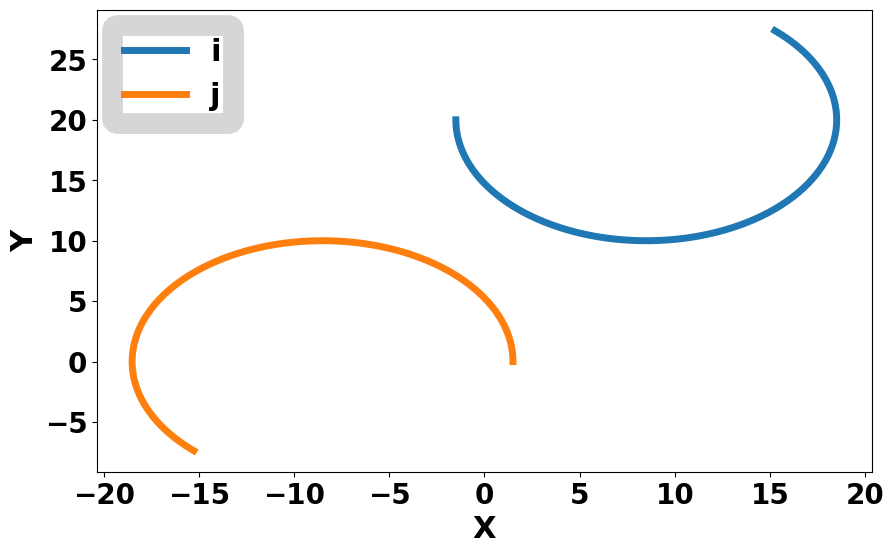}
    \includegraphics[width=0.49\linewidth]{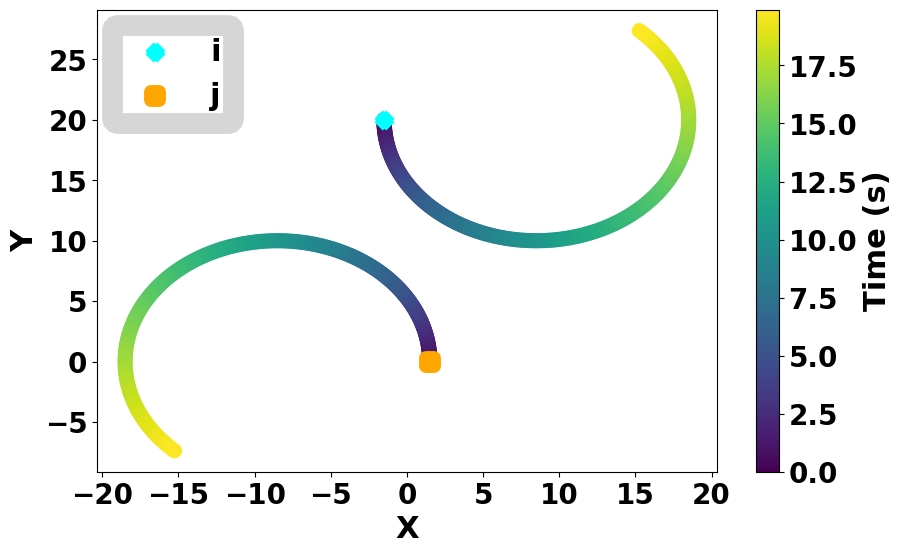}
    \caption{Estimated Trajectories of Vehicles with Initial conditions (\edit{Position: $\+{p}_{\circ,i} = \elements{-1.5,20}$, $\+{p}_{\circ,j} = \elements{1.5,0}$; Velocity: $\+{v}_{\circ,i} = \elements{0,-1}$, $\+{v}_{\circ,j} = \elements{0,1}$; Acceleration: $\+{a}_{\circ,i} = \elements{0.1,-0.1}$, $\+{a}_{\circ,j} = \elements{-0.1,0.1}$}) in Computation of 2D-TTC for Scenario 1. \edit{The l}eft side of the figure shows the overall trajectories of each vehicle\edit{, estimated as circular trajectories for both vehicles at time $t_\circ$ under} the second-order TTC model. \edit{The right side of the figure shows the progression of the vehicles' location from time $t_\circ$ on the estimated trajectories under the second-order TTC scheme through the color gradient.} The cyan mark and orange mark represent the initial position of vehicle $i$ and vehicle $j$ respectively.}
    \label{fig:scen1_2dttc}
\end{figure}

With the initial condition of two vehicles for the scenario, we now estimate trajectories of vehicles, find $\textit{d}_{ij}-\phi$ \edit{along} estimated trajectories, and compute the TTC under both first-order TTC scheme and second-order TTC scheme. The estimated trajectories of vehicles \edit{at time $t_\circ$} in this scenario for the computation of 2D-TTC are shown in Figure \ref{fig:scen1_2dttc}.

Based on the estimated trajectories in Figure \ref{fig:scen1_2dttc}, the second-order TTC scheme is capable of predicting the turning movement of both vehicles at the intersection for computing the TTC and identify\edit{ing} that a collision is not expected.
\begin{figure} [ht!]
    \centering
    \includegraphics[width=0.49\linewidth]{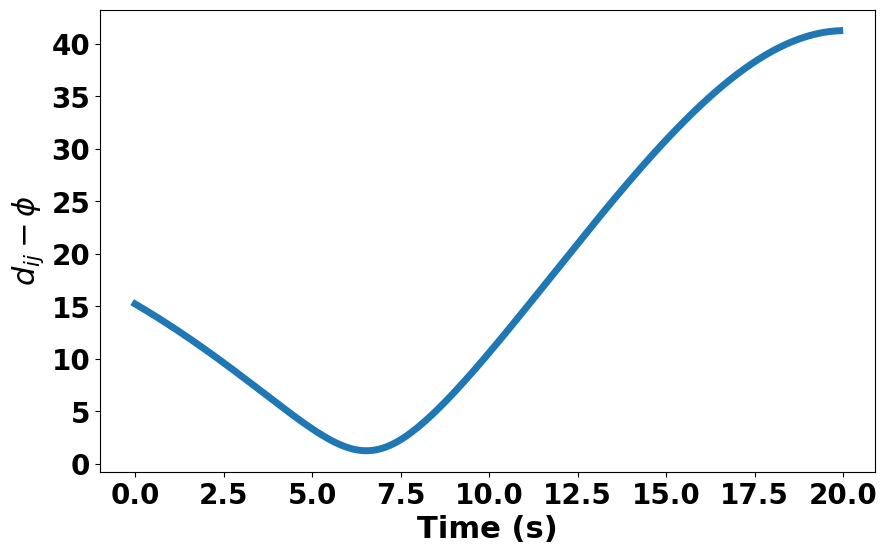}
    \caption{Predicted $\textit{d}_{ij}-\phi$ based on the initial conditions as time progresses in the computation of 2D-TTC \edit{at time $t_\circ$} in  Scenario 1\edit{, estimating no collision.}}
    \label{fig:scen1_2dttc_d}
\end{figure}

Figure \ref{fig:scen1_2dttc_d} shows $\textit{d}_{ij}-\phi$ \edit{along} the estimated trajectories of vehicles shown in Figure \ref{fig:scen1_2dttc}. \edit{T}he vehicles \edit{first} approach each other as $\textit{d}_{ij}-\phi$ is decreasing. Then vehicles move further apart as $\textit{d}_{ij}-\phi$ is increasing back. A collision is not expected to happen \edit{with estimated trajectories of vehicles based on initial conditions at $t_\circ$} as $\textit{d}_{ij}-\phi$ is always above 0.

\begin{figure} [ht!]
    \centering
    \includegraphics[width=0.49\linewidth]{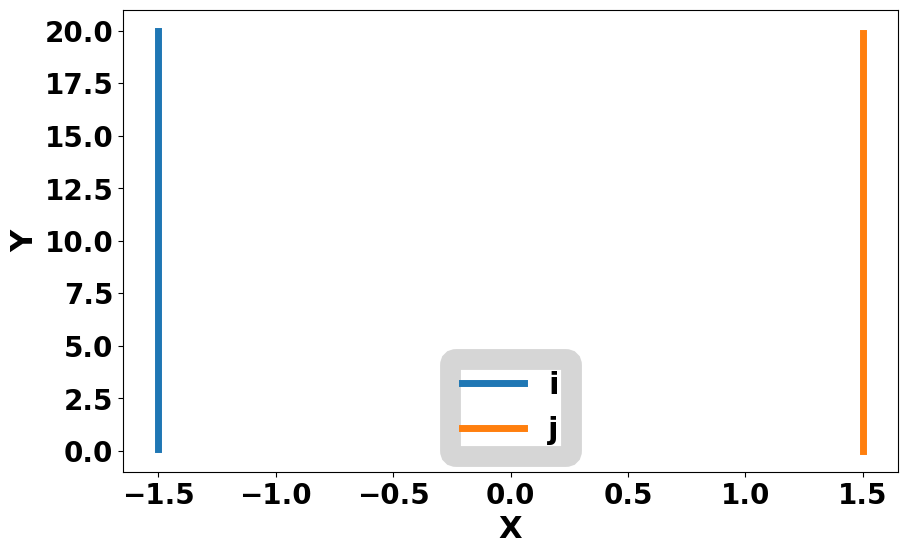}
    \includegraphics[width=0.49\linewidth]{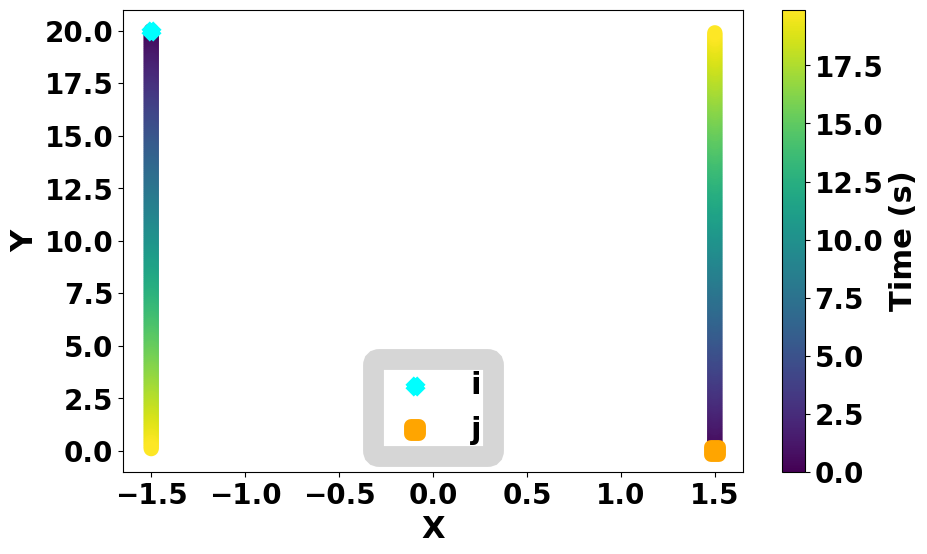}
    \caption{Estimated Trajectories of Vehicles with Initial conditions (same \edit{$\+{p}_{\circ, i}$, $\+{p}_{\circ, j}$, $\+{v}_{\circ, i}$, $\+{v}_{\circ, j}$} as Fig. \ref{fig:scen1_2dttc} except Acceleration: \edit{$\+{a}_{\circ, i} = \elements{0,0}$, $\+{a}_{\circ, j} = \elements{0,0}$}) in Computation of 1D-TTC for Scenario 1. \edit{The l}eft side of the figure shows the overall trajectories of each vehicle\edit{, estimated as straight line trajectories for both vehicles at time $t_\circ$ under} the first-order TTC model. \edit{The right side of the figure shows the progression of the vehicles' location from time $t_\circ$ on the estimated trajectories under the first-order TTC scheme through the color gradient.} The cyan mark and orange mark represent the initial position of vehicle $i$ and vehicle $j$ respectively.} 
    \label{fig:scen1_1dttc}
\end{figure}

Similar to second-order TTC in Figure \ref{fig:scen1_2dttc}, the estimated trajectories of vehicles \edit{based on $\+{p}_{\circ, i}$, $\+{p}_{\circ, j}$, $\+{v}_{\circ, i}$, $\+{v}_{\circ, j}$, $\+{a}_{\circ, i}$, and $\+{a}_{\circ, j}$} in Scenario 1 for the computation of 1D-TTC are shown in Figure \ref{fig:scen1_1dttc}.

Unlike the estimated trajectories of vehicles under the second-order TTC scheme in Figure \ref{fig:scen1_2dttc}, the first-order TTC cannot capture the turning behavior of vehicles as it predicts straight line trajectories of vehicles as shown in Figure \ref{fig:scen1_1dttc}. 

\begin{figure} [ht!]
    \centering
    \includegraphics[width=0.49\linewidth]{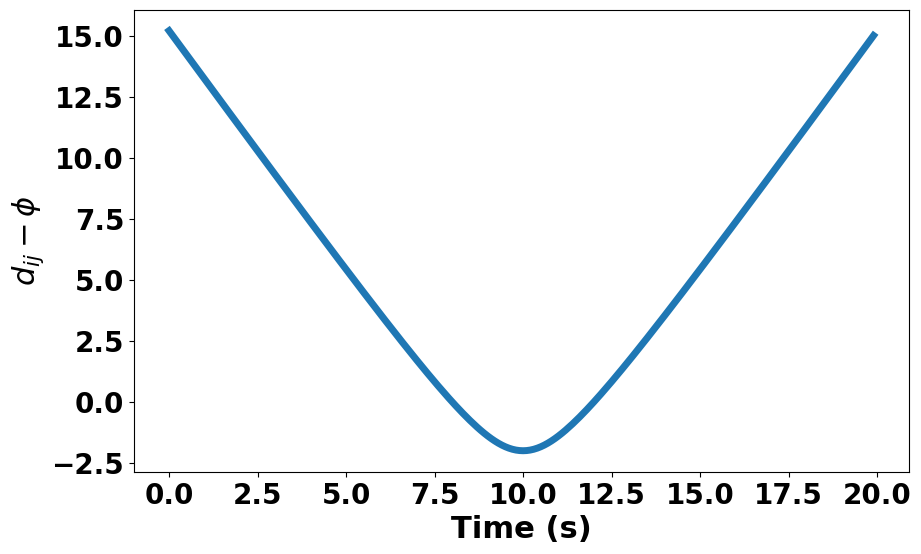}
    \caption{Predicted $\textit{d}_{ij}-\phi$ based on the initial conditions as time progresses in the computation of 1D-TTC \edit{at time $t_\circ$} in  Scenario 1\edit{, estimating a collision.}}
    \label{fig:scen1_1dttc_d}
\end{figure}

Figure \ref{fig:scen1_1dttc_d} shows $\textit{d}_{ij}-\phi$ \edit{along} the estimated trajectories of vehicles \edit{based on initial conditions under 1D-TTC scheme as} shown in Figure \ref{fig:scen1_1dttc}. 


Similar to the $\textit{d}_{ij}-\phi$ under the second-order TTC scheme, the vehicles approach each other initially and further apart afterward as shown as \edit{a} V-shaped pattern in Figure \ref{fig:scen1_1dttc_d}. However, the collision is expected to happen at 8 seconds \edit{from time $t_\circ$} under the first-order scheme as $\textit{d}_{ij}-\phi$ goes below 0. 

In summary, we expect no collision between two vehicles as they are turning at the intersection based on the design of the scenario as shown in Figure \ref{fig:scen1}. While 2D-TTC can capture the turning movement of the vehicles in \edit{estimated trajectories as shown in} Figure \ref{fig:scen1_2dttc}, 1D-TTC cannot and assumes they travel in straight trajectories \edit{as drawn} in Figure \ref{fig:scen1_1dttc}. In addition, \edit{the safety assessment at time $t_\circ$ expects a} collision to happen for 1D-TTC while \edit{not} for 2D-TTC from Figure \ref{fig:scen1_2dttc_d} and \ref{fig:scen1_1dttc_d}. 

\begin{figure} [ht!]
    \centering
  \includegraphics[width=0.5\linewidth]{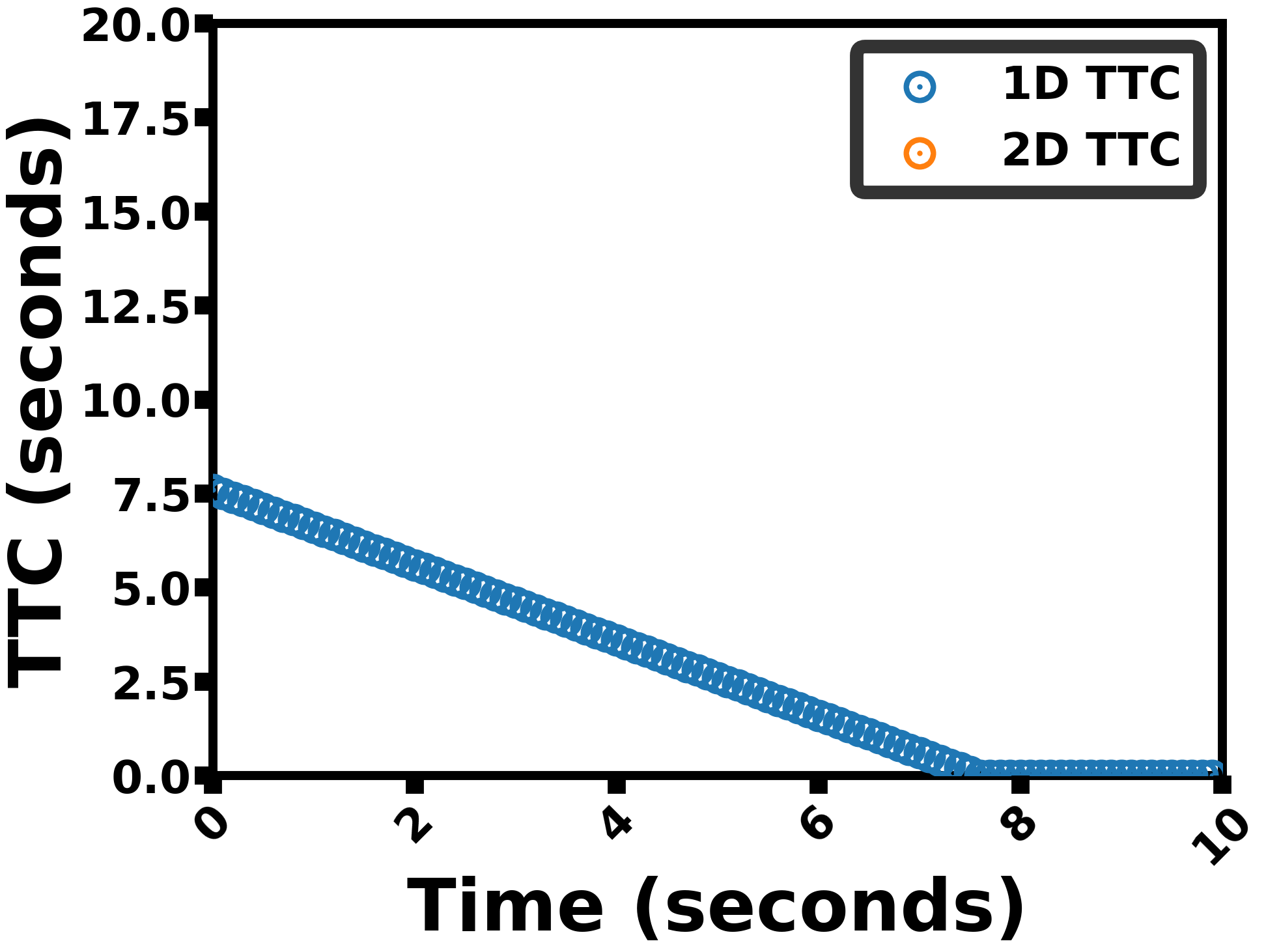}
    \caption{\edit{TTC values} for each timestep in Scenario 1. \edit{On the plot, only 1D-TTC exists as the collision between two vehicles is predicted continuously over time under the first-order TTC model while not under the second-order TTC model.}}
    \label{fig:traj_compare_scenario1}
\end{figure}

We continue computing both the first-order TTC and the second-order TTC \edit{at} each timestep \edit{for numerical simulation over 10 seconds} in Scenario 1. Figure \ref{fig:traj_compare_scenario1} is the plot comparing the first-order TTC and the second-order TTC in Scenario 1.

In Figure \ref{fig:traj_compare_scenario1}, the first-order TTC scheme continuously expects a collision to happen. TTC value in the first-order TTC scheme decreases as simulation time progresses. The second-order TTC scheme, on the other hand, do\edit{es} not expect a collision within an entire simulation time horizon. 

\edit{Scenario 1 is set up such that both vehicles are turning at an intersection without collision.} As a collision between two vehicles is not designed to happen in Scenario 1, the second-order TTC scheme performs better safety assessment than the first-order TTC scheme for Scenario 1. 

\subsection{Scenario 2: One Vehicle Turning Right into a Slow/Unmoving Vehicle}

\begin{figure} [ht!]
    \centering
  \includegraphics[width=0.9\linewidth]{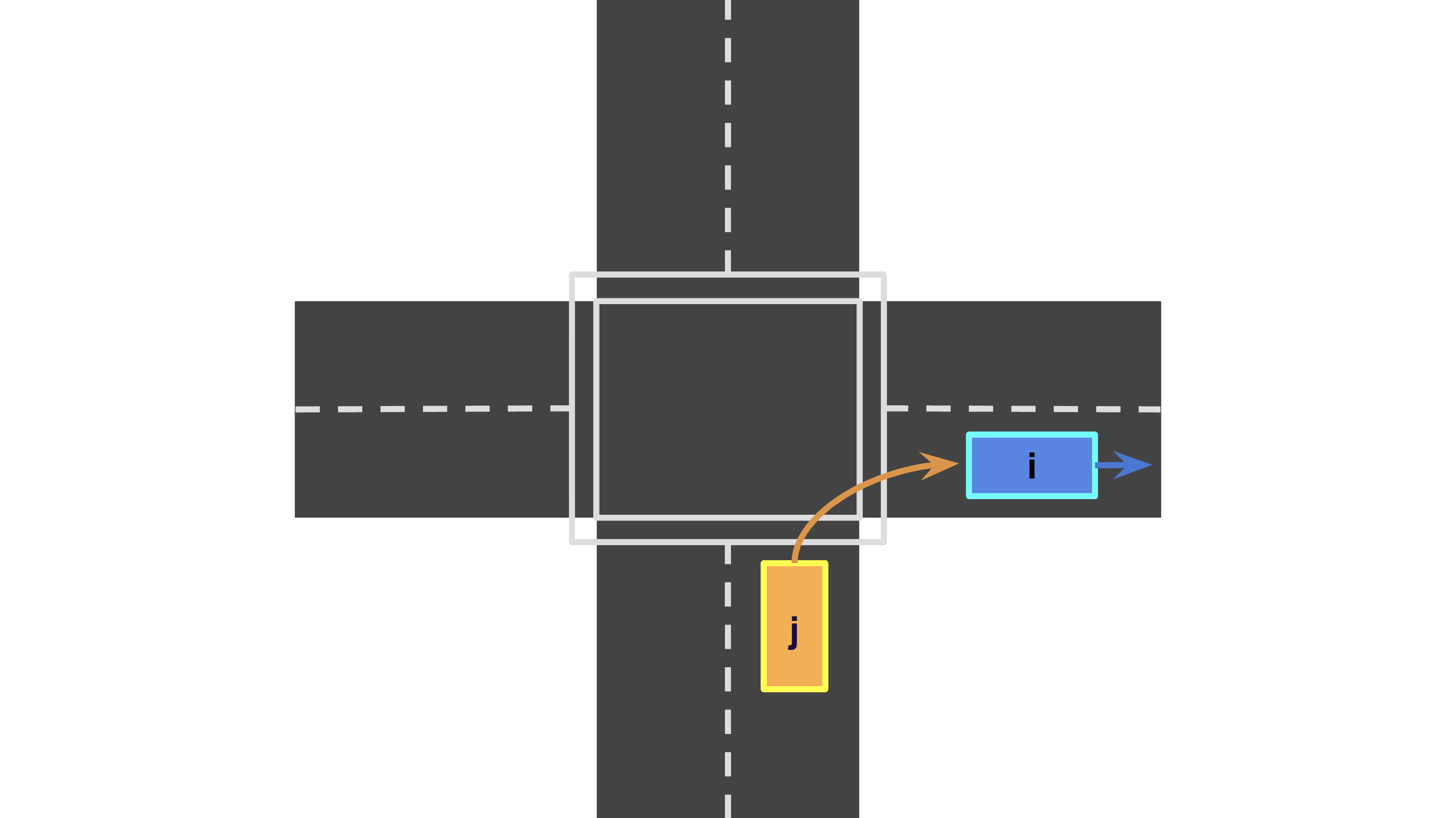}
    \caption{Scenario 2 illustration}
    \label{fig:scen2}
\end{figure}

In this scenario, one vehicle is stopped at an intersection about to turn right into the lane the other vehicle occupies. When the traffic light turns green, the first vehicle turns right into the other vehicle and collides, as seen in Figure \ref{fig:scen2}. 

Based on how scenario 2 is set up, we expect the vehicles to approach close to each other and the TTC value to reflect the scenario as \edit{a} near-crash scenario. 

\begin{figure} [ht!]
    \centering
    \includegraphics[width=0.49\linewidth]{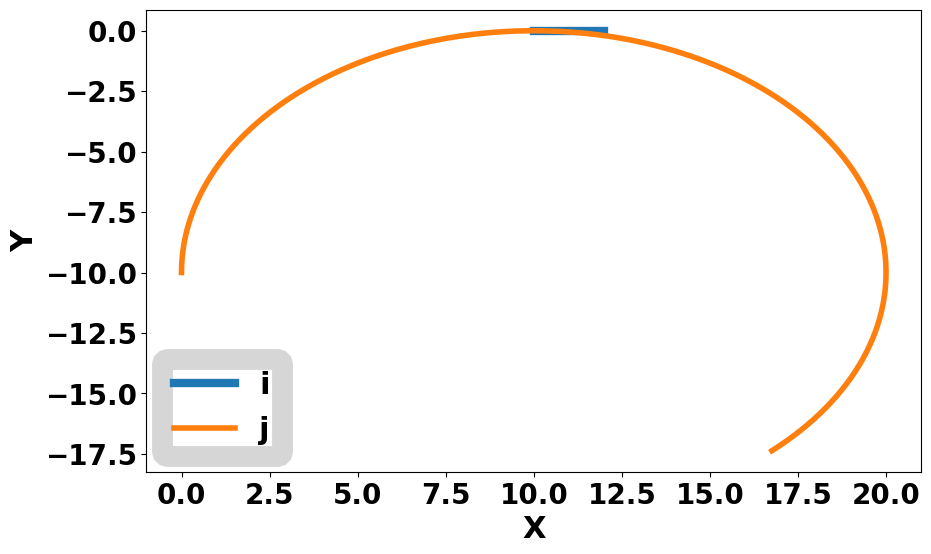}
    \includegraphics[width=0.49\linewidth]{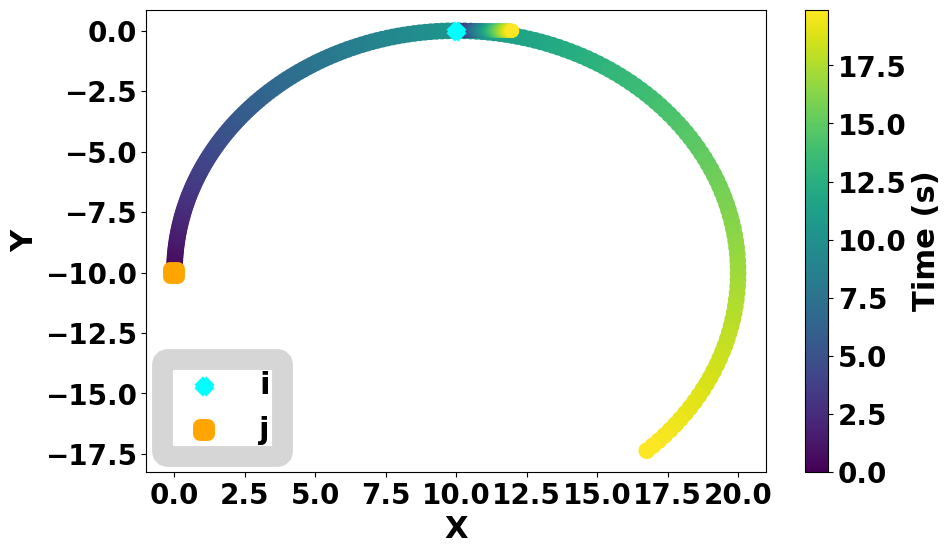}
    \caption{Estimated Trajectories of Vehicles with Initial conditions (\edit{Position: $\+{p}_{\circ,i} = \elements{10,0}$, $\+{p}_{\circ,j} = \elements{0,-10}$; Velocity: $\+{v}_{\circ,i} = \elements{0.1,0}$, $\+{v}_{\circ,j} = \elements{0,1}$; Acceleration: $\+{a}_{\circ,i} = \elements{0,0}$, $\+{a}_{\circ,j} = \elements{0.1,-0.1}$}) in Computation of 2D-TTC for Scenario 2. \edit{The l}eft side of the figure shows the overall trajectories of each vehicle\edit{, estimated as a straight line trajectory for vehicles $i$ and a circular trajectory for vehicle $j$ at time $t_\circ$ under} the second-order TTC model. \edit{The right side of the figure shows the progression of the vehicles' location from time $t_\circ$ on the estimated trajectories under the second-order TTC scheme through the color gradient.} The cyan mark and orange mark represent the initial position of vehicle $i$ and vehicle $j$ respectively.}
    \label{fig:scen2_2dttc}
\end{figure}

We first investigate the performance of the second-order TTC with the initial condition of both vehicles. The estimated trajectories of vehicles \edit{based on $\+{p}_{\circ, i}$, $\+{p}_{\circ, j}$, $\+{v}_{\circ, i}$, $\+{v}_{\circ, j}$, $\+{a}_{\circ, i}$, and $\+{a}_{\circ, j}$} in scenario 2 for the computation of 2D-TTC are shown in Figure \ref{fig:scen2_2dttc}.

From Figure \ref{fig:scen2_2dttc}, the second-order TTC well reflect\edit{s} the turning of the vehicle $j$. In addition, the right side of Figure \ref{fig:scen2_2dttc} visually shows that the collision of two vehicles to happen \edit{along $\+{p}_i(t)$ and $\+{p}_j(t)$ determined at time $t_\circ$} as the second-order TTC scheme estimates both vehicles to be very close to each other based on the color \edit{gradients on} estimated trajectories. 

\begin{figure} [ht!]
    \centering
    \includegraphics[width=0.49\linewidth]{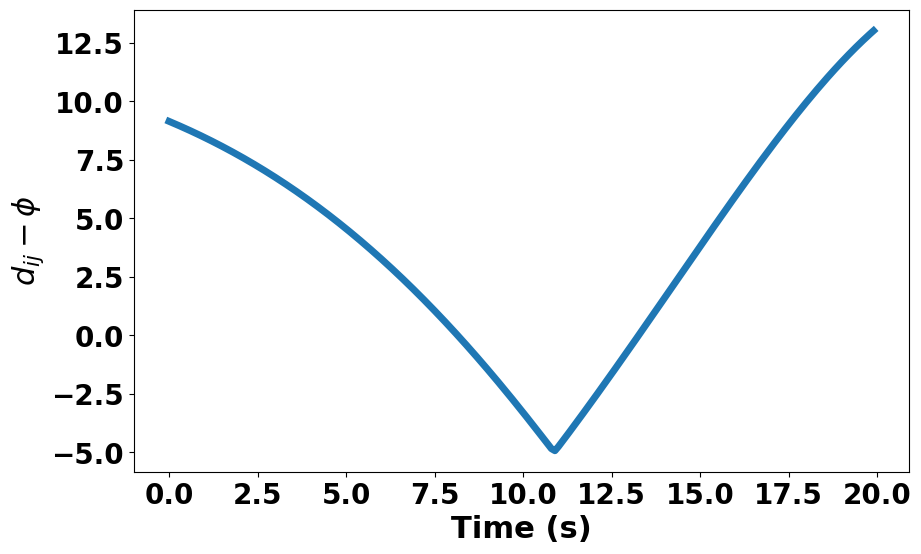}
    \caption{Predicted $\textit{d}_{ij}-\phi$ based on the initial conditions as time progresses in the computation of 2D-TTC \edit{at time $t_\circ$} in  Scenario 2\edit{, estimating a collision.}}
    \label{fig:scen2_2dttc_d}
\end{figure}

Based on the estimated trajectories of two vehicles \edit{in Figure \ref{fig:scen2_2dttc}}, Figure \ref{fig:scen2_2dttc_d} shows $\textit{d}_{ij}-\phi$ \edit{along the estimated trajectories based on initial conditions}. 

\begin{figure} [ht!]
    \centering
    \includegraphics[width=0.49\linewidth]{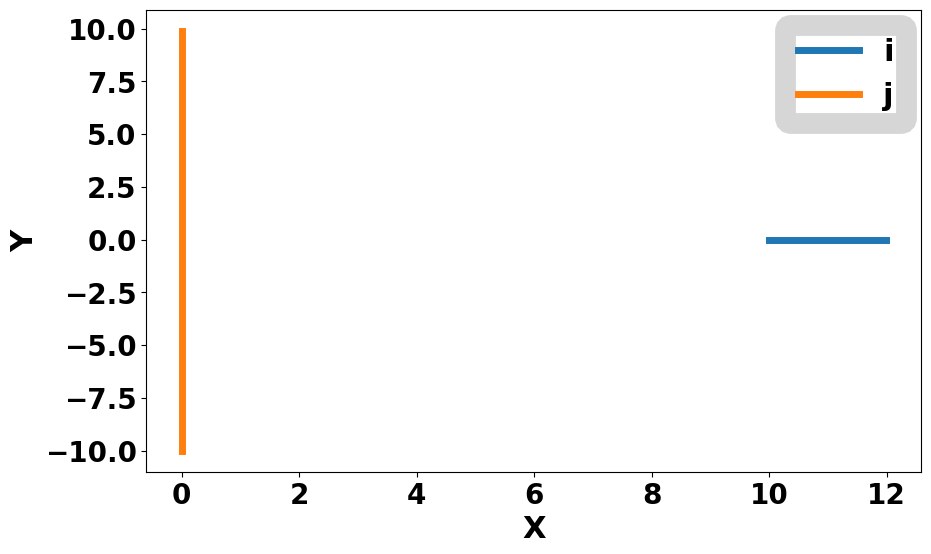}
    \includegraphics[width=0.49\linewidth]{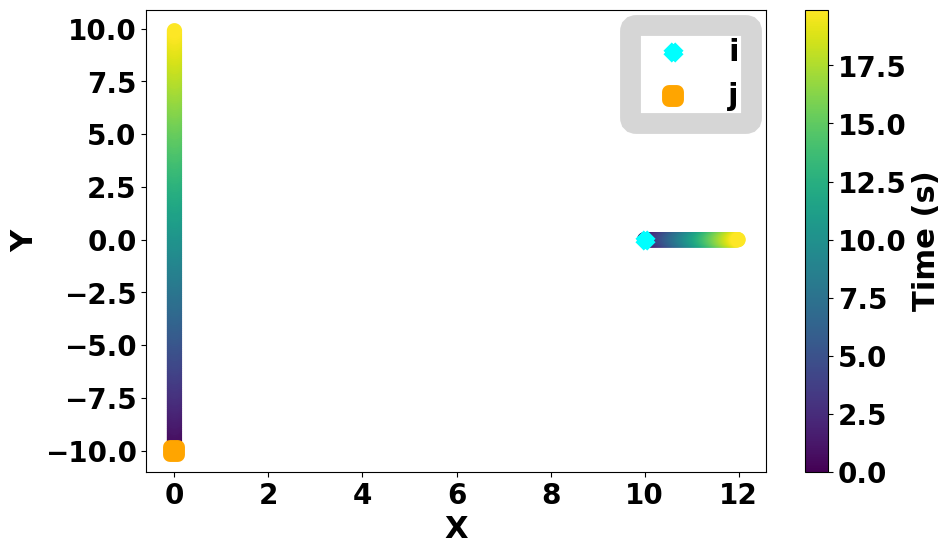}
    \caption{Estimated Trajectories of Vehicles with Initial conditions (same \edit{$\+{p}_{\circ, i}$, $\+{p}_{\circ, j}$, $\+{v}_{\circ, i}$, $\+{v}_{\circ, j}$} as Fig. \ref{fig:scen2_2dttc} except Acceleration: \edit{$\+{a}_{\circ, i} = \elements{0,0}$, $\+{a}_{\circ, j} = \elements{0,0}$}) in Computation of 1D-TTC for Scenario 2. \edit{The l}eft side of the figure shows the overall trajectories of each vehicle\edit{, estimated as straight line trajectories for both vehicles at time $t_\circ$ under} the first-order TTC model. \edit{The right side of the figure shows the progression of the vehicles' location from time $t_\circ$ on the estimated trajectories under the first-order TTC scheme through the color gradient.} The cyan mark and orange mark represent the initial position of vehicle $i$ and vehicle $j$ respectively.}
    \label{fig:scen2_1dttc}
\end{figure}

\begin{figure} [ht!]
    \centering
    \includegraphics[width=0.49\linewidth]{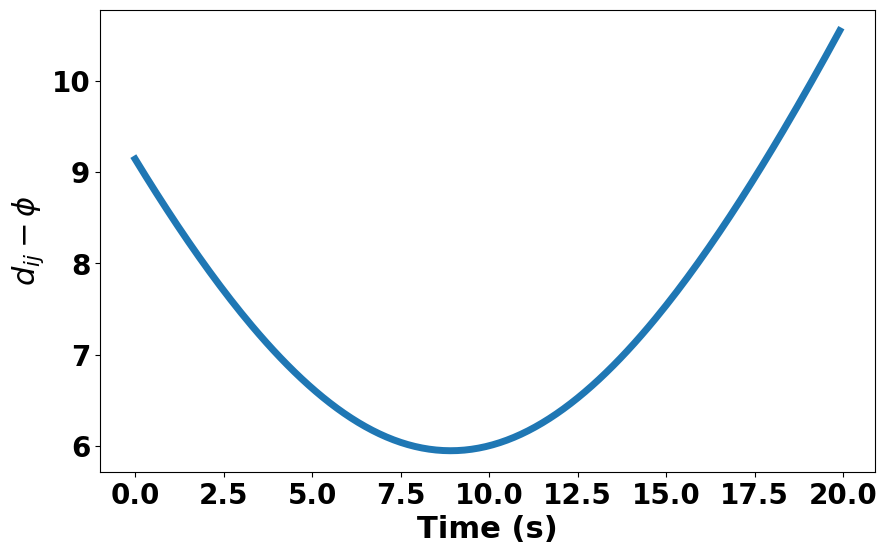}
    \caption{Predicted $\textit{d}_{ij}-\phi$ based on the initial conditions as time progresses in the computation of 1D-TTC \edit{at time $t_\circ$} in  Scenario 2\edit{, estimating no collision.}}
    \label{fig:scen2_1dttc_d}
\end{figure}

From the pattern of $\textit{d}_{ij}-\phi$ in Figure \ref{fig:scen2_2dttc_d}, not only vehicles approach close to each other but also \edit{vehicles expect to collide} 8.15 seconds \edit{after time $t_\circ$} under the second-order scheme. \\
We now move to the first-order TTC for scenario 2. The estimated trajectories of vehicles in scenario 2 for the computation of 1D-TTC are shown in Figure \ref{fig:scen2_1dttc}.\\
Unlike \edit{the} estimated trajectories of vehicle $j$ in Figure \ref{fig:scen2_2dttc}, vehicle $j$ is estimated to move in \edit{a} straight line trajectory in Figure \ref{fig:scen2_1dttc}. Therefore, estimated trajectories are not overlapped under the first-order TTC scheme. \\
Figure \ref{fig:scen2_1dttc_d} shows $\textit{d}_{ij}-\phi$ \edit{along} the estimated trajectories of vehicles shown in Figure \ref{fig:scen2_1dttc}. 
\edit{The pattern of} $\textit{d}_{ij}-\phi$ in Figure \ref{fig:scen2_1dttc_d} supports the observation made in Figure \ref{fig:scen2_1dttc} as $\textit{d}_{ij}-\phi$ never goes below 0. Therefore, the collision is estimated not to happen \edit{based on initial positions, velocities, and accelerations of vehicles} under the first-order TTC scheme.

\begin{figure} [ht!]
    \centering
  \includegraphics[width=0.5\linewidth]{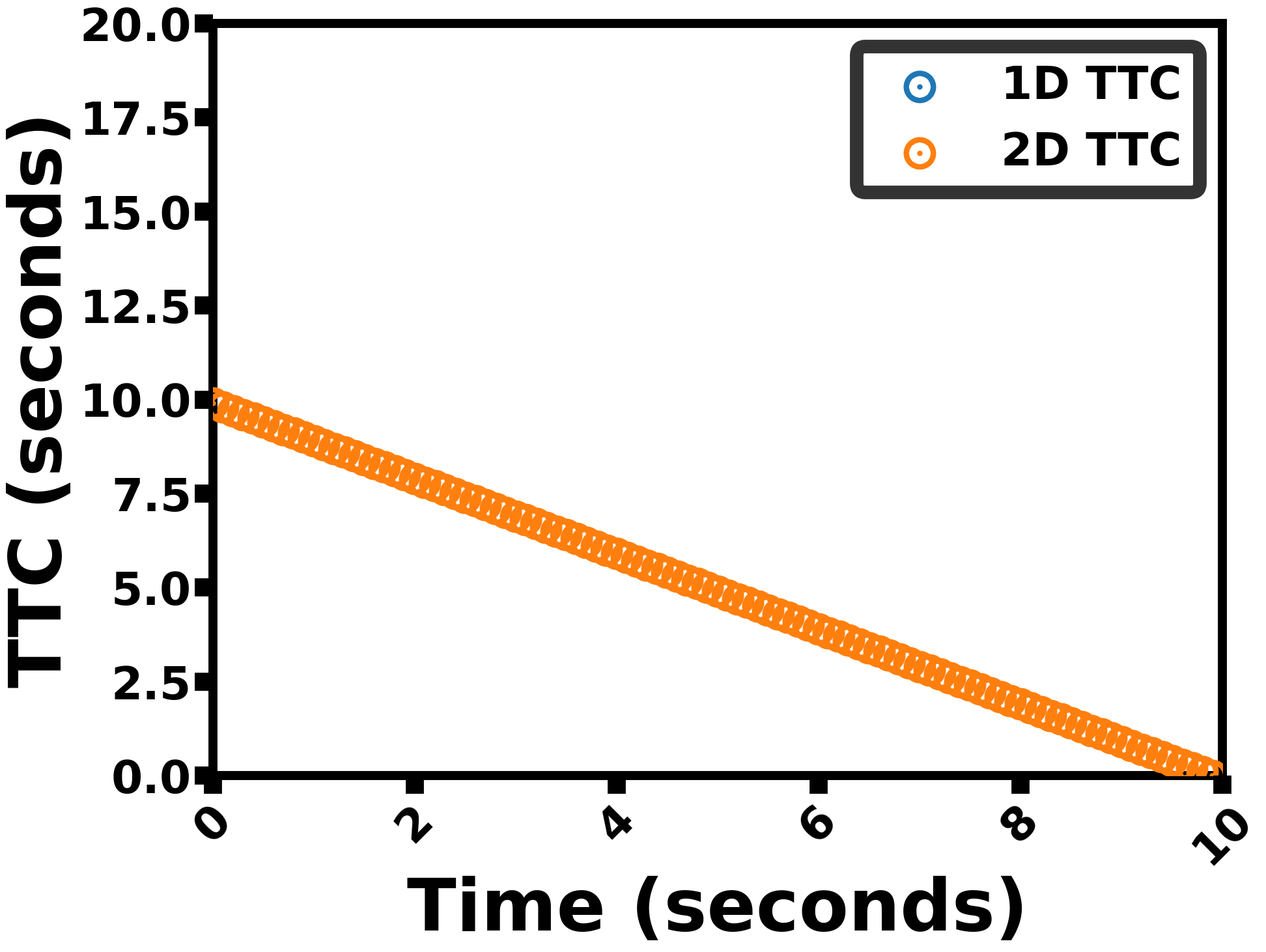}
    \caption{\edit{TTC values} for each timestep in Scenario 2. \edit{On the plot, only 2D-TTC exists as the collision between two vehicles is predicted continuously over time under the second-order TTC model while not under the first-order TTC model.}}
    \label{fig:traj_compare_scenario2}
\end{figure}

For scenario 2, we expect at least \edit{a} near-collision of the turning vehicle into the unmoving vehicle in the 2D-TTC scenario as shown in Figure \ref{fig:scen2}. From Figure \ref{fig:scen2_2dttc} and Figure \ref{fig:scen2_2dttc_d}, collision is estimated to happen under 2D-TTC. As \edit{the} 1D-TTC scheme cannot predict the turning of vehicle $j$ as shown in Figure \ref{fig:scen2_1dttc}, $\textit{d}_{ij}-\phi$ does not approach 0 in Figure \ref{fig:scen2_1dttc_d}. Therefore, \edit{the} 1D-TTC scheme cannot capture \edit{the} near-crash between \edit{the} two vehicles. \\
We now compare the first-order TTC and the second-order TTC for Scenario 2 at each timestep in Figure \ref{fig:traj_compare_scenario2}. 

\noindent In Scenario 2, the second-order TTC scheme keeps expecting a collision between vehicle $i$ and $j$ to happen. In addition, the value of TTC decreases as the simulation of Scenario 2 proceeds. In contrast, the first-order TTC scheme do\edit{es} not estimate the collision throughout Scenario 2. \\
In conclusion, 2D-TTC better reflect\edit{s} Scenario 2 as it correctly estimates a collision between two vehicles unlike 1D-TTC scheme \edit{based on Scenario 2 setup, where a right-turning vehicle is on a collision course with a slowly moving vehicle.} 

\subsection{Scenario 3: Both Vehicles Turning from Perpendicular Lanes}

\begin{figure} [ht!]
    \centering
  \includegraphics[width=0.9\linewidth]{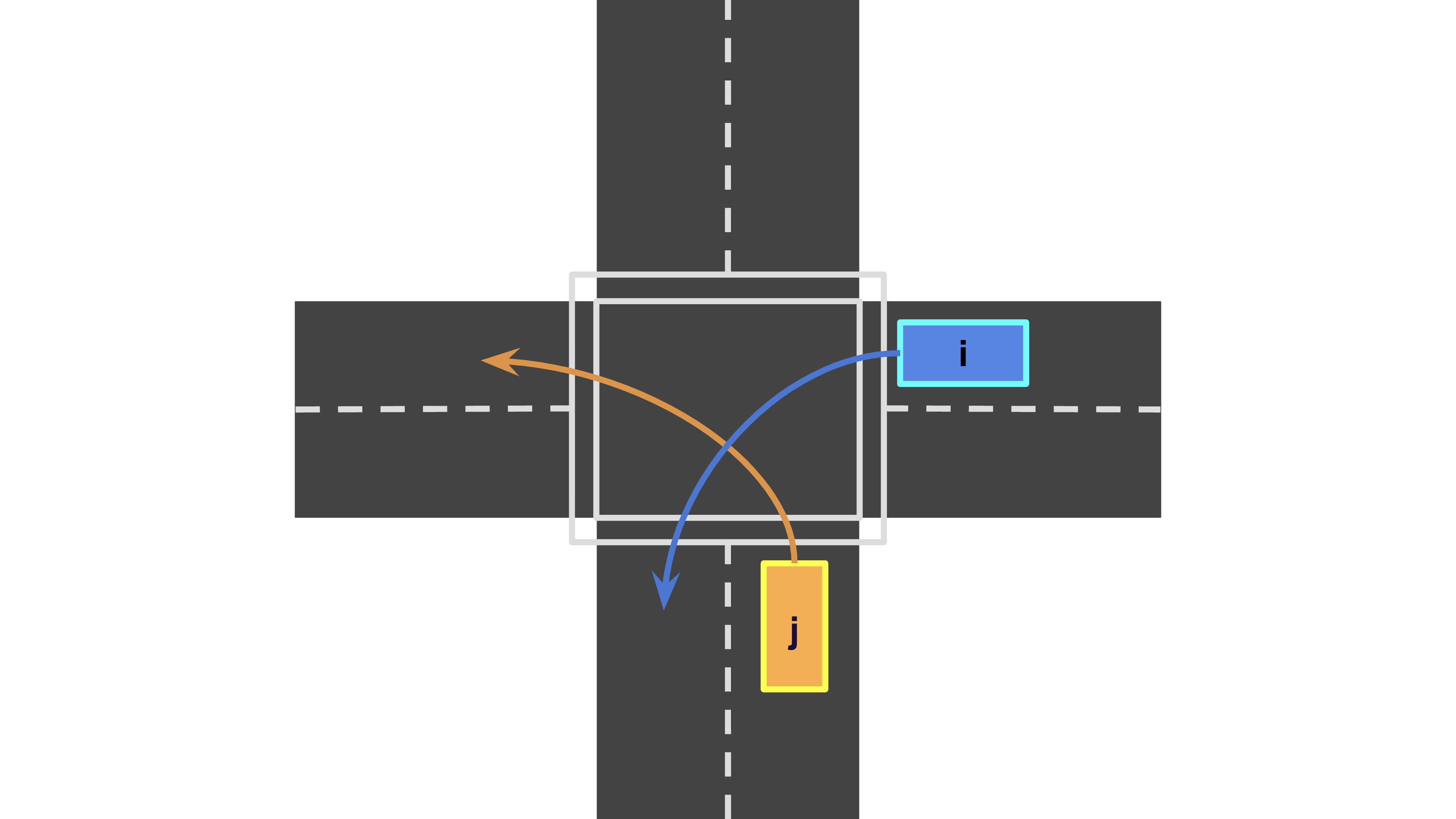}
    \caption{Scenario 3 illustration}
    \label{fig:scen3}
\end{figure}

In this scenario, both vehicles are about to turn left, but one vehicle is turning into the other vehicle's lane. When the traffic light turns green, both start turning but they go past each other without collision, as seen in Figure \ref{fig:scen3}. 

\begin{figure} [ht!]
    \centering
    \includegraphics[width=0.49\linewidth]{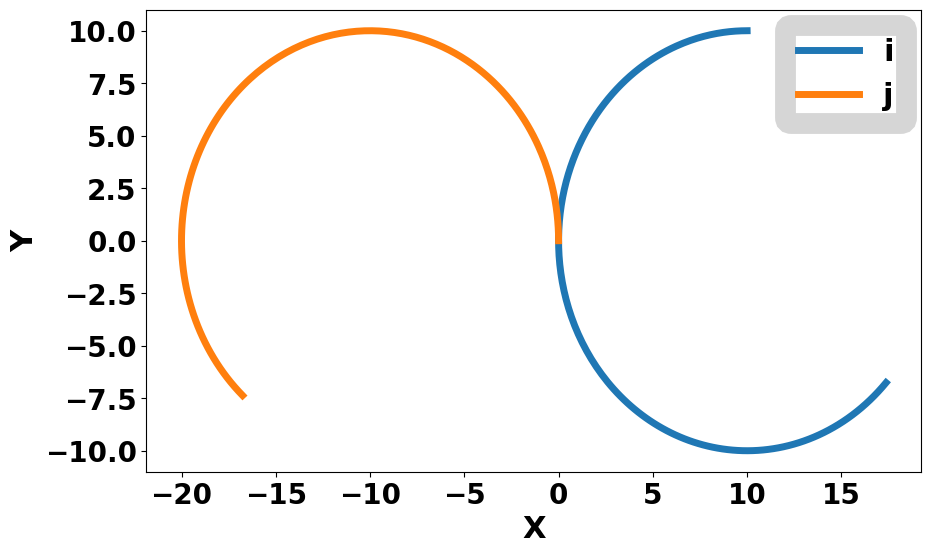}
    \includegraphics[width=0.49\linewidth]{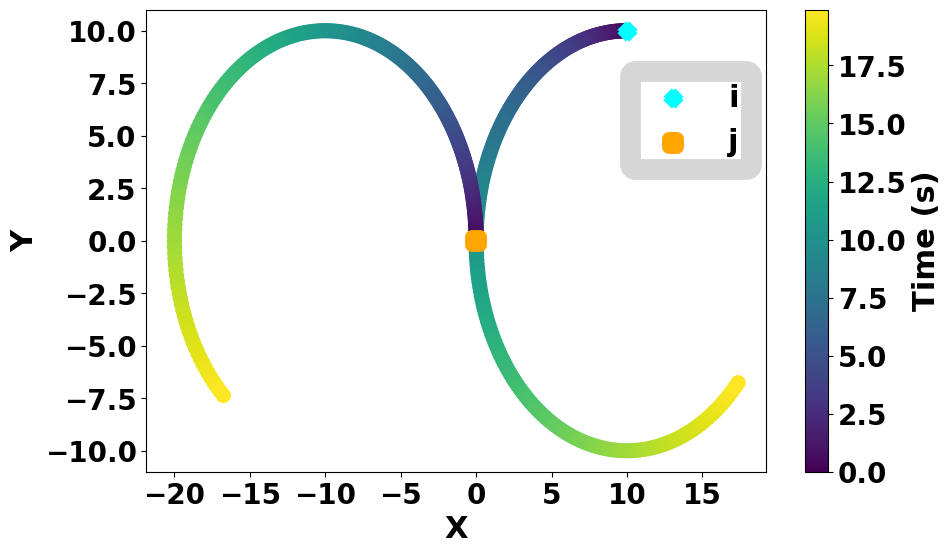}
    \caption{Estimated Trajectories of Vehicles with Initial conditions (\edit{Position: $\+{p}_{\circ,i} = \elements{10,10}$, $\+{p}_{\circ,j} = \elements{0,0}$; Velocity: $\+{v}_{\circ,i} = \elements{-1,0}$, $\+{v}_{\circ,j} = \elements{0,1}$; Acceleration: $\+{a}_{\circ,i} = \elements{-0.1,-0.1}$, $\+{a}_{\circ,j} = \elements{-0.1,0.1}$}) in Computation of 2D-TTC for Scenario 3. \edit{The l}eft side of the figure shows the overall trajectories of each vehicle\edit{, estimated as circular trajectories for both vehicles at time $t_\circ$ under} the second-order TTC model. \edit{The right side of the figure shows the progression of the vehicles' location from time $t_\circ$ on the estimated trajectories under the second-order TTC scheme through the color gradient.} The cyan mark and orange mark represent the initial position of vehicle $i$ and vehicle $j$ respectively.}
    \label{fig:scen3_2dttc}
\end{figure}

\begin{figure} [ht!]
    \centering
    \includegraphics[width=0.49\linewidth]{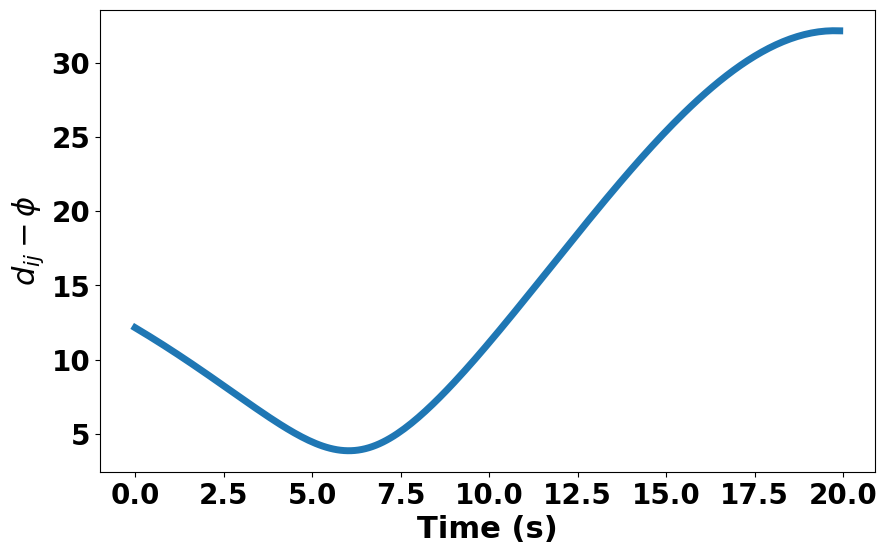}
    \caption{Predicted $\textit{d}_{ij}-\phi$ based on the initial conditions as time progresses in the computation of 2D-TTC \edit{at time $t_\circ$} in  Scenario 3\edit{, estimating no collision.}}
    \label{fig:scen3_2dttc_d}
\end{figure}

Although a collision seems likely to occur as the trajectories of vehicles overlap in Figure \ref{fig:scen3}, the scenario and initial conditions are designed so that vehicles make turn\edit{s} and approach close to each other but without an actual collision. Therefore, we want \edit{the} TTC value to reflect near-crash but no collision. \\
The estimated trajectories of vehicles \edit{from $t_\circ$} in scenario 3 for the computation of 2D-TTC are shown in Figure \ref{fig:scen3_2dttc}.
The estimated trajectories of vehicles under the second-order TTC are well-estimating the actual trajectories of the scenario. The estimated trajectories show that both vehicles are estimated to turn left at the intersection. The color on the estimated trajectories \edit{on} the right side of Figure \ref{fig:scen3_2dttc} depicts at least the near collision case where two vehicles are approaching close to each other. \\
Through Figure \ref{fig:scen3_2dttc_d}, we investigate further by looking at $\textit{d}_{ij}-\phi$ \edit{along} the estimated trajectories of vehicles shown in Figure \ref{fig:scen3_2dttc}. 

\begin{figure} [ht!]
    \centering
    \includegraphics[width=0.49\linewidth]{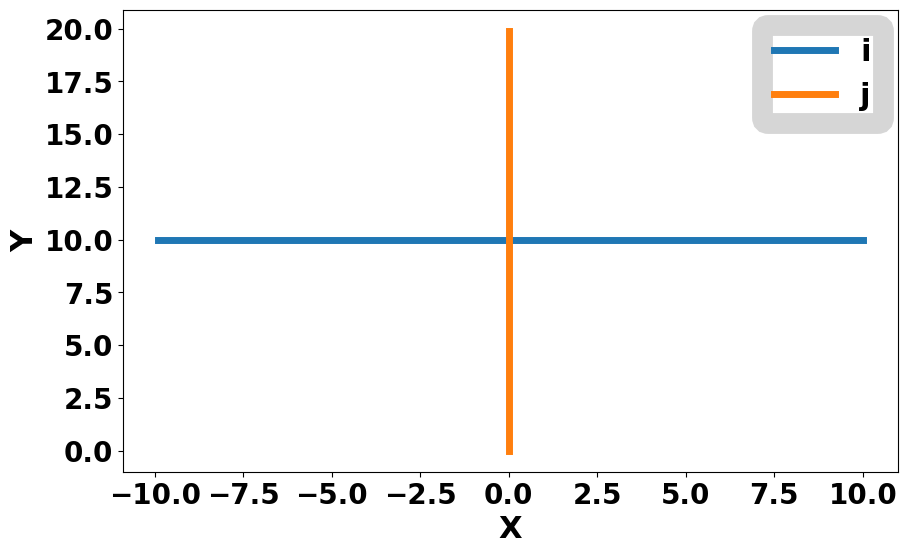}
    \includegraphics[width=0.49\linewidth]{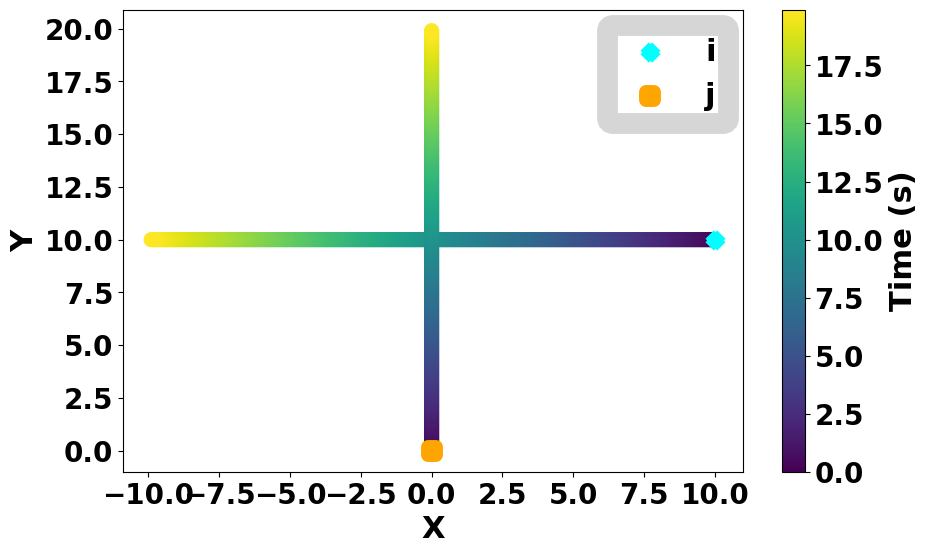}
    \caption{Estimated Trajectories of Vehicles with Initial conditions (same \edit{$\+{p}_{\circ, i}$, $\+{p}_{\circ, j}$, $\+{v}_{\circ, i}$, $\+{v}_{\circ, j}$} as Fig. \ref{fig:scen3_2dttc} except Acceleration: \edit{$\+{a}_{\circ, i} = \elements{0,0}$, $\+{a}_{\circ, j} = \elements{0,0}$}) in Computation of 1D-TTC for Scenario 3. \edit{The l}eft side of the figure shows the overall trajectories of each vehicle\edit{, estimated as straight line trajectories for both vehicles at time $t_\circ$ under} the first-order TTC model. \edit{The right side of the figure shows the progression of the vehicles' location from time $t_\circ$ on the estimated trajectories under the first-order TTC scheme through the color gradient.} The cyan mark and orange mark represent the initial position of vehicle $i$ and vehicle $j$ respectively.}
    \label{fig:scen3_1dttc}
\end{figure}

\noindent From Figure \ref{fig:scen3_2dttc_d}, $\textit{d}_{ij}-\phi$ approaches close to 0 but does not fall below 0. In other words, the second-order TTC can capture the near-collision of two vehicles but no actual collision, which well aligns with the design of the scenario. \\
Figure \ref{fig:scen3_1dttc} shows the estimated trajectories of vehicles \edit{based on initial conditions} in scenario 3 for the computation of 1D-TTC. 

\noindent Due to the limitation of the first-order TTC, first-order TTC cannot estimate the turning of vehicles in the scenario and incorrectly produce \edit{straight-line} trajectories of the vehicles. 

\begin{figure} [ht!]
    \centering
    \includegraphics[width=0.49\linewidth]{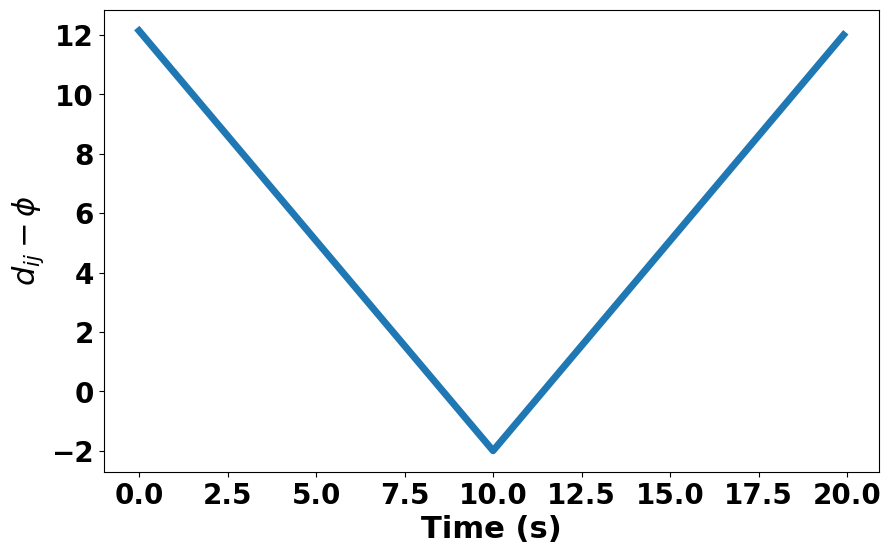}
    \caption{Predicted $\textit{d}_{ij}-\phi$ based on the initial conditions as time progresses in the computation of 1D-TTC \edit{at time $t_\circ$} in  Scenario 3\edit{, estimating a collision.}}
    \label{fig:scen3_1dttc_d}
\end{figure}

Based on the estimated trajectories of vehicles shown in Figure \ref{fig:scen3_1dttc}, Figure \ref{fig:scen3_1dttc_d} shows $\textit{d}_{ij}-\phi$ \edit{based on $\+{p}_{i}$ and $\+{p}_{j}$ for $t \geq t_\circ$}. 

Similar to $\textit{d}_{ij}-\phi$ for 2D-TTC, two vehicles approach close to each other. However, the collision between two vehicles is expected at 6.46 seconds \edit{after moving along estimated trajectories} under the first-order TTC scheme. 

For each timestep in the simulation of Scenario 3, the first-order TTC and the second-order TTC are continuously computed and compared. 
\begin{figure} [ht!]
    \centering
  \includegraphics[width=0.5\linewidth]{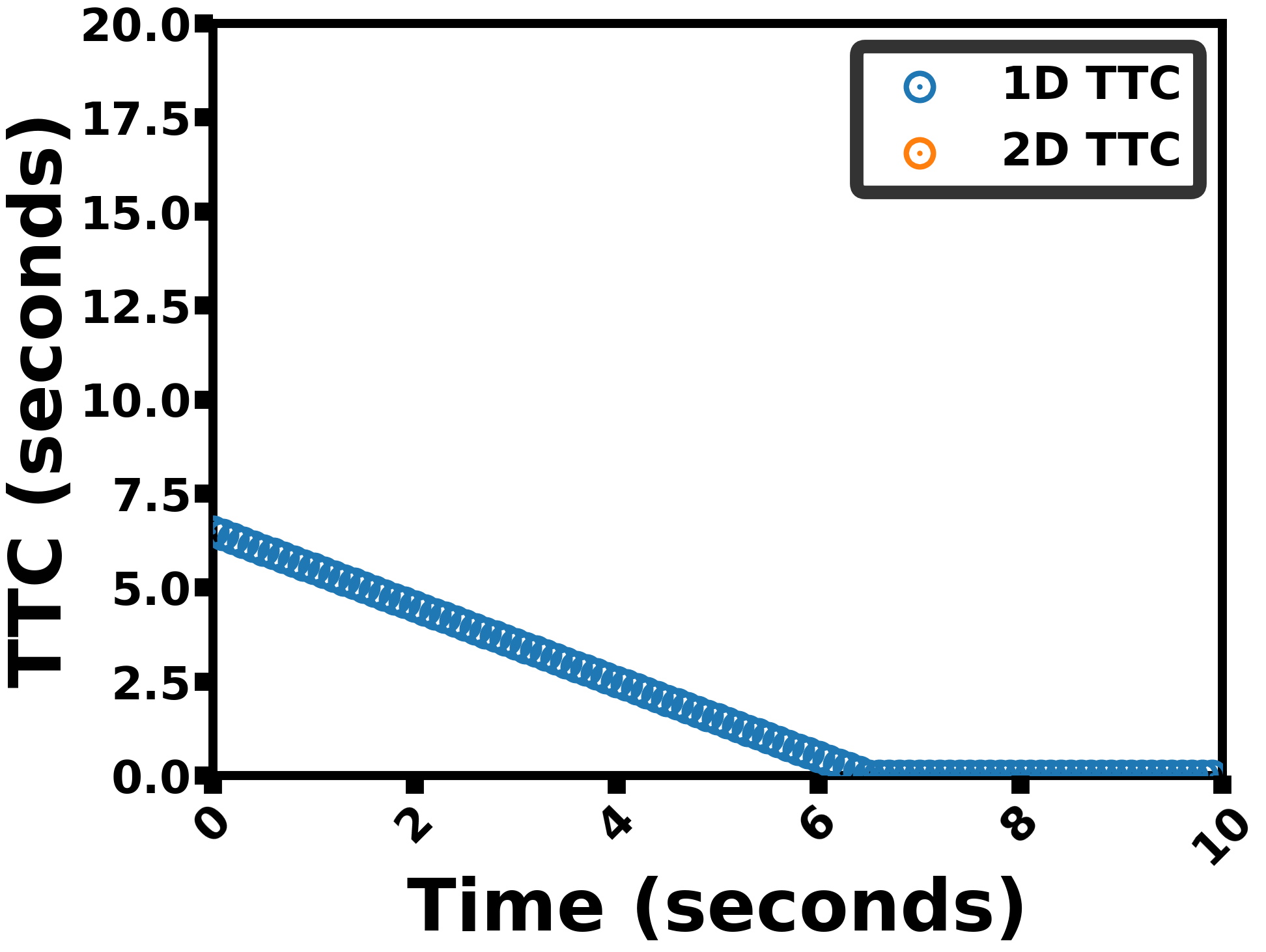}
    \caption{\edit{TTC values} for each timestep in Scenario 3. \edit{On the plot, only 1D-TTC exists as the collision between two vehicles is predicted continuously over time under the first-order TTC model while not under the second-order TTC model.}}
    \label{fig:traj_compare_scenario3}
\end{figure}

In Figure \ref{fig:traj_compare_scenario3}, the first-order TTC scheme and the second-order TTC scheme estimate differently. While the collision is not expected in the second-order TTC scheme for Scenario 3, the collision is continuously estimated for Scenario 3 under the first-order TTC. Based on the design of Scenario 3 \edit{where two vehicles are expected to collide when turning at the intersection}, the second-order TTC performs better in estimating the collision for Scenario 3. 



\subsection{Scenario 4: Vehicle Turning Left into Vehicle \edit{G}oing Straight}

\begin{figure} [ht!]
    \centering
  \includegraphics[width=0.9\linewidth]{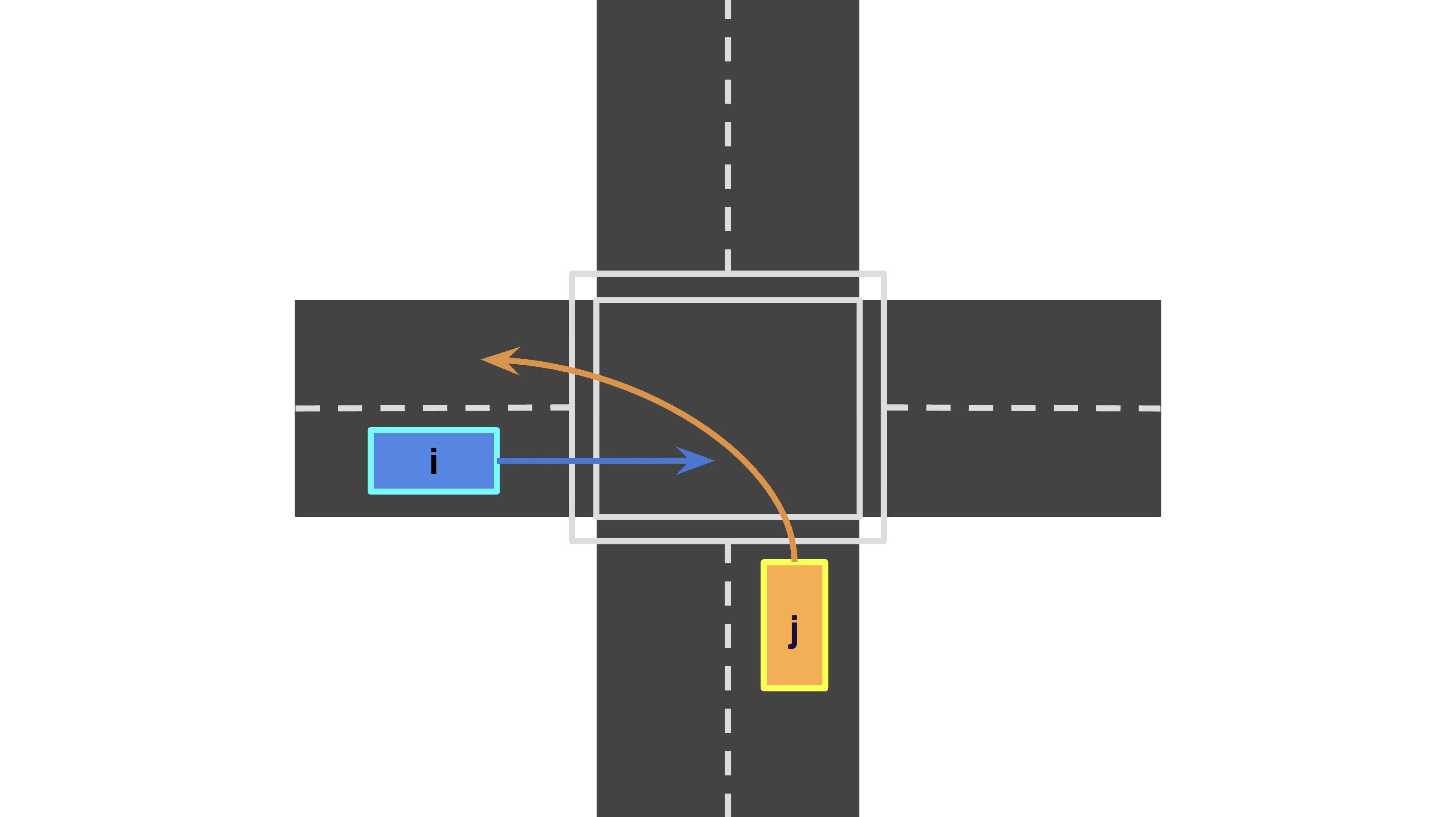}
    \caption{Scenario 4 illustration}
    \label{fig:scen4}
\end{figure}

In this scenario, one vehicle is going straight and another vehicle is about to turn left into the first vehicle's lane, as seen in Figure \ref{fig:scen4}. 

As shown in Figure \ref{fig:scen4}, the scenario is designed so that the collision between two vehicles happens. Therefore, the ideal TTC value needs to not only estimate the turning of vehicle $j$ but also predict the collision between two vehicles. 

\begin{figure} [ht!]
    \centering
    \includegraphics[width=0.49\linewidth]{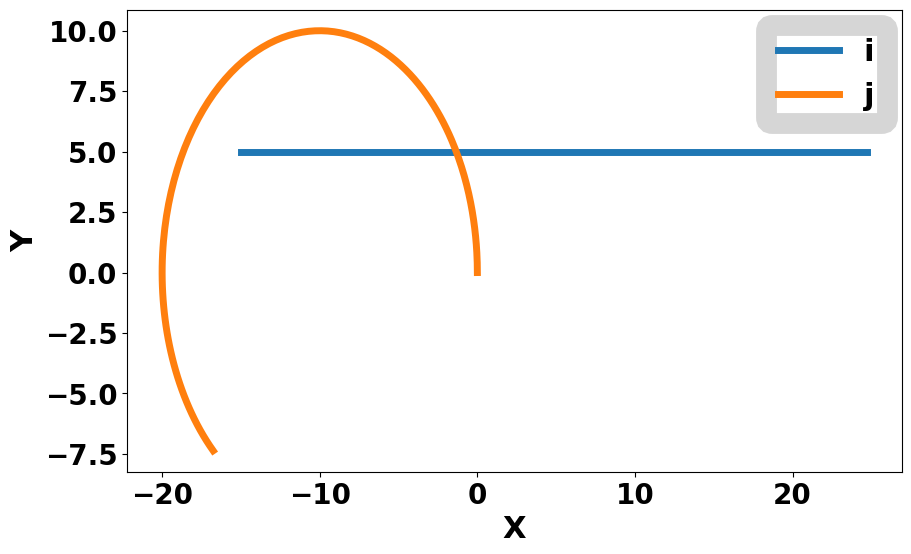}
    \includegraphics[width=0.49\linewidth]{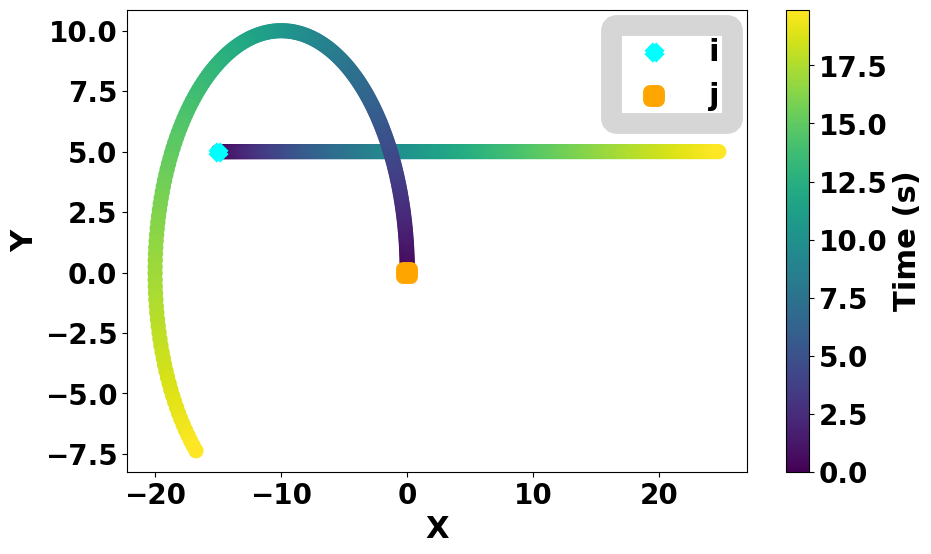}
    \caption{Estimated Trajectories of Vehicles with Initial conditions (\edit{Position: $\+{p}_{\circ,i} = \elements{-15,5}$, $\+{p}_{\circ,j} = \elements{0,0}$; Velocity: $\+{v}_{\circ,i} = \elements{1,0}$, $\+{v}_{\circ,j} = \elements{0,1}$; Acceleration: $\+{a}_{\circ,i} = \elements{0.1,0}$, $\+{a}_{\circ,j} = \elements{-0.1,0.1}$}) in Computation of 2D-TTC for Scenario 4. \edit{The l}eft side of the figure shows the overall trajectories of each vehicle\edit{, estimated as a straight line trajectory for vehicle $i$ and a circular trajectory for vehicle $j$ at time $t_\circ$ under} the second-order TTC model. \edit{The right side of the figure shows the progression of the vehicles' location from time $t_\circ$ on the estimated trajectories under the second-order TTC scheme through the color gradient.} The cyan mark and orange mark represent the initial position of vehicle $i$ and vehicle $j$ respectively.}
    \label{fig:scen4_2dttc}
\end{figure}

\edit{Figure \ref{fig:scen4_2dttc} shows $\+{p}_{i}(t;t_\circ)$ and $\+{p}_{j}(t; t_\circ)$ in scenario 4 for the computation of 2D-TTC.}

As the estimated trajectories in Figure \ref{fig:scen4_2dttc} depict, 2D-TTC is capable of estimating the turning of vehicle $j$. In addition, the left side of the figure show\edit{s} that \edit{the} two vehicles are approaching closely and might be colliding. 

To further investigate whether the collision is expected to happen between two vehicles under the second-order TTC scheme, we plot $\textit{d}_{ij}-\phi$ \edit{along} the estimated trajectories of vehicles. 
\begin{figure} [ht!]
    \centering
    \includegraphics[width=0.49\linewidth]{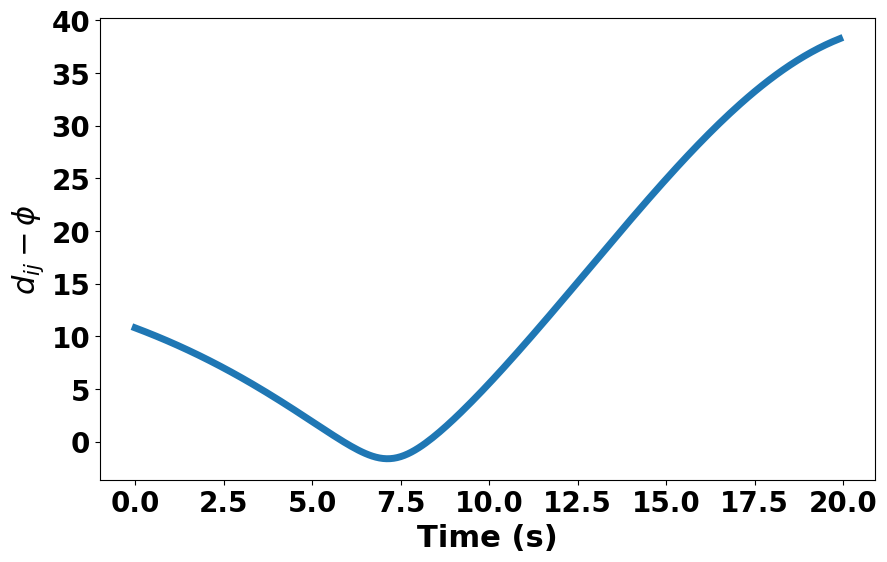}
    \caption{Predicted $\textit{d}_{ij}-\phi$ based on the initial conditions as time progresses in the computation of 2D-TTC \edit{at time $t_\circ$} in  Scenario 4\edit{, estimating a collision.}}
    \label{fig:scen4_2dttc_d}
\end{figure}

Based on the $\textit{d}_{ij}-\phi$ \edit{along} estimated trajectories, \edit{the} two vehicles are estimated to approach each other \edit{and} to collide as $\textit{d}_{ij}-\phi$ decreases and goes below 0 at 5.88 seconds \edit{after time $t_\circ$}. 

\begin{figure} [ht!]
    \centering
    \includegraphics[width=0.49\linewidth]{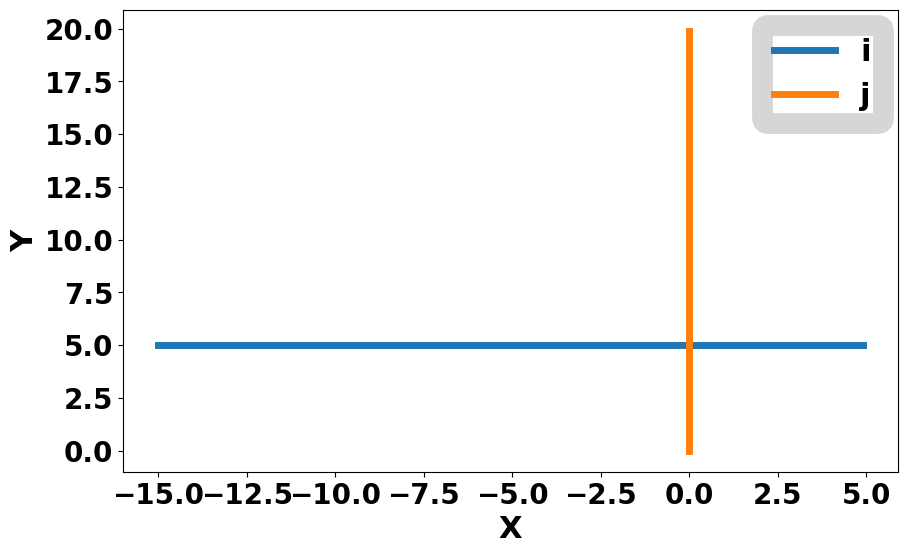}
    \includegraphics[width=0.49\linewidth]{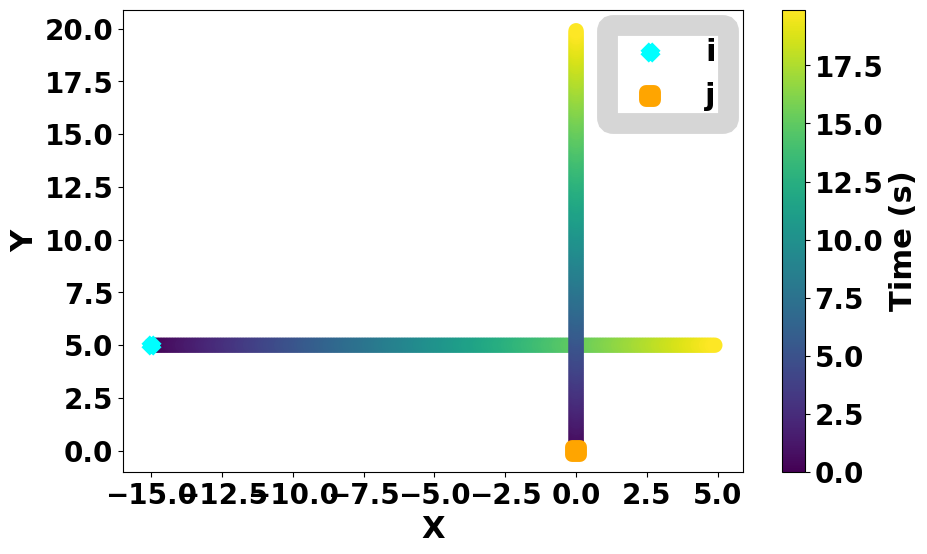}
    \caption{Estimated Trajectories of Vehicles with Initial conditions (same \edit{$\+{p}_{\circ, i}$, $\+{p}_{\circ, j}$, $\+{v}_{\circ, i}$, $\+{v}_{\circ, j}$} as Fig. \ref{fig:scen4_2dttc} except Acceleration: \edit{$\+{a}_{\circ, i} = \elements{0,0}$, $\+{a}_{\circ, j} = \elements{0,0}$}) in Computation of 1D-TTC for Scenario 4. \edit{The l}eft side of the figure shows the overall trajectories of each vehicle\edit{, estimated as straight line trajectories for both vehicles at time $t_\circ$ under} the first-order TTC model. \edit{The right side of the figure shows the progression of the vehicles' location from time $t_\circ$ on the estimated trajectories under the first-order TTC scheme through the color gradient.} The cyan mark and orange mark represent the initial position of vehicle $i$ and vehicle $j$ respectively.}
    \label{fig:scen4_1dttc}
\end{figure}

The trajectories of vehicles \edit{estimated at time $t_\circ$} in scenario 4 for the computation of 1D-TTC are shown in Figure \ref{fig:scen4_1dttc}.

Both vehicles are estimated to move in straight\edit{-}line trajectories under \edit{the} 1D-TTC scheme. \edit{Under these predicted trajectories}, the collision is unlikely to happen based on the color gradient on both estimated trajectories \edit{o}n the left side of Figure \ref{fig:scen4_1dttc}. 

\begin{figure} [ht!]
    \centering
    \includegraphics[width=0.49\linewidth]{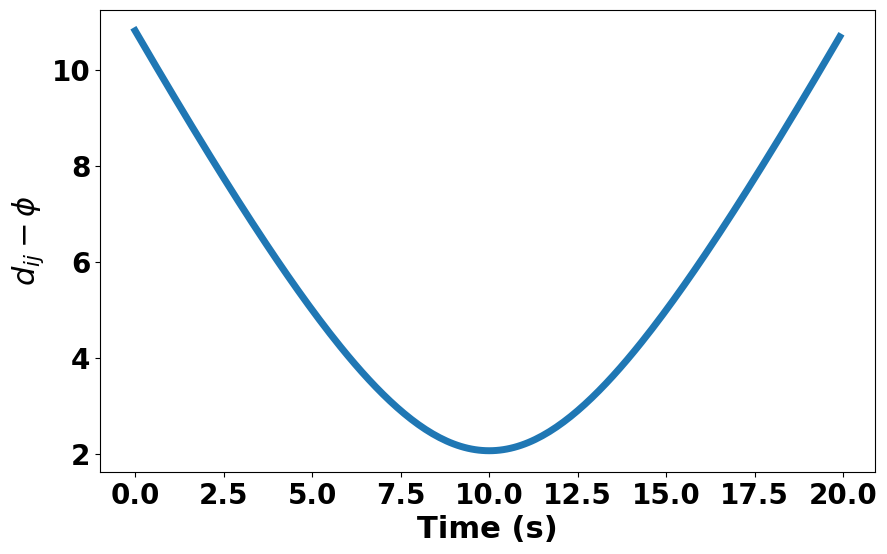}
    \caption{Predicted $\textit{d}_{ij}-\phi$ based on the initial conditions as time progresses in the computation of 1D-TTC \edit{at time $t_\circ$} in  Scenario 4\edit{, estimating no collision.}}
    \label{fig:scen4_1dttc_d}
\end{figure}

Figure \ref{fig:scen4_1dttc_d} shows $\textit{d}_{ij}-\phi$ \edit{along} the estimated trajectories of vehicles shown in Figure~\ref{fig:scen4_1dttc}. 

\edit{In Figure \ref{fig:scen4_1dttc_d},} $\textit{d}_{ij}-\phi$ decreases initially as two vehicles are approaching toward each other. However, the collision is estimated not to occur, and two vehicles further apart afterward. 

To summarize observations from Figure \ref{fig:scen4_2dttc}, \ref{fig:scen4_2dttc_d}, \ref{fig:scen4_1dttc}, and \ref{fig:scen4_1dttc_d}, \edit{based on estimated trajectories from initial conditions,} we expect a collision of both vehicles in the 2D-TTC scenario since one vehicle turns into the other as shown in Figure \ref{fig:scen4}. Unlike how well the 2D-TTC in Figure \ref{fig:scen4_2dttc}, the 1D-TTC scenario does not capture a collision as shown in Figure \ref{fig:scen4_1dttc}. In Figure \ref{fig:scen4_1dttc_d} and Figure \ref{fig:scen4_2dttc_d}, the collision is expected to happen under 2D-TTC while the collision expects not to happen under 1D-TTC. 

\begin{figure} [ht!]
    \centering
  \includegraphics[width=0.5\linewidth]{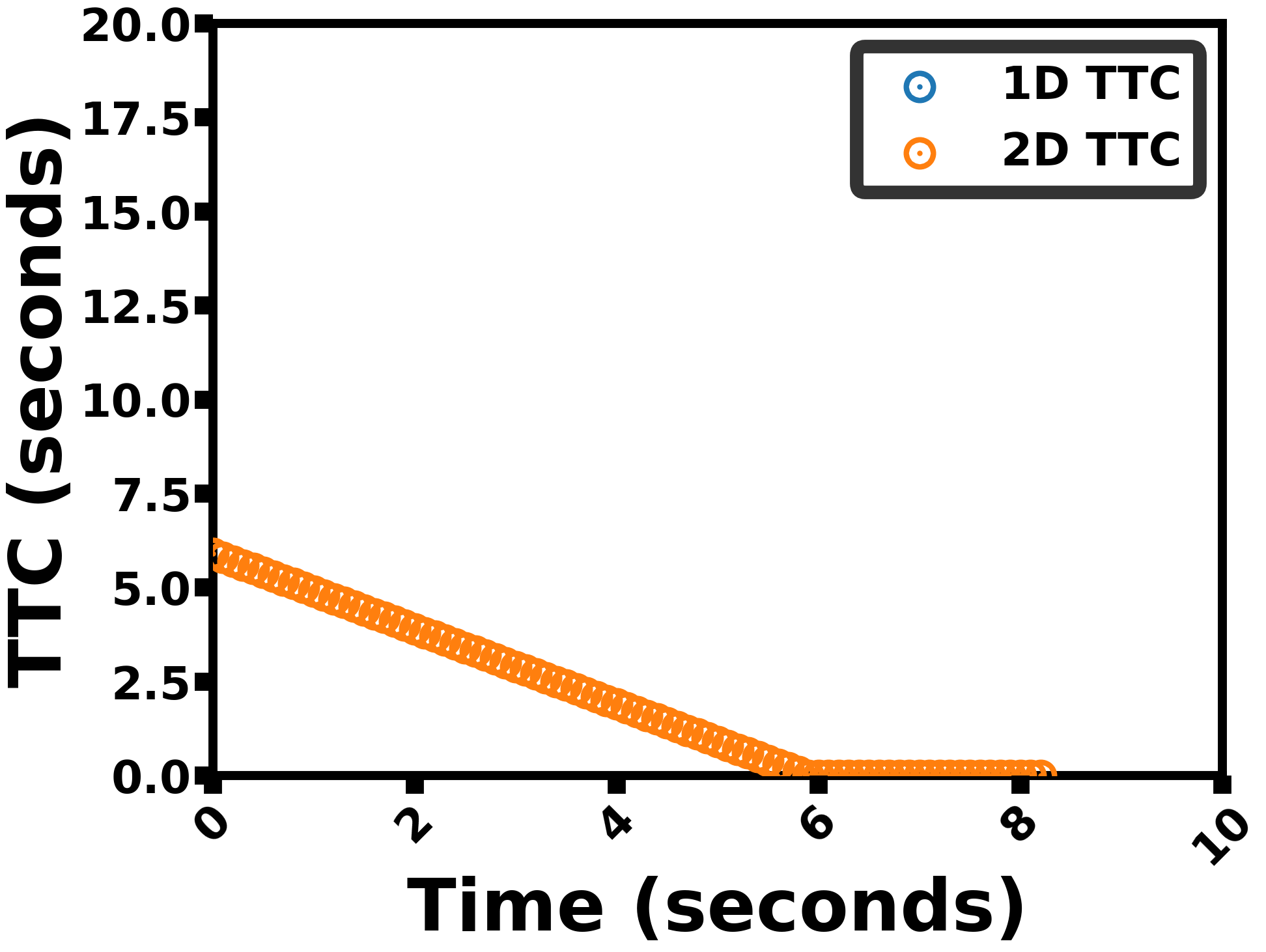}
    \caption{\edit{TTC values} for each timestep in Scenario 4. \edit{On the plot, only 2D-TTC exists as the collision between two vehicles is predicted continuously over time under the second-order TTC model while not under the first-order TTC model.}}
    \label{fig:traj_compare_scenario4}
\end{figure}

Now we compare the performance of the first-order TTC and the second-order TTC throughout the simulation horizon of Scenario 4. Figure \ref{fig:traj_compare_scenario4} is the plot of TTC values at each timestep of the simulation for Scenario 4 under both the first-order TTC and the second-order TTC. 

The second-order TTC scheme, unlike the first-order TTC scheme, expects a collision between vehicles to happen continuously in Figure \ref{fig:traj_compare_scenario4}. As Scenario 4 is designed so that a collision between two vehicles to happen, the second-order TTC scheme performs better. 

\subsection{Scenario 5: One Vehicle Turning Left \& Another Vehicle going Straight on Opposite Direction}

\begin{figure} [ht!]
    \centering
  \includegraphics[width=0.9\linewidth]{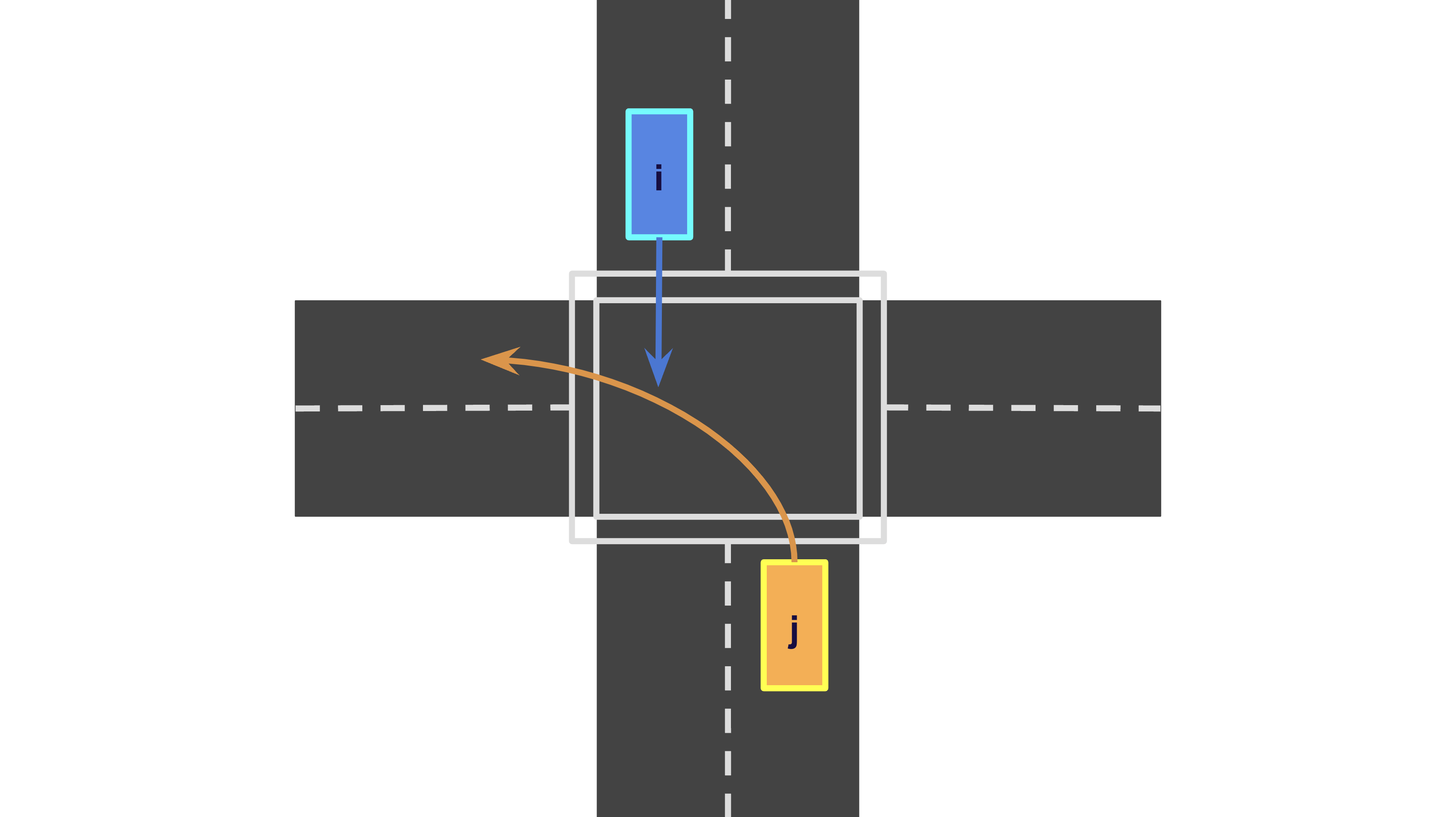}
    \caption{Scenario 5 illustration}
    \label{fig:scen5}
\end{figure}

In this scenario, one vehicle is about to turn left while another vehicle in the opposite direction goes straight, as seen in Figure \ref{fig:scen5}. 

The scenario in Figure \ref{fig:scen5} is designed so that the collision happen\edit{s} as vehicle $j$ almost complete\edit{s} the left turn and as vehicle $i$ enters the intersection. 

\begin{figure} [ht!]
    \centering
    \includegraphics[width=0.49\linewidth]{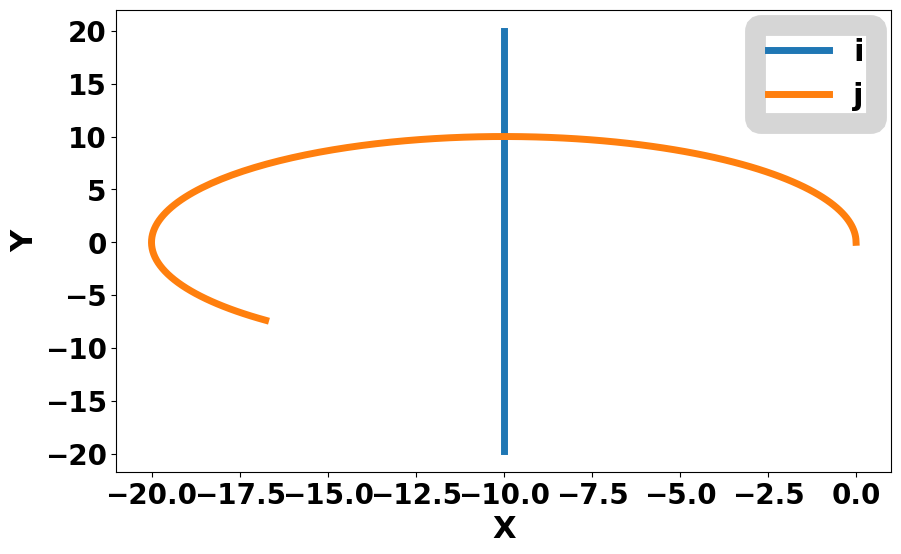}
    \includegraphics[width=0.49\linewidth]{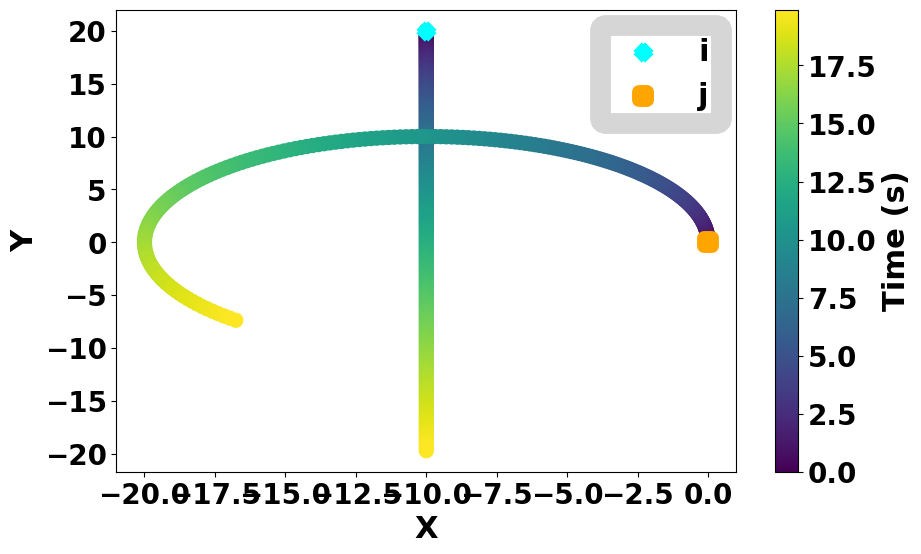}
    \caption{Estimated Trajectories of Vehicles with Initial conditions (\edit{Position: $\+{p}_{\circ,i} = \elements{-15,5}$, $\+{p}_{\circ,j} = \elements{0,0}$; Velocity: $\+{v}_{\circ,i} = \elements{1,0}$, $\+{v}_{\circ,j} = \elements{0,1}$; Acceleration: $\+{a}_{\circ,i} = \elements{0.1,0}$, $\+{a}_{\circ,j} = \elements{-0.1,0.1}$}) in Computation of 2D-TTC for Scenario 5. \edit{The l}eft side of the figure shows the overall trajectories of each vehicle\edit{, estimated as a straight line trajectory for vehicle $i$ and a circular trajectory for vehicle $j$ at time $t_\circ$ under} the second-order TTC model. \edit{The right side of the figure shows the progression of the vehicles' location from time $t_\circ$ on the estimated trajectories under the second-order TTC scheme through the color gradient.} The cyan mark and orange mark represent the initial position of vehicle $i$ and vehicle $j$ respectively. }
    \label{fig:scen5_2dttc}
\end{figure}

To investigate which TTC scheme is better for depicting the scenario, we find the estimated trajectories \edit{based on positions, velocities, and accelerations at time $t_\circ$}, $\textit{d}_{ij}-\phi$ \edit{along the estimated trajectories}, and TTC value under the second-order TTC. The estimated trajectories of vehicles in scenario 5 for the computation of 2D-TTC are shown in Figure \ref{fig:scen5_2dttc}.

As seen in all previous scenarios, Figure \ref{fig:scen5_2dttc} demonstrate\edit{s} the capability of the second-order TTC scheme in predicting the left turning of vehicle $j$.

To check if the collision is estimated to happen \edit{with $\+{p}_{i}(t)$ and $\+{p}_{j}(t)$ estimated at time $t_\circ$} under the second-order TTC scheme, we plot $\textit{d}_{ij}-\phi$ along the estimated trajectories as shown in Figure \ref{fig:scen5_2dttc}. 
\begin{figure} [ht!]
    \centering
    \includegraphics[width=0.49\linewidth]{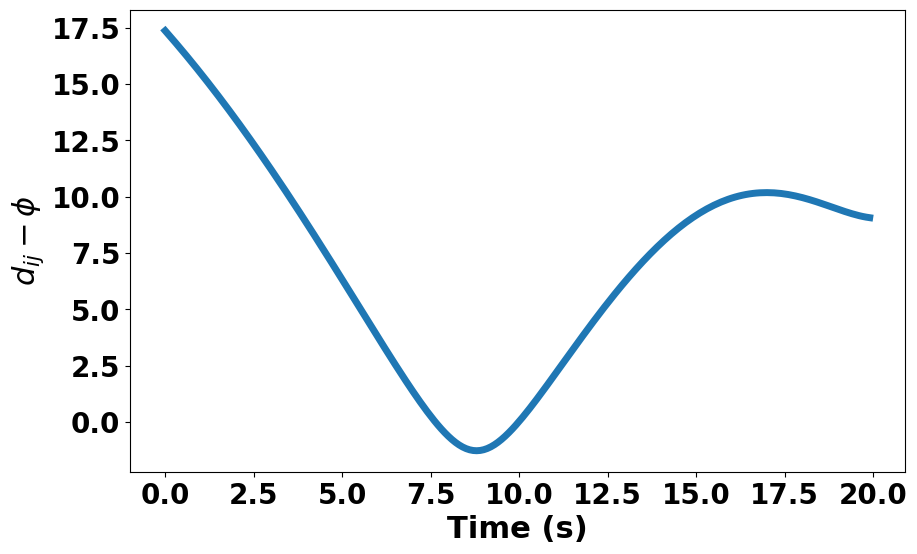}
    \caption{Predicted $\textit{d}_{ij}-\phi$ based on the initial conditions as time progresses in the computation of 2D-TTC \edit{at time $t_\circ$} in  Scenario 5\edit{, estimating a collision.}}
    \label{fig:scen5_2dttc_d}
\end{figure}

Based on Figure \ref{fig:scen5_2dttc_d}, the collision is expected to happen at 7.7 seconds \edit{after moving along the estimated trajectories} as $\textit{d}_{ij}-\phi$ goes below 0 after being decreased from the start. 

\begin{figure} [ht!]
    \centering
    \includegraphics[width=0.49\linewidth]{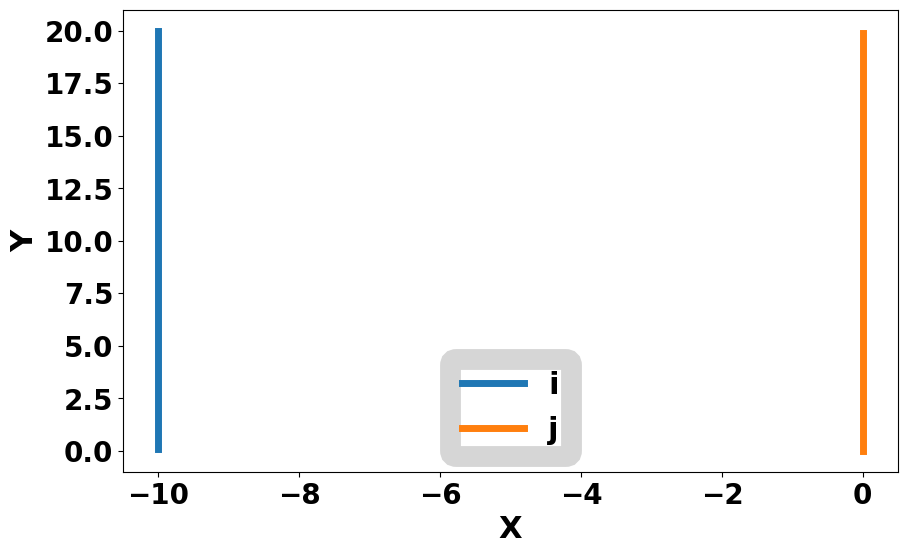}
    \includegraphics[width=0.49\linewidth]{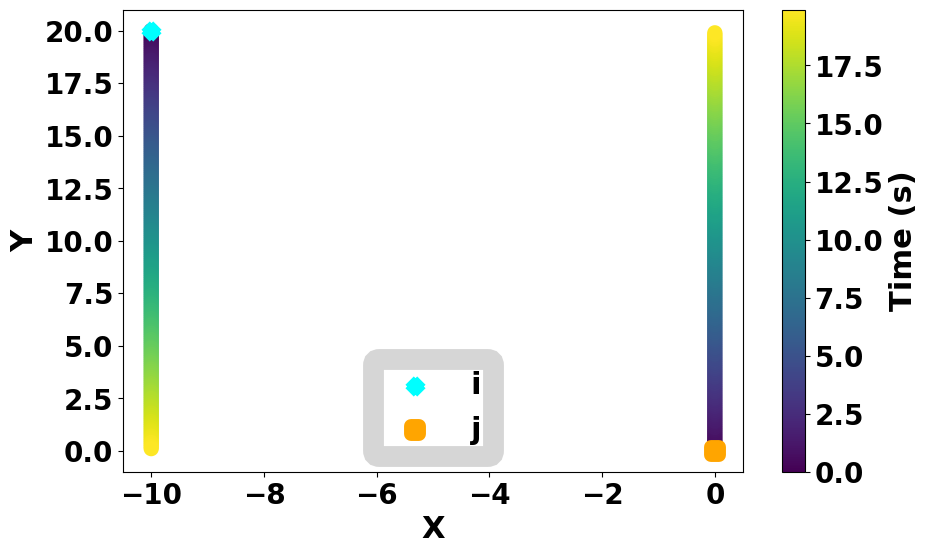}
    \caption{Estimated Trajectories of Vehicles with Initial conditions (same \edit{$\+{p}_{\circ, i}$, $\+{p}_{\circ, j}$, $\+{v}_{\circ, i}$, $\+{v}_{\circ, j}$} as Fig. \ref{fig:scen5_2dttc} except Acceleration: \edit{$\+{a}_{\circ, i} = \elements{0,0}$, $\+{a}_{\circ, j} = \elements{0,0}$}) in Computation of 1D-TTC for Scenario 5. \edit{The l}eft side of the figure shows the overall trajectories of each vehicle\edit{, estimated as straight line trajectories for both vehicles at time $t_\circ$ under} the first-order TTC model. \edit{The right side of the figure shows the progression of the vehicles' location from time $t_\circ$ on the estimated trajectories under the first-order TTC scheme through the color gradient.} The cyan mark and orange mark represent the initial position of vehicle $i$ and vehicle $j$ respectively.}
    \label{fig:scen5_1dttc}
\end{figure}

Now that we check how the first-order TTC performs for the scenario. The estimated trajectories of vehicles \edit{based on initial conditions} in scenario 5 for the computation of 1D-TTC are shown in Figure \ref{fig:scen5_1dttc}.

Unlike how the second-order TTC is capable of well-estimating the turning of vehicle $j$, the first-order TTC cannot estimate the turning of vehicle $j$ as it estimate\edit{s} vehicle $j$ to move in a straight line trajectory. As \edit{a} consequence, the first-order TTC does not expect the collision between two vehicles to happen, which does not align with the design of the scenario. 

\begin{figure} [ht!]
    \centering
    \includegraphics[width=0.49\linewidth]{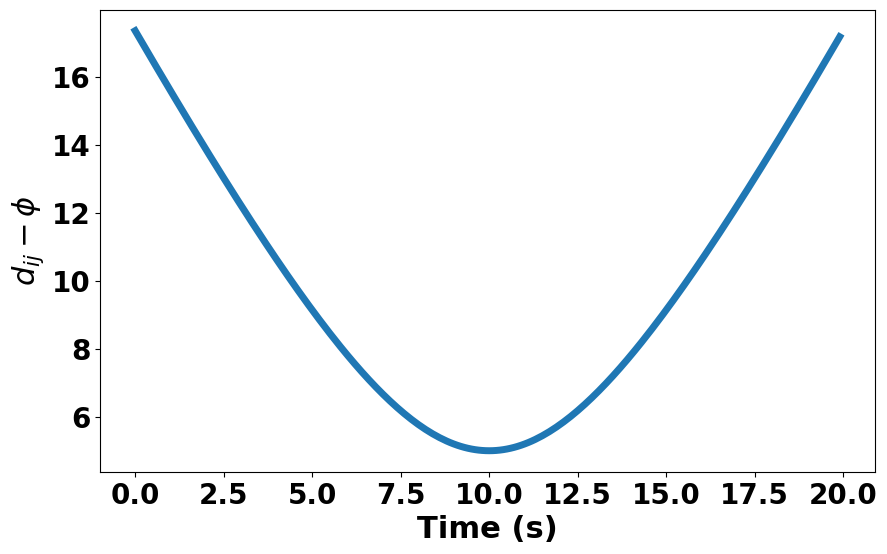}
    \caption{Predicted $\textit{d}_{ij}-\phi$ based on the initial conditions as time progresses in the computation of 1D-TTC \edit{at time $t_\circ$} in  Scenario 5\edit{, estimating no collision.}}
    \label{fig:scen5_1dttc_d}
\end{figure}

Figure \ref{fig:scen5_1dttc_d} shows $\textit{d}_{ij}-\phi$ \edit{along} the estimated trajectories of vehicles shown in Figure \ref{fig:scen5_1dttc}. 

\edit{As consistent with the previous observation, the v}ehicles initially approach close to each other, and vehicles move away from each other afterward\edit{with no collision happening}. 

We additionally compute TTC for each timestep under both the first-order TTC scheme and the second-order TTC scheme for comparison. 
\begin{figure} [ht!]
    \centering
  \includegraphics[width=0.5\linewidth]{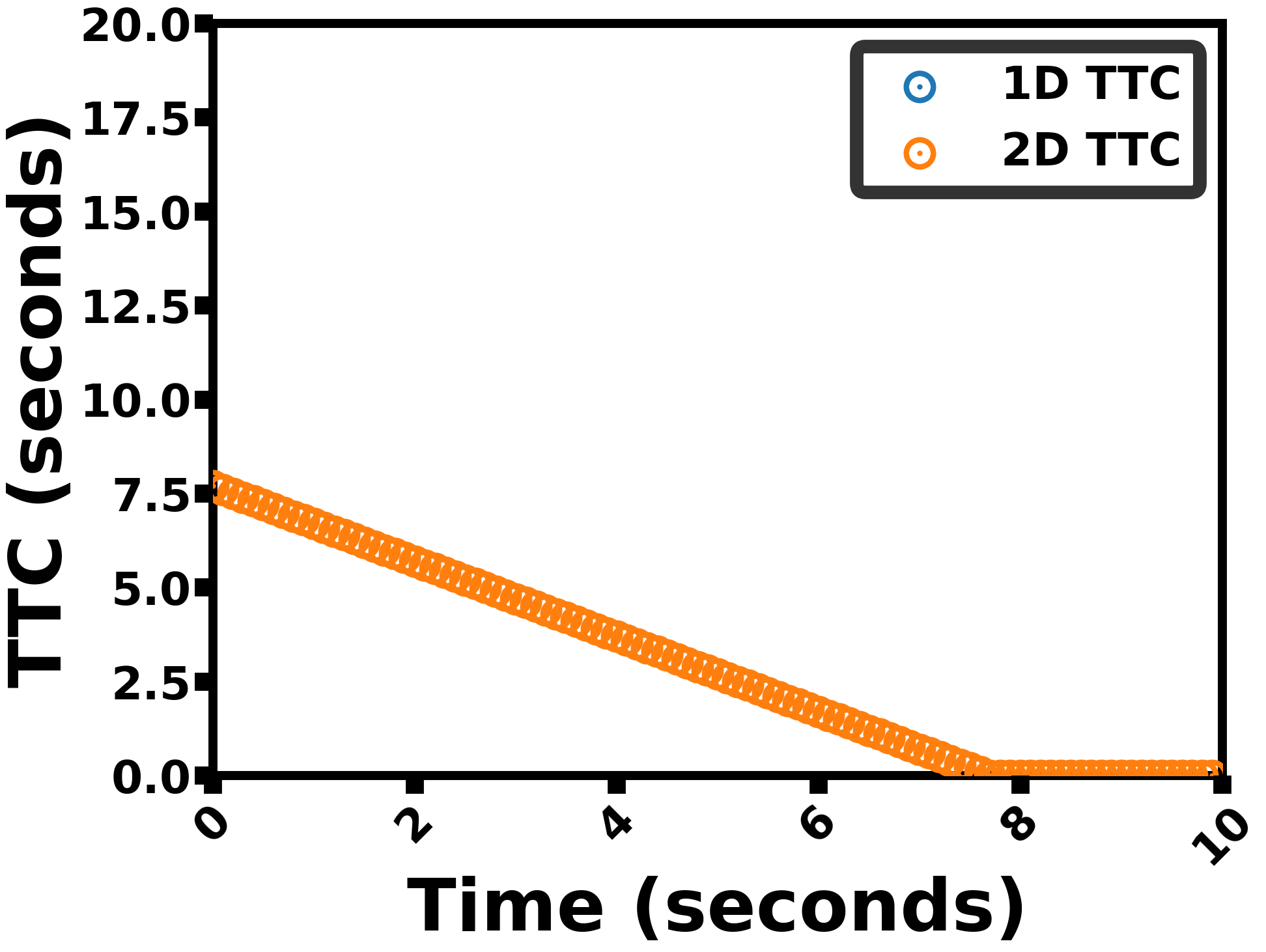}
    \caption{\edit{TTC values} for each timestep in Scenario 5. \edit{On the plot, only 2D-TTC exists as the collision between two vehicles is predicted continuously over time under the second-order TTC model while not under the first-order TTC model.}}
    \label{fig:traj_compare_scenario5}
\end{figure}

Based on Scenario 5 setup, we expect a collision of both vehicles in the 2D-TTC scenario as shown in Figure \ref{fig:scen5}, \ref{fig:scen5_2dttc}, and \ref{fig:traj_compare_scenario5}. However, the 1D-TTC scenario does not capture a collision as shown in Figure \ref{fig:scen5_1dttc} and \ref{fig:traj_compare_scenario5}. 
\\
\edit{\subsection{Summary of Numerical Simulation Analysis}}

As these various simulations show, 2D-TTC and 1D-TTC yield very different results when applied in the same scenarios. These differences determine whether a collision happens or not, which is crucial when considering intersection safety. These simulations prove that 2D-TTC is superior to 1D-TTC, as it captures a more accurate time to collision when turning is involved in either vehicle.

\section{Evaluation with Real-World Data} \label{sec: real_data}

In this section, we compare the effectiveness of 2D-TTC in terms of 1D-TTC with publicly accessible real-world data. 

We utilize the Argoverse 2 dataset, the real-world data released to the public for further advancing autonomous vehicle (AV) technologies \cite{Argoverse2}. Argoverse 2 dataset containing the trajectories of objects, including not only automobiles but also pedestrians and cyclists, in six different cities in the US \cite{Argoverse2, TrustButVerify}. Argoverse 2 dataset has more than 10,000 real-world cases with the details of road geometry data \cite{Argoverse2, TrustButVerify}. Almost all of the scenarios in \edit{the} Argoverse 2 dataset are scenarios with no collision between objects. 

\subsection{Selected Scenario Evaluation} \label{sec: specific_scenario}

\begin{figure} [ht!]
    \centering
  \includegraphics[width=0.8\linewidth]{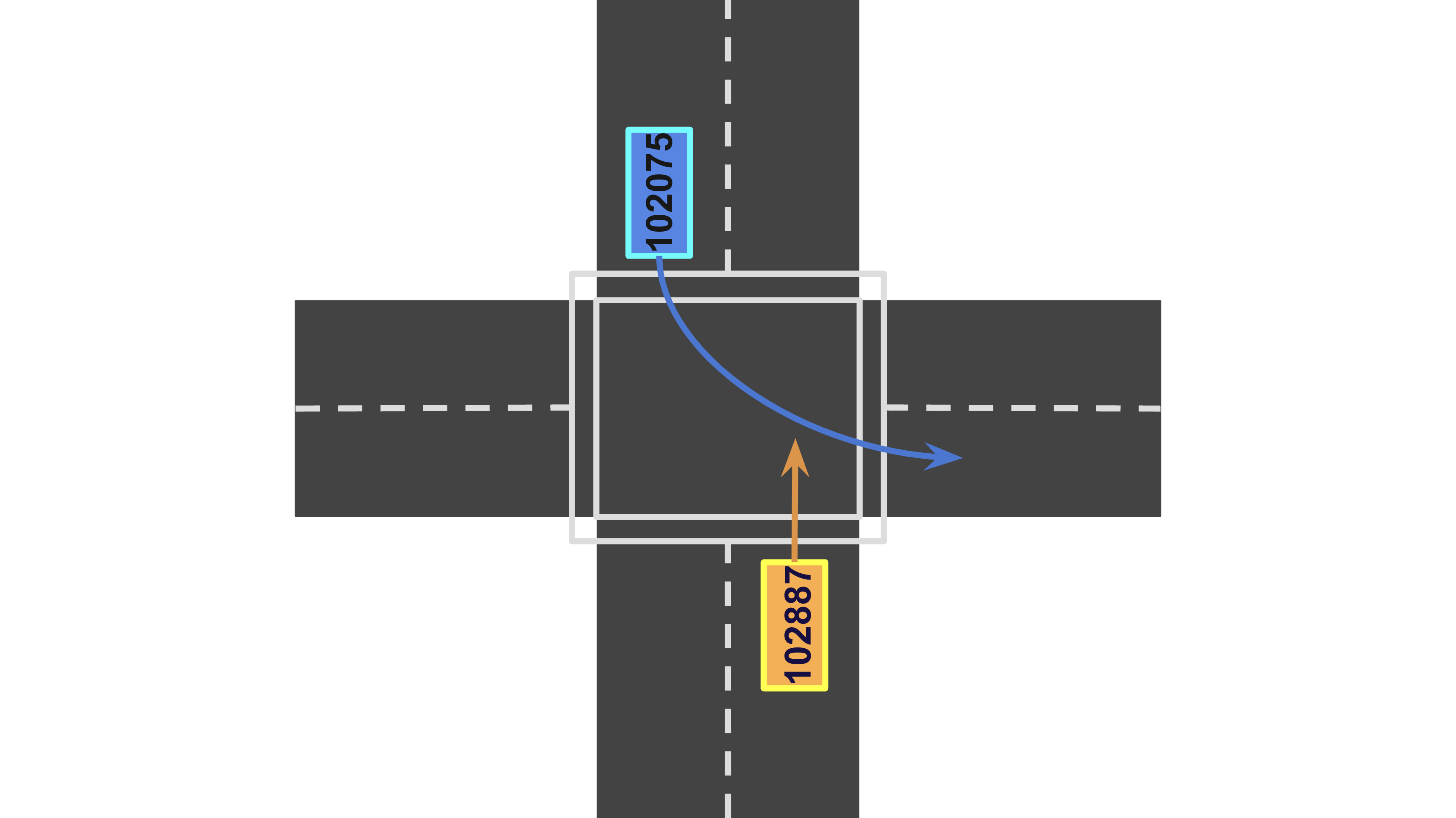}
    \caption{Movement of the vehicle 102075 and the vehicle 102887 on the selected scenario.}
    \label{fig:real_traj}
\end{figure}

We specifically use the trajectory data of vehicle 102075 and vehicle 102887 for scenario ID 00d9f613-2637-4665-b640-15ca2219929a. Figure \ref{fig:real_traj} depicts the movements of two vehicles in the selected scenario.


In the selected scenario, vehicle 102075 is turning left \edit{at} the intersection while vehicle 102887 is going straight. In the scenario with no collision, we examine the first-order TTC and second-order TTC. The orange arrow and the blue arrow represent the direction of vehicle 102075 and vehicle 102887, respectively. 

Among different features, we specifically use the position and velocity of vehicles in the $x$ and $y$ directions. The timestep of the dataset is 0.1s. As the dataset does not contain the acceleration in the $x$ and $y$ direction, we estimate the acceleration by computing the instantaneous acceleration\edit{, $\vec{\dot \zeta}$, along the actual trajectory defined as $\vec \xi$ and $\vec \zeta$, which a vehicle traveled in a real-world}.
\begin{equation*}
     \vec{\dot \zeta}(t) = \frac{\vec \zeta(t+0.1s) -\vec \zeta(t)}{0.1 s}  
\end{equation*}

Now that we have \edit{$\vec \xi$, $\vec \zeta$ and $\vec{\dot \zeta}$} for both vehicle $i$ and vehicle $j$, we can compute the first-order TTC and the second-order TTC. For each timestep, \edit{we determine $\+{p}(t)$, $\+{v}(t)$, and $\+{a}(t)$ for both vehicles under each TTC scheme based on  $\vec \xi$, $\vec \zeta$ and $\vec{\dot \zeta}$ at a corresponding timestep. Then, we compute first-order TTC and second-order TTC at that timestep} based on \eqref{E:collision_first} and Algorithm \ref{alg:cap} respectively. \edit{We continue computing both TTCs for every timestep to evaluate and compare the performance of both metrics for the collision risk assessment.}
\begin{figure} [ht!]
    \centering
  \includegraphics[width=0.5\linewidth]{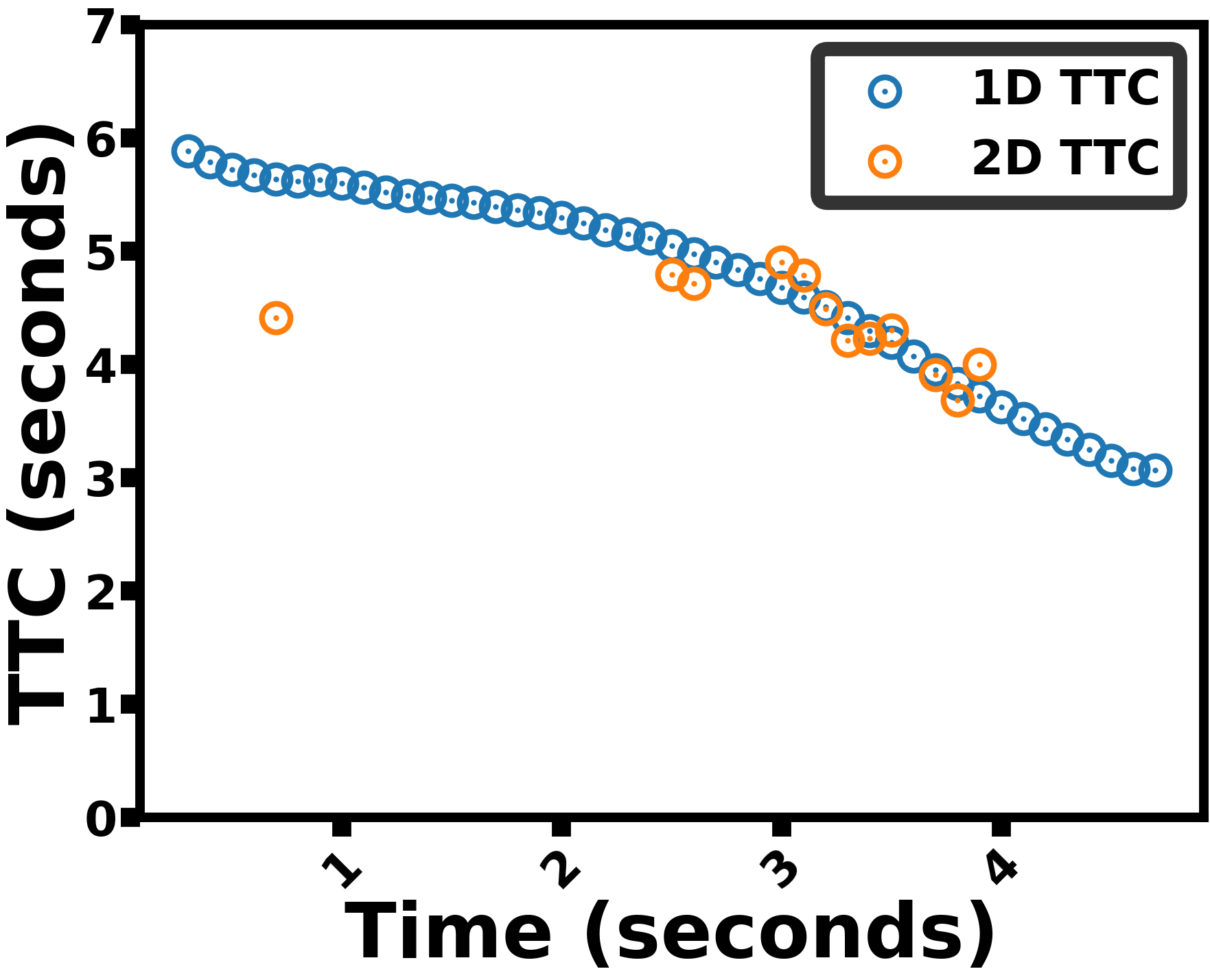}
    \caption{Comparison between the first-order TTC and the second-order TTC for each timestep in the selected real-world scenario.}
    \label{fig:real_traj_compare}
\end{figure}

From the comparison between the first-order TTC and the second-order TTC \edit{in} Figure \ref{fig:real_traj_compare}, the first-order TTC detects many more potential collisions compared to the second-order TTC and with smaller TTC values. We know that the vehicles safely traverse the intersection for which 2D-TTC produces fewer false alarm\edit{s} for collision risk compared to 1D-TTC. In other words, since a collision does not happen in the selected scenario, it is possible to conclude that the second-order TTC performs better than the first-order TTC in evaluating the collision time and assessing risk. 

\edit{It is important to note that the second-order TTC model is an estimation model and, as such, may generate a few false alarms. Additionally, depending on the positions, velocities, and accelerations of vehicles, the TTC values produced by the second-order TTC model can be smaller than the first-order TTC model at certain timesteps. However, overall, we can conclude that the second-order TTC model provides a more reliable safety assessment compared to the first-order TTC model as it generally produces fewer false alarms.}

\subsection{Multiple Scenarios Evaluation}

We now further validate the superiority of the second-order TTC to the first-order TTC \edit{for} accurate collision risk assessment with multiple real-world scenarios from the same dataset. We choose the scenarios where any of \edit{the} two vehicles are turning either right or left. 

From multiple scenarios, we check \edit{a} total number of non-infinite value\edit{s} from the first-order TTC and the second-order TTC. We then look at the distribution of TTC less than 5 seconds for the first-order TTC and the second-order TTC, similar to Section \ref{sec: specific_scenario}. 


To assess how accurately the second-order TTC captures the collision risk compared to the first-order TTC, we consider the critical value of TTC. The critical value of TTC differentiates safe situations or dangerous situations in terms of collision \cite{minderhoud2001extended}. Different past works propose the critical value in the range of 2 seconds to 5 seconds for accurately determining risky encounters \cite{van1991time, hogema1996effects, hirst2020format}. As the critical value of TTC is 2–5 seconds, we only check TTC values less than 5 seconds to check if 1D-TTC or 2D-TTC determines the situation to be risky or not in real-world scenarios.

Out of 1819 TTC computations in 21 scenarios, the second order-TTC produce\edit{s} 68 non-infinite TTC value\edit{s} less than 5 seconds while the first order-TTC produce\edit{s} 91 non-infinite TTC values less than 5 seconds. As scenarios in Argoverse 2 have no collision, this suggests that the second-order TTC performs more accurate collision risk assessments and better distinguishes between safe vehicle interactions and conflict situations. 

\begin{figure}[ht!]
    \centering
    \begin{subfigure}[b]{0.4\textwidth}
        \centering
        \includegraphics[width=\linewidth]{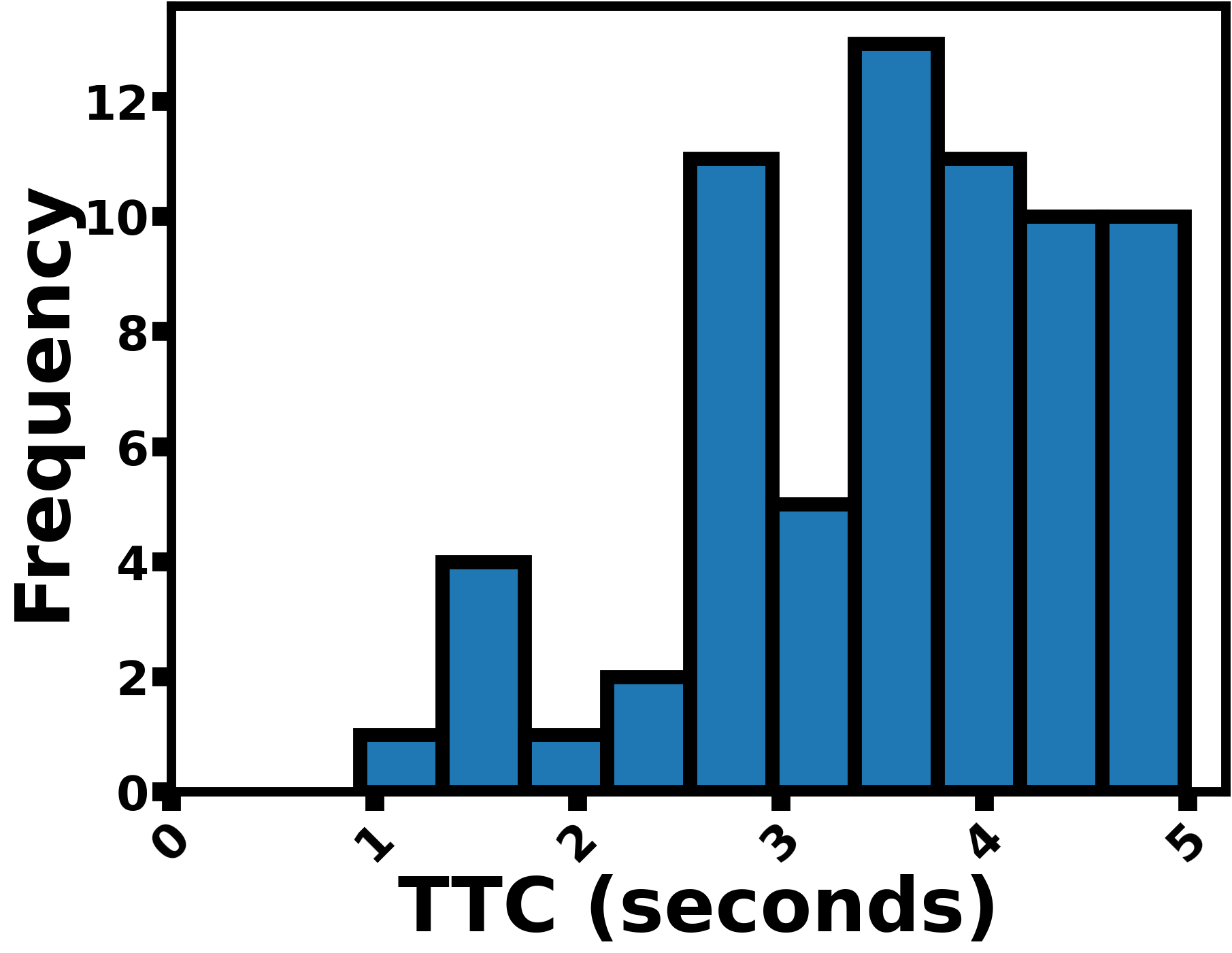}
        \caption{\edit{2D-TTC}}
        \label{fig:ttc_dist_2d}
    \end{subfigure}
    \hspace{2ex}
    \begin{subfigure}[b]{0.4\textwidth}
        \centering
        \includegraphics[width=\linewidth]{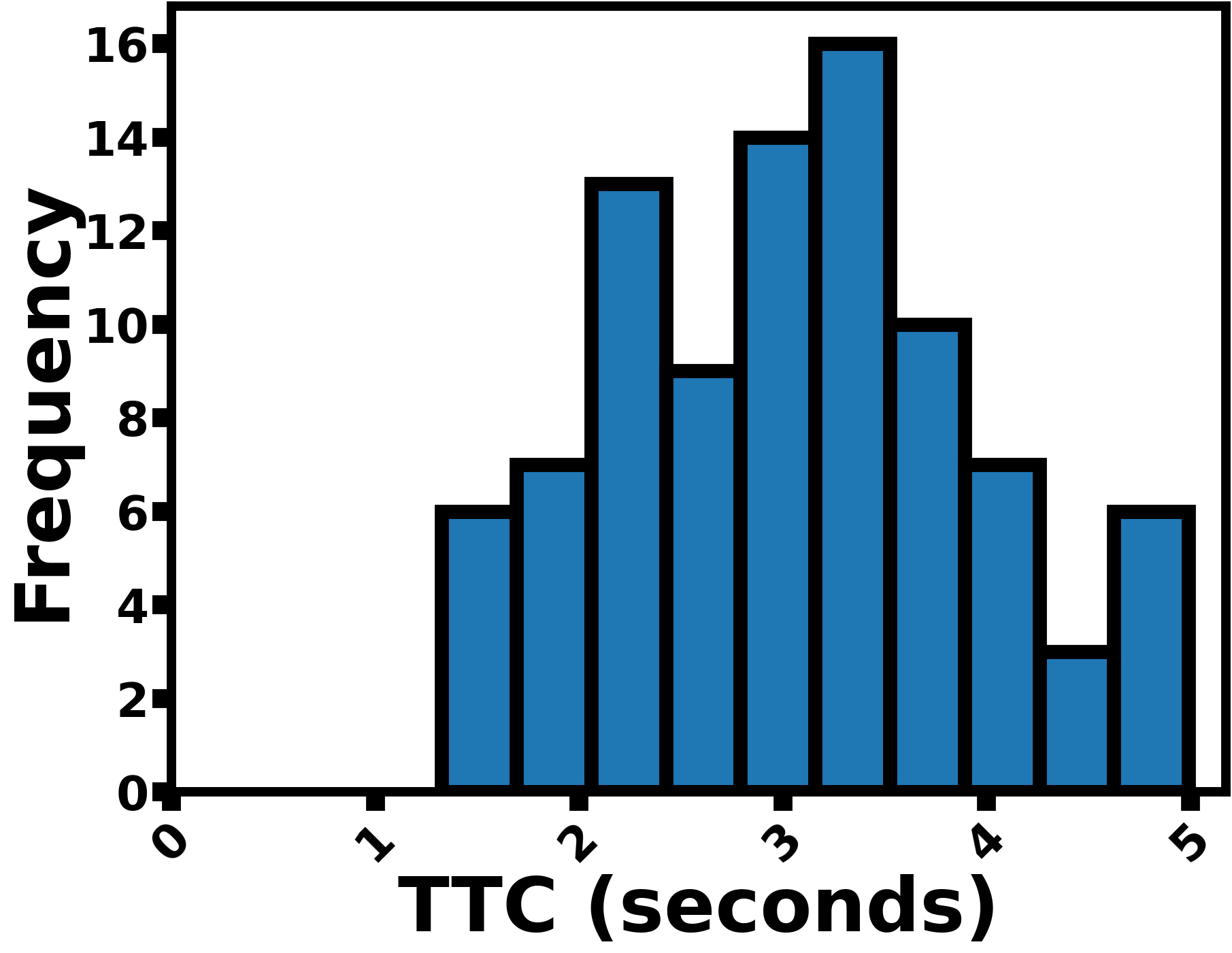}
        \caption{\edit{1D-TTC}}
        \label{fig:ttc_dist_1d}
    \end{subfigure}
    \caption{Distribution of TTC less than 5 seconds from the second-order TTC (Left) and the first-order TTC (Right).}
    \label{fig:1d_2dttc_distribution}
\end{figure}

Figure \ref{fig:1d_2dttc_distribution} shows the overall distribution of non-infinite TTC based on the second-order TTC (Left) and first-order TTC (Right). 







From Figure \ref{fig:1d_2dttc_distribution}, the first-order TTC produces TTC values smaller than 5 seconds much more frequently than the second-order TTC as we see in the selected scenario in Section \ref{sec: specific_scenario}. In addition, Figure \ref{fig:1d_2dttc_distribution} shows that the distribution of TTC less than 5 seconds is left-skewed for the second-order TTC while right-skewed for the first-order TTC. In other words, the first-order TTC identifies more potential collision moments. As we had observed in Section~\ref{sec: specific_scenario} many of the first-order TTC potential collisions are spurious simply due to the oversimplified kinematics model. Moreover, the left-skewness of the first-order TTC contradicts the expected frequency of traffic events from previous studies~\cite{chin1997measurement}. Since collisions do not happen in scenarios from \edit{the} Argoverse 2 Dataset, these results support that the second-order TTC provide\edit{s} a more accurate characterization of collision risk with fewer false alarms. 

\section{Conclusions/Future Work}

This work presents a new second-order TTC formulation that considers non-state acceleration and turning. The formulation is based on the assumptions corresponding to the well-known ``bicycle model," with its constant turn rate and acceleration \edit{magnitude} kinematics~\cite{lefevre2014survey}. Provided that there is no slip, as is to be expected under normal driving conditions, this kinematic model results in more realistic predictions of vehicle trajectories for detecting possible collisions. As a result, the second-order TTC more accurately reflects the safety risk in vehicular interactions, as demonstrated in Section~\ref{sec: analysis}. It appropriately characterizes situations in which vehicles might come close but safe\edit{l}y negotiate their trajectories and identifies potential collisions much earlier than the first-order TTC when turning is involve\edit{d}, such as in roundabouts. Comparisons in representative scenarios and scenarios from real-world data demonstrate that the proposed second-order TTC is a more representative safety measure, resulting in better characterizations of safety risk than the first-order TTC.

The estimation of second-order TTC is somewhat more involved due to the nonlinear vehicle trajectories. To address this we presented an algorithm to quickly and accurately estimate the second-order TTC (cf.~Section~\ref{sec: algorithm}). Although standard simulation-based time-stepping could be used, as Section~\ref{sec: evaluation} shows, the proposed algorithm is equally accurate and much faster.

There are several possibilities for future work. It could include the further assessment of the second-order TTC in additional real-world trajectories and its impact \edit{on} transportation safety studies or evaluation of AV planning algorithms. Additionally, kinematic models that further relax the constant turn rate and acceleration assumptions could be explored, as well as algorithms that take into account the actual vehicle geometry and dimensions. \edit{Consideration of human reactions and \edit{the} stochastic nature of driving behavior would also improve the validity of the measures and lead to more accurate safety assessment.}








\bmhead{Acknowledgements}
This work is a part of the Berkeley DeepDrive Project ``Collision Indeterminacy Prediction via Stochastic Trajectory Generation." Yuneil Yeo is partially supported by the Dwight David Eisenhower Transportation Fellowship Program Under Grant No. 693JJ32445085. 

\bibliography{bibliography} 


\begin{thebibliography}{42}
\ifx \bisbn   \undefined \def \bisbn  #1{ISBN #1}\fi
\ifx \binits  \undefined \def \binits#1{#1}\fi
\ifx \bauthor  \undefined \def \bauthor#1{#1}\fi
\ifx \batitle  \undefined \def \batitle#1{#1}\fi
\ifx \bjtitle  \undefined \def \bjtitle#1{#1}\fi
\ifx \bvolume  \undefined \def \bvolume#1{\textbf{#1}}\fi
\ifx \byear  \undefined \def \byear#1{#1}\fi
\ifx \bissue  \undefined \def \bissue#1{#1}\fi
\ifx \bfpage  \undefined \def \bfpage#1{#1}\fi
\ifx \blpage  \undefined \def \blpage #1{#1}\fi
\ifx \burl  \undefined \def \burl#1{\textsf{#1}}\fi
\ifx \doiurl  \undefined \def \doiurl#1{\url{https://doi.org/#1}}\fi
\ifx \betal  \undefined \def \betal{\textit{et al.}}\fi
\ifx \binstitute  \undefined \def \binstitute#1{#1}\fi
\ifx \binstitutionaled  \undefined \def \binstitutionaled#1{#1}\fi
\ifx \bctitle  \undefined \def \bctitle#1{#1}\fi
\ifx \beditor  \undefined \def \beditor#1{#1}\fi
\ifx \bpublisher  \undefined \def \bpublisher#1{#1}\fi
\ifx \bbtitle  \undefined \def \bbtitle#1{#1}\fi
\ifx \bedition  \undefined \def \bedition#1{#1}\fi
\ifx \bseriesno  \undefined \def \bseriesno#1{#1}\fi
\ifx \blocation  \undefined \def \blocation#1{#1}\fi
\ifx \bsertitle  \undefined \def \bsertitle#1{#1}\fi
\ifx \bsnm \undefined \def \bsnm#1{#1}\fi
\ifx \bsuffix \undefined \def \bsuffix#1{#1}\fi
\ifx \bparticle \undefined \def \bparticle#1{#1}\fi
\ifx \barticle \undefined \def \barticle#1{#1}\fi
\bibcommenthead
\ifx \bconfdate \undefined \def \bconfdate #1{#1}\fi
\ifx \botherref \undefined \def \botherref #1{#1}\fi
\ifx \url \undefined \def \url#1{\textsf{#1}}\fi
\ifx \bchapter \undefined \def \bchapter#1{#1}\fi
\ifx \bbook \undefined \def \bbook#1{#1}\fi
\ifx \bcomment \undefined \def \bcomment#1{#1}\fi
\ifx \oauthor \undefined \def \oauthor#1{#1}\fi
\ifx \citeauthoryear \undefined \def \citeauthoryear#1{#1}\fi
\ifx \endbibitem  \undefined \def \endbibitem {}\fi
\ifx \bconflocation  \undefined \def \bconflocation#1{#1}\fi
\ifx \arxivurl  \undefined \def \arxivurl#1{\textsf{#1}}\fi
\csname PreBibitemsHook\endcsname

\bibitem[\protect\citeauthoryear{Nikolaou et~al.}{2023}]{nikolaou2023exploiting}
\begin{barticle}
\bauthor{\bsnm{Nikolaou}, \binits{D.}},
\bauthor{\bsnm{Dragomanovits}, \binits{A.}},
\bauthor{\bsnm{Ziakopoulos}, \binits{A.}},
\bauthor{\bsnm{Deliali}, \binits{A.}},
\bauthor{\bsnm{Handanos}, \binits{I.}},
\bauthor{\bsnm{Karadimas}, \binits{C.}},
\bauthor{\bsnm{Kostoulas}, \binits{G.}},
\bauthor{\bsnm{Frantzola}, \binits{E.K.}},
\bauthor{\bsnm{Yannis}, \binits{G.}}:
\batitle{Exploiting surrogate safety measures and road design characteristics towards crash investigations in motorway segments}.
\bjtitle{Infrastructures}
\bvolume{8}(\bissue{3}),
\bfpage{40}
(\byear{2023})
\end{barticle}
\endbibitem

\bibitem[\protect\citeauthoryear{Wang et~al.}{2021}]{wang2021review}
\begin{barticle}
\bauthor{\bsnm{Wang}, \binits{C.}},
\bauthor{\bsnm{Xie}, \binits{Y.}},
\bauthor{\bsnm{Huang}, \binits{H.}},
\bauthor{\bsnm{Liu}, \binits{P.}}:
\batitle{A review of surrogate safety measures and their applications in connected and automated vehicles safety modeling}.
\bjtitle{Accident Analysis \& Prevention}
\bvolume{157},
\bfpage{106157}
(\byear{2021})
\end{barticle}
\endbibitem

\bibitem[\protect\citeauthoryear{Johnsson et~al.}{2018}]{johnsson2018search}
\begin{barticle}
\bauthor{\bsnm{Johnsson}, \binits{C.}},
\bauthor{\bsnm{Laureshyn}, \binits{A.}},
\bauthor{\bsnm{De~Ceunynck}, \binits{T.}}:
\batitle{In search of surrogate safety indicators for vulnerable road users: a review of surrogate safety indicators}.
\bjtitle{Transport reviews}
\bvolume{38}(\bissue{6}),
\bfpage{765}--\blpage{785}
(\byear{2018})
\end{barticle}
\endbibitem

\bibitem[\protect\citeauthoryear{Arun et~al.}{2021}]{arun2021systematic}
\begin{barticle}
\bauthor{\bsnm{Arun}, \binits{A.}},
\bauthor{\bsnm{Haque}, \binits{M.M.}},
\bauthor{\bsnm{Washington}, \binits{S.}},
\bauthor{\bsnm{Sayed}, \binits{T.}},
\bauthor{\bsnm{Mannering}, \binits{F.}}:
\batitle{A systematic review of traffic conflict-based safety measures with a focus on application context}.
\bjtitle{Analytic methods in accident research}
\bvolume{32},
\bfpage{100185}
(\byear{2021})
\end{barticle}
\endbibitem

\bibitem[\protect\citeauthoryear{Hayward}{1971}]{hayward1971near}
\begin{botherref}
\oauthor{\bsnm{Hayward}, \binits{J.}}:
Near misses as a measure of safety at urban intersections.
Pennsylvania Transportation and Traffic Safety Center
(1971)
\end{botherref}
\endbibitem

\bibitem[\protect\citeauthoryear{Hou et~al.}{2015}]{hou2015new}
\begin{barticle}
\bauthor{\bsnm{Hou}, \binits{J.}},
\bauthor{\bsnm{List}, \binits{G.F.}},
\bauthor{\bsnm{Guo}, \binits{X.}}:
\batitle{New algorithms for computing the time-to-collision in freeway traffic simulation models}.
\bjtitle{Computational intelligence and neuroscience}
\bvolume{2014},
\bfpage{57}--\blpage{57}
(\byear{2015})
\end{barticle}
\endbibitem

\bibitem[\protect\citeauthoryear{Tak et~al.}{2018}]{tak2018comparison}
\begin{botherref}
\oauthor{\bsnm{Tak}, \binits{S.}},
\oauthor{\bsnm{Kim}, \binits{S.}},
\oauthor{\bsnm{Lee}, \binits{D.}},
\oauthor{\bsnm{Yeo}, \binits{H.}}, et al.:
A comparison analysis of surrogate safety measures with car-following perspectives for advanced driver assistance system.
Journal of Advanced Transportation
\textbf{2018}
(2018)
\end{botherref}
\endbibitem

\bibitem[\protect\citeauthoryear{Goudarzi and Hassanzadeh}{2024}]{goudarzi2024collision}
\begin{barticle}
\bauthor{\bsnm{Goudarzi}, \binits{P.}},
\bauthor{\bsnm{Hassanzadeh}, \binits{B.}}:
\batitle{Collision risk in autonomous vehicles: classification, challenges, and open research areas}.
\bjtitle{Vehicles}
\bvolume{6}(\bissue{1}),
\bfpage{157}--\blpage{190}
(\byear{2024})
\end{barticle}
\endbibitem

\bibitem[\protect\citeauthoryear{Zhang et~al.}{2021}]{zhang2021evaluating}
\begin{barticle}
\bauthor{\bsnm{Zhang}, \binits{H.}},
\bauthor{\bsnm{Hou}, \binits{N.}},
\bauthor{\bsnm{Zhang}, \binits{J.}},
\bauthor{\bsnm{Li}, \binits{X.}},
\bauthor{\bsnm{Huang}, \binits{Y.}}:
\batitle{Evaluating the safety impact of connected and autonomous vehicles with lane management on freeway crash hotspots using the surrogate safety assessment model}.
\bjtitle{Journal of advanced transportation}
\bvolume{2021},
\bfpage{1}--\blpage{14}
(\byear{2021})
\end{barticle}
\endbibitem

\bibitem[\protect\citeauthoryear{Wachenfeld et~al.}{2016}]{wachenfeld2016worst}
\begin{bchapter}
\bauthor{\bsnm{Wachenfeld}, \binits{W.}},
\bauthor{\bsnm{Junietz}, \binits{P.}},
\bauthor{\bsnm{Wenzel}, \binits{R.}},
\bauthor{\bsnm{Winner}, \binits{H.}}:
\bctitle{The worst-time-to-collision metric for situation identification}.
In: \bbtitle{2016 IEEE Intelligent Vehicles Symposium (IV)},
pp. \bfpage{729}--\blpage{734}
(\byear{2016}).
\bcomment{IEEE}
\end{bchapter}
\endbibitem

\bibitem[\protect\citeauthoryear{Chin et~al.}{1991}]{chin1991traffic}
\begin{barticle}
\bauthor{\bsnm{Chin}, \binits{H.C.}},
\bauthor{\bsnm{Quek}, \binits{S.}},
\bauthor{\bsnm{Cheu}, \binits{R.}}:
\batitle{Traffic conflicts in expressway merging}.
\bjtitle{Journal of transportation engineering}
\bvolume{117}(\bissue{6}),
\bfpage{633}--\blpage{643}
(\byear{1991})
\end{barticle}
\endbibitem

\bibitem[\protect\citeauthoryear{Chin and Quek}{1997}]{chin1997measurement}
\begin{barticle}
\bauthor{\bsnm{Chin}, \binits{H.-C.}},
\bauthor{\bsnm{Quek}, \binits{S.-T.}}:
\batitle{Measurement of traffic conflicts}.
\bjtitle{Safety Science}
\bvolume{26}(\bissue{3}),
\bfpage{169}--\blpage{185}
(\byear{1997})
\end{barticle}
\endbibitem

\bibitem[\protect\citeauthoryear{Lakhal et~al.}{2019}]{lakhal2019risk}
\begin{barticle}
\bauthor{\bsnm{Lakhal}, \binits{N.M.B.}},
\bauthor{\bsnm{Adouane}, \binits{L.}},
\bauthor{\bsnm{Nasri}, \binits{O.}},
\bauthor{\bsnm{Slama}, \binits{J.B.H.}}:
\batitle{Risk management for intelligent vehicles based on interval analysis of ttc}.
\bjtitle{IFAC-PapersOnLine}
\bvolume{52}(\bissue{8}),
\bfpage{338}--\blpage{343}
(\byear{2019})
\end{barticle}
\endbibitem

\bibitem[\protect\citeauthoryear{Venthuruthiyil and Chunchu}{2022}]{venthuruthiyil2022anticipated}
\begin{barticle}
\bauthor{\bsnm{Venthuruthiyil}, \binits{S.P.}},
\bauthor{\bsnm{Chunchu}, \binits{M.}}:
\batitle{Anticipated collision time (act): A two-dimensional surrogate safety indicator for trajectory-based proactive safety assessment}.
\bjtitle{Transportation research part C: emerging technologies}
\bvolume{139},
\bfpage{103655}
(\byear{2022})
\end{barticle}
\endbibitem

\bibitem[\protect\citeauthoryear{Wakabayashi et~al.}{2003}]{wakabayashi2003traffic}
\begin{barticle}
\bauthor{\bsnm{Wakabayashi}, \binits{H.}},
\bauthor{\bsnm{Takahashi}, \binits{Y.}},
\bauthor{\bsnm{Niimi}, \binits{S.}},
\bauthor{\bsnm{Renge}, \binits{K.}}:
\batitle{Traffic conflict analysis using vehicle tracking system/digital vcr and proposal of a new conflict indicator}.
\bjtitle{Infrastructure Planning Review}
\bvolume{20},
\bfpage{949}--\blpage{956}
(\byear{2003})
\end{barticle}
\endbibitem

\bibitem[\protect\citeauthoryear{Brown}{2005}]{brown2005adjusted}
\begin{bchapter}
\bauthor{\bsnm{Brown}, \binits{T.L.}}:
\bctitle{Adjusted minimum time-to-collision (ttc): A robust approach to evaluating crash scenarios}.
In: \bbtitle{Proceedings of the Driving Simulation Conference North America},
vol. \bseriesno{40},
pp. \bfpage{40}--\blpage{48}
(\byear{2005}).
\bcomment{https://criticality-metrics.readthedocs.io/en/latest/references.html\#wachenfeld2016}
\end{bchapter}
\endbibitem

\bibitem[\protect\citeauthoryear{Kiefer et~al.}{2005}]{kiefer2005developing}
\begin{barticle}
\bauthor{\bsnm{Kiefer}, \binits{R.J.}},
\bauthor{\bsnm{LeBlanc}, \binits{D.J.}},
\bauthor{\bsnm{Flannagan}, \binits{C.A.}}:
\batitle{Developing an inverse time-to-collision crash alert timing approach based on drivers’ last-second braking and steering judgments}.
\bjtitle{Accident Analysis \& Prevention}
\bvolume{37}(\bissue{2}),
\bfpage{295}--\blpage{303}
(\byear{2005})
\end{barticle}
\endbibitem

\bibitem[\protect\citeauthoryear{Ward et~al.}{2015}]{ward2015extending}
\begin{barticle}
\bauthor{\bsnm{Ward}, \binits{J.R.}},
\bauthor{\bsnm{Agamennoni}, \binits{G.}},
\bauthor{\bsnm{Worrall}, \binits{S.}},
\bauthor{\bsnm{Bender}, \binits{A.}},
\bauthor{\bsnm{Nebot}, \binits{E.}}:
\batitle{Extending time to collision for probabilistic reasoning in general traffic scenarios}.
\bjtitle{Transportation Research Part C: Emerging Technologies}
\bvolume{51},
\bfpage{66}--\blpage{82}
(\byear{2015})
\end{barticle}
\endbibitem

\bibitem[\protect\citeauthoryear{Bando et~al.}{1995}]{bando1995dynamical}
\begin{barticle}
\bauthor{\bsnm{Bando}, \binits{M.}},
\bauthor{\bsnm{Hasebe}, \binits{K.}},
\bauthor{\bsnm{Nakayama}, \binits{A.}},
\bauthor{\bsnm{Shibata}, \binits{A.}},
\bauthor{\bsnm{Sugiyama}, \binits{Y.}}:
\batitle{Dynamical model of traffic congestion and numerical simulation}.
\bjtitle{Physical review E}
\bvolume{51}(\bissue{2}),
\bfpage{1035}
(\byear{1995})
\end{barticle}
\endbibitem

\bibitem[\protect\citeauthoryear{Wilson and Ward}{2011}]{wilson2011car}
\begin{barticle}
\bauthor{\bsnm{Wilson}, \binits{R.E.}},
\bauthor{\bsnm{Ward}, \binits{J.A.}}:
\batitle{Car-following models: fifty years of linear stability analysis--a mathematical perspective}.
\bjtitle{Transportation Planning and Technology}
\bvolume{34}(\bissue{1}),
\bfpage{3}--\blpage{18}
(\byear{2011})
\end{barticle}
\endbibitem

\bibitem[\protect\citeauthoryear{Pipes}{1953}]{pipes1953operational}
\begin{barticle}
\bauthor{\bsnm{Pipes}, \binits{L.A.}}:
\batitle{An operational analysis of traffic dynamics}.
\bjtitle{Journal of applied physics}
\bvolume{24}(\bissue{3}),
\bfpage{274}--\blpage{281}
(\byear{1953})
\end{barticle}
\endbibitem

\bibitem[\protect\citeauthoryear{Newell}{1961}]{newell1961nonlinear}
\begin{barticle}
\bauthor{\bsnm{Newell}, \binits{G.F.}}:
\batitle{Nonlinear effects in the dynamics of car following}.
\bjtitle{Operations research}
\bvolume{9}(\bissue{2}),
\bfpage{209}--\blpage{229}
(\byear{1961})
\end{barticle}
\endbibitem

\bibitem[\protect\citeauthoryear{Nick Zinat~Matin and Sowers}{2022}]{nick2022near}
\begin{barticle}
\bauthor{\bsnm{Nick Zinat~Matin}, \binits{H.}},
\bauthor{\bsnm{Sowers}, \binits{R.B.}}:
\batitle{Near-collision dynamics in a noisy car-following model}.
\bjtitle{SIAM Journal on Applied Mathematics}
\bvolume{82}(\bissue{6}),
\bfpage{2080}--\blpage{2110}
(\byear{2022})
\end{barticle}
\endbibitem

\bibitem[\protect\citeauthoryear{Chou et~al.}{2024}]{chou2024stability}
\begin{barticle}
\bauthor{\bsnm{Chou}, \binits{F.-C.}},
\bauthor{\bsnm{Keimer}, \binits{A.}},
\bauthor{\bsnm{Bayen}, \binits{A.M.}}:
\batitle{Stability of ring roads and string stability of car following models}.
\bjtitle{Mathematical Control and Related Fields}
\bvolume{14}(\bissue{4}),
\bfpage{1752}--\blpage{1775}
(\byear{2024})
\end{barticle}
\endbibitem

\bibitem[\protect\citeauthoryear{Matin and Delle~Monache}{2023}]{10384086}
\begin{bchapter}
\bauthor{\bsnm{Matin}, \binits{H.N.Z.}},
\bauthor{\bsnm{Delle~Monache}, \binits{M.L.}}:
\bctitle{Near collision and controllability analysis of nonlinear optimal velocity follow-the-leader dynamical model in traffic flow}.
In: \bbtitle{2023 62nd IEEE Conference on Decision and Control (CDC)},
pp. \bfpage{8057}--\blpage{8062}
(\byear{2023}).
\doiurl{10.1109/CDC49753.2023.10384086}
\end{bchapter}
\endbibitem

\bibitem[\protect\citeauthoryear{Zinat~Matin and Sowers}{2020a}]{9147363}
\begin{bchapter}
\bauthor{\bsnm{Zinat~Matin}, \binits{H.N.}},
\bauthor{\bsnm{Sowers}, \binits{R.B.}}:
\bctitle{Nonlinear optimal velocity car following dynamics (i): Approximation in presence of deterministic and stochastic perturbations}.
In: \bbtitle{2020 American Control Conference (ACC)},
pp. \bfpage{410}--\blpage{415}
(\byear{2020}).
\doiurl{10.23919/ACC45564.2020.9147363}
\end{bchapter}
\endbibitem

\bibitem[\protect\citeauthoryear{Zinat~Matin and Sowers}{2020b}]{9147244}
\begin{bchapter}
\bauthor{\bsnm{Zinat~Matin}, \binits{H.N.}},
\bauthor{\bsnm{Sowers}, \binits{R.B.}}:
\bctitle{Nonlinear optimal velocity car following dynamics (ii): Rate of convergence in the presence of fast perturbation}.
In: \bbtitle{2020 American Control Conference (ACC)},
pp. \bfpage{416}--\blpage{421}
(\byear{2020}).
\doiurl{10.23919/ACC45564.2020.9147244}
\end{bchapter}
\endbibitem

\bibitem[\protect\citeauthoryear{Matin et~al.}{2024}]{matin2024analytical}
\begin{botherref}
\oauthor{\bsnm{Matin}, \binits{H.N.Z.}},
\oauthor{\bsnm{Yeo}, \binits{Y.}},
\oauthor{\bsnm{Gong}, \binits{X.}},
\oauthor{\bsnm{Delle~Monache}, \binits{M.L.}}:
On the analytical properties of a nonlinear microscopic dynamical model for connected and automated vehicles.
IEEE Control Systems Letters
(2024)
\end{botherref}
\endbibitem

\bibitem[\protect\citeauthoryear{Delle~Monache et~al.}{2019}]{delle2019feedback}
\begin{botherref}
\oauthor{\bsnm{Delle~Monache}, \binits{M.L.}},
\oauthor{\bsnm{Liard}, \binits{T.}},
\oauthor{\bsnm{Rat}, \binits{A.}},
\oauthor{\bsnm{Stern}, \binits{R.}},
\oauthor{\bsnm{Bhadani}, \binits{R.}},
\oauthor{\bsnm{Seibold}, \binits{B.}},
\oauthor{\bsnm{Sprinkle}, \binits{J.}},
\oauthor{\bsnm{Work}, \binits{D.B.}},
\oauthor{\bsnm{Piccoli}, \binits{B.}}:
Feedback control algorithms for the dissipation of traffic waves with autonomous vehicles.
Computational Intelligence and Optimization Methods for Control Engineering,
275--299
(2019)
\end{botherref}
\endbibitem

\bibitem[\protect\citeauthoryear{Wu et~al.}{2017}]{wu2017multi}
\begin{bchapter}
\bauthor{\bsnm{Wu}, \binits{C.}},
\bauthor{\bsnm{Vinitsky}, \binits{E.}},
\bauthor{\bsnm{Kreidieh}, \binits{A.}},
\bauthor{\bsnm{Bayen}, \binits{A.}}:
\bctitle{Multi-lane reduction: A stochastic single-lane model for lane changing}.
In: \bbtitle{2017 IEEE 20th International Conference on Intelligent Transportation Systems (ITSC)},
pp. \bfpage{1}--\blpage{8}
(\byear{2017}).
\bcomment{IEEE}
\end{bchapter}
\endbibitem

\bibitem[\protect\citeauthoryear{Kesting et~al.}{2007}]{kesting2007general}
\begin{barticle}
\bauthor{\bsnm{Kesting}, \binits{A.}},
\bauthor{\bsnm{Treiber}, \binits{M.}},
\bauthor{\bsnm{Helbing}, \binits{D.}}:
\batitle{General lane-changing model mobil for car-following models}.
\bjtitle{Transportation Research Record}
\bvolume{1999}(\bissue{1}),
\bfpage{86}--\blpage{94}
(\byear{2007})
\end{barticle}
\endbibitem

\bibitem[\protect\citeauthoryear{Khelfa et~al.}{2023}]{khelfa2023predicting}
\begin{barticle}
\bauthor{\bsnm{Khelfa}, \binits{B.}},
\bauthor{\bsnm{Ba}, \binits{I.}},
\bauthor{\bsnm{Tordeux}, \binits{A.}}:
\batitle{Predicting highway lane-changing maneuvers: A benchmark analysis of machine and ensemble learning algorithms}.
\bjtitle{Physica A: Statistical Mechanics and its Applications}
\bvolume{612},
\bfpage{128471}
(\byear{2023})
\end{barticle}
\endbibitem

\bibitem[\protect\citeauthoryear{Lef{\`e}vre et~al.}{2014}]{lefevre2014survey}
\begin{barticle}
\bauthor{\bsnm{Lef{\`e}vre}, \binits{S.}},
\bauthor{\bsnm{Vasquez}, \binits{D.}},
\bauthor{\bsnm{Laugier}, \binits{C.}}:
\batitle{A survey on motion prediction and risk assessment for intelligent vehicles}.
\bjtitle{ROBOMECH journal}
\bvolume{1},
\bfpage{1}--\blpage{14}
(\byear{2014})
\end{barticle}
\endbibitem

\bibitem[\protect\citeauthoryear{Fehlberg}{1969}]{fehlberg1969low}
\begin{botherref}
\oauthor{\bsnm{Fehlberg}, \binits{E.}}:
Low-order classical runge-kutta formulas with stepsize control and their application to some heat transfer problems.
vol. 315.
National aeronautics and space administration
(1969)
\end{botherref}
\endbibitem

\bibitem[\protect\citeauthoryear{Hairer et~al.}{}]{hairer10solving}
\begin{botherref}
\oauthor{\bsnm{Hairer}, \binits{E.}},
\oauthor{\bsnm{N{\o}rsett}, \binits{S.P.}},
\oauthor{\bsnm{Wanner}, \binits{G.}}:
Solving ordinary differential equations i, nonstiff problems. 1993.
Springer-Verlag, Berlin, DOI
\textbf{10},
978--3
\end{botherref}
\endbibitem

\bibitem[\protect\citeauthoryear{Ince}{1956}]{ince1956integration}
\begin{botherref}
\oauthor{\bsnm{Ince}, \binits{E.L.}}:
Integration of ordinary differential equations.
(No Title)
(1956)
\end{botherref}
\endbibitem

\bibitem[\protect\citeauthoryear{Wilson et~al.}{2021}]{Argoverse2}
\begin{bchapter}
\bauthor{\bsnm{Wilson}, \binits{B.}},
\bauthor{\bsnm{Qi}, \binits{W.}},
\bauthor{\bsnm{Agarwal}, \binits{T.}},
\bauthor{\bsnm{Lambert}, \binits{J.}},
\bauthor{\bsnm{Singh}, \binits{J.}},
\bauthor{\bsnm{Khandelwal}, \binits{S.}},
\bauthor{\bsnm{Pan}, \binits{B.}},
\bauthor{\bsnm{Kumar}, \binits{R.}},
\bauthor{\bsnm{Hartnett}, \binits{A.}},
\bauthor{\bsnm{Pontes}, \binits{J.K.}},
\bauthor{\bsnm{Ramanan}, \binits{D.}},
\bauthor{\bsnm{Carr}, \binits{P.}},
\bauthor{\bsnm{Hays}, \binits{J.}}:
\bctitle{Argoverse 2: Next generation datasets for self-driving perception and forecasting}.
In: \bbtitle{Proceedings of the Neural Information Processing Systems Track on Datasets and Benchmarks (NeurIPS Datasets and Benchmarks 2021)}
(\byear{2021})
\end{bchapter}
\endbibitem

\bibitem[\protect\citeauthoryear{Lambert and Hays}{2021}]{TrustButVerify}
\begin{bchapter}
\bauthor{\bsnm{Lambert}, \binits{J.}},
\bauthor{\bsnm{Hays}, \binits{J.}}:
\bctitle{Trust, but verify: Cross-modality fusion for hd map change detection}.
In: \bbtitle{Proceedings of the Neural Information Processing Systems Track on Datasets and Benchmarks (NeurIPS Datasets and Benchmarks 2021)}
(\byear{2021})
\end{bchapter}
\endbibitem

\bibitem[\protect\citeauthoryear{Minderhoud and Bovy}{2001}]{minderhoud2001extended}
\begin{barticle}
\bauthor{\bsnm{Minderhoud}, \binits{M.M.}},
\bauthor{\bsnm{Bovy}, \binits{P.H.}}:
\batitle{Extended time-to-collision measures for road traffic safety assessment}.
\bjtitle{Accident Analysis \& Prevention}
\bvolume{33}(\bissue{1}),
\bfpage{89}--\blpage{97}
(\byear{2001})
\end{barticle}
\endbibitem

\bibitem[\protect\citeauthoryear{Van~der Horst}{1991}]{van1991time}
\begin{botherref}
\oauthor{\bsnm{Horst}, \binits{A.R.A.}}:
A time-based analysis of road user behaviour in normal and critical encounters.
(1991)
\end{botherref}
\endbibitem

\bibitem[\protect\citeauthoryear{Hogema and Janssen}{1996}]{hogema1996effects}
\begin{bchapter}
\bauthor{\bsnm{Hogema}, \binits{J.}},
\bauthor{\bsnm{Janssen}, \binits{W.}}:
\bctitle{Effects of intelligent cruise control on driving behavior: a simulator study}.
In: \bbtitle{Intelligent Transportation: Realizing the Future. Abstracts of the Third World Congress on Intelligent Transport Systems ITS America}
(\byear{1996})
\end{bchapter}
\endbibitem

\bibitem[\protect\citeauthoryear{Hirst and Graham}{2020}]{hirst2020format}
\begin{bchapter}
\bauthor{\bsnm{Hirst}, \binits{S.}},
\bauthor{\bsnm{Graham}, \binits{R.}}:
\bctitle{The format and presentation of collision warnings}.
In: \bbtitle{Ergonomics and Safety of Intelligent Driver Interfaces},
pp. \bfpage{203}--\blpage{219}.
\bpublisher{CRC Press},
\blocation{Boca Raton, FL}
(\byear{2020})
\end{bchapter}
\endbibitem

\end{thebibliography}
\end{document}